\documentclass[reqno]{amsart}

\pdfoutput=1

\usepackage{graphicx}
\usepackage{amsfonts}
\usepackage{amssymb}
\usepackage{enumitem}
\usepackage{Baskervaldx}
\usepackage{pdfpages}
\usepackage{tempora}

\usepackage[ukrainian,english]{babel}

\newtheorem*{theorem*}{Theorem}
 
\newtheorem*{lemma*}{Lemma}

\newtheorem*{corollary*}{Corollary}

\newtheorem{innercustomgeneric}{\customgenericname}
\providecommand{\customgenericname}{}
\newcommand{\newcustomtheorem}[2]{%
  \newenvironment{#1}[1]
  {%
   \renewcommand\customgenericname{#2}%
   \renewcommand\theinnercustomgeneric{##1}%
   \innercustomgeneric
  }
  {\endinnercustomgeneric}
}
\newcustomtheorem{customthm}{Theorem}
\newcustomtheorem{customlemma}{Lemma}

\newtheorem*{mainlemma*}{Main Lemma}

\newtheorem*{adianRabin*}{The Adian-Rabin Theorem}

\theoremstyle{definition}
\newcustomtheorem{customremark}{Remark}

\newtheorem*{definition*}{Definition}

\newcommand{\pres}[3]{\textnormal{#1} \langle #2 \mid #3 \rangle}

\begin{document}

\iffalse
%\thispagestyle{empty}
{

	\vfill
	\vspace{3cm}
	{\Large \textsc{Translated into English from the original Ukrainian}\par}
	{\Large \textsc{by} \\ \Large \textsc{C.-F. Nyberg-Brodda}\par}
	\vspace{2.0cm}
	{\Large 2023} \\

}
\clearpage 

\fi

%\thispagestyle{empty}

{
\centering

	{\scshape\LARGE Three Articles on One-relator Groups \\ by Wilhelm Magnus \par} 
	\vspace{1cm}
	{\Large \textsc{Translated into English}\par}
	{\Large \textsc{by} \\ \Large \textsc{C.-F. Nyberg-Brodda}\par}
	}
	\vspace{0.5cm}
\begin{center}
\today\\
\rule{\textwidth/2}{1.0pt}
\end{center}	
	
	\vspace{1cm}
%\begin{center}
\noindent{\large\textbf{Translator's Preface}}
%\end{center}

\vspace{0.5cm}

\noindent The present work consists of an English translation of three articles, originally written in German, by Wilhelm Magnus (1907--1990). The bibliographic details for the articles are as follows; the page numbers stated next to the naming of each article is the page number for the translation in the present work. 

\

\begin{center}
\makebox[\linewidth]{\hfill \textbf{[Mag1930]} \hfill \llap{(pp.\ 6--31)}}
\end{center}
W.\ Magnus, \textit{Über diskontinuierliche Gruppen mit einer definierenden Relation (Der Freiheitssatz)}, J.\ f\"ur die reine und angewandte Mathematik, \textbf{163}:3 (1930), pp.\ 141--165. \\

\begin{center}
\makebox[\linewidth]{\hfill \textbf{[Mag1931]} \hfill \llap{(pp.\ 32--49)}}
\end{center}
W.\ Magnus, \textit{Untersuchungen über einige unendliche diskontinuierliche Gruppen}, Math.\ Ann.\ \textbf{105} (1931), pp.\ 52--74. \\

\begin{center}
\makebox[\linewidth]{\hfill \textbf{[Mag1932]} \hfill \llap{(pp.\ 50--58)}}
\end{center}
W.\ Magnus, \textit{Das Identitätsproblem für Gruppen mit einer definierenden Relation}, Math.\ Ann.\ \textbf{106} (1932), pp.\ 295--307. \\

\

\noindent The subject matter of these articles is combinatorial group theory; but to summarize it in this way is almost an injustice to the profound significance these articles, particularly [Mag1930] and [Mag1932], would come to play in shaping that very area. Indeed, it would not be an exaggeration to say that these articles created the area that is today known as the theory of one-relator groups. This area has recently been extensively surveyed in its grand entirety, both from the point of view of modern progress and the history of mathematics, by myself and M.\ Linton. Thus, in this short preface, I could not, nor am I consequently under any obligation to, provide even a semblance of the full richness of this history, and gladly redirect the reader to that survey. 

Instead, I am left with ample space to discuss the actual content of the three articles above; and some of this content is as exciting as it is poorly known today. I will begin by presenting the two key results proved by Magnus on one-relator groups, which form the central pillars around which the three articles above are shaped. Let us set the stage for these theorems with some basic definitions. A \textit{one-relator group} $G = \pres{}{A}{r=1}$ is the quotient of the free group $F_A$ on some alphabet $A$ by the normal closure of some word $r \in F_A$. We may and will assume without loss of generality that the word $r$ is non-trivial and \textit{cyclically reduced}, i.e.\ that it is not of the form $xsx^{-1}$ for some $x \in F_A$ (were this the case, then of course $G$ would be identical to the one-relator group $\pres{}{A}{s=1}$). The \textit{word problem} (or \textit{identity problem}) in $G$ asks for an effective procedure to decide for a given word $w \in F_A$ whether or not $w =_G 1$, i.e.\ if it lies in the normal closure of $r$ or not. With this light dusting of a background, it is now easy to state the two key results of Magnus' articles above: 

\begin{theorem*}[W.\ Magnus]
Let $G = \pres{}{A}{r=1}$ be a one-relator group as above. Then:
\begin{enumerate}
\item (The \textit{Freiheitssatz}) If $A_0 \subset A$ is such that not every letter appearing in $r$ lies in $A_0$, then the subgroup of $G$ generated by $A_0$ is a free group, with basis $A_0$. 
\item The word problem in $G$ is decidable. 
\end{enumerate}
\end{theorem*}

The German word \textit{Freiheitssatz} (the ``Freeness Theorem'') has entered the dictionary also of English-speaking group theorists, and is always referred to as such. These theorems, when combined, form the backbone of one-relator group theory, and are often the first two results to be presented to a novice in the subject. Furthermore, the \textit{Freiheitssatz} is a key component in proving the decidability of the word problem (see below).  However, they are not the only two results to appear in the three articles on the subject, nor was there an immediate step from one to the other. I will therefore clarify the sequence of events by giving a brief summary of the three articles above. 

The article [Mag1930] is divided in two parts, with a short appendix. Magnus begins by mentioning that his doctoral supervisor M.\ Dehn (1878--1952) had given a series of seminars in Leipzig, in which a proof of the \textit{Freiheitssatz} had been sketched. Magnus was subsequently assigned the task of proving it as part of his doctoral thesis, and the remainder of the first part is thus easy to summarize: it is a proof of the \textit{Freiheitssatz}. However, in the course of doing this it also introduces the reader to all of combinatorial group theory, and all necessary setup is provided therein, including the word problem in free groups, normal forms for elements equal to $1$ in a finitely presented groups, and other elementary concepts. The proof of the \textit{Freiheitssatz} passes via what is today called the \textit{Magnus hierarchy}, which decomposes a one-relator group as the union of highly controlled amalgamated free products (see below) of one-relator groups in which the defining relation word has a shorter length. The high degree of control over these amalgamated free products then gives the result by induction, with the base case essentially being the case of a free group, for which the \textit{Freiheitssatz} is a special case of the Nielsen--Schreier Theorem (all subgroups of free groups are free). 

In the above overview of the idea of the proof of the \textit{Freiheitssatz}, there is one important anachronism: Magnus does not phrase his result in terms of amalgamated free products, and instead proves the required results ``manually'' in terms of three lemmas (Lemma~1--3). These lemmas may appear very clumsy and unmanageable to the modern reader, especially when compared to the elegance of amalgamated free products. The reason for Magnus' choice of formulation is simple: amalgamated free products had only appeared (in work by Schreier) three years earlier, and Magnus was unaware of their existence during the writing of his article. In a note added at the very end, he does however remark that his method has an easy reformulation in terms of Schreier's amalgamated free products, and that he will explore this further in future articles.

Before moving onto these future articles, I would like mention a few of the consequences of the \textit{Freiheitssatz} that already appear in [Mag1930]. The first concerns the (normal) \textit{root problem}, which Magnus introduces already in the same paragraph as the word problem. The root problem asks, for a given word $r \in F_A$, to find all words $s \in F_A$ such that $r =_G 1$ in the one-relator group $\pres{}{A}{s=1}$. Words $s$ of this form are called \textit{roots} of $r$. This problem is distinct from the three fundamental problems introduced by Dehn in 1911 -- the word problem, the conjugacy problem, and the isomorphism problem -- but it nevertheless is of great importance for one-relator groups. For any given pair of words $r, s$ one can of course determine if $s$ is a root of $r$ if one has a solution to the word problem in $\pres{}{A}{s=1}$. However, in general, determining \textit{all} roots of a given word is significantly more complicated than the word problem, and is even today only solved in a few cases. In the second part of [Mag1930], Magnus deals largely with the root problem. He proves two types of results in this line. The first is the following theorem, today usually known as the \textit{Conjugacy Theorem}:

\begin{theorem*}[Magnus, 1930]
Let $r_1, r_2 \in F_A$. If $r_1$ is a root of $r_2$, and $r_2$ is a root of $r_1$, then $r_1$ is conjugate to $r_2$ or $r_2^{\pm 1}$. 
\end{theorem*}

That is, if the normal closures of two words $r_1, r_2 \in F_A$ in a free group coincide, then $r_1$ is conjugate to $r_2$ or $r_2^{\pm 1}$ (and vice versa). This gives a form of uniqueness of one-relator presentations, but it is not in general sufficient to solve the isomorphism problem in one-relator groups, which remains an open problem. Nevertheless, the Conjugacy Theorem plays a key role in many of the cases that the isomorphism problem is known to be decidable.

The next application of the \textit{Freiheitssatz} in [Mag1930] is the determination of the roots of some particular words. In [\S7, Mag1930] Magnus determines all roots of the commutator $r = aba^{-1}b^{-1}$, and shows that the roots are either primitive words or a conjugate of $r$. This shows, for example, that $\pres{}{a,b}{aba^{-1}b^{-1} = 1}$ is the only one-relator presentation with a cyclically reduced relator for the free abelian group $\mathbb{Z}^2$, up to renaming generators. Magnus also determines the roots of $a^2 b^p, a^2 b^{2^k}$, and $a^p b^{p^k}$, where $p$ is a prime and $k \in \mathbb{Z}$. Already in the case of the word $ab^6a^{-1}b^{-6}$, which defines the Baumslag--Solitar group $\operatorname{BS}(6,6)$, Magnus notes that the problem of finding all its roots seems difficult, and is unable to do so. The problem of finding all roots of $ab^6a^{-1}b^{-6}$ would remain unsolved until 2000, when it was solved by McCool. 

Finally, [Mag1930] ends with a discussion of two-relator groups. Here, he notes that the situation is much more difficult than in the case of a single relation. In this context, he formulates and conjectures a natural form of the Conjugacy Theorem for two-relator groups, involving beyond conjugacy also other transformations akin to Nielsen transformations between the defining relations. This conjecture appears to still be open today, and is a very early version of a form of the Andrews--Curtis Conjecture formulated in 1965. He also considers some finite two-relator groups, notes their diverse nature, and finally asks whether or not is decidable if a two-relator group is ``essentially two-relator'', i.e.\ whether it does not admit a presentation with fewer than two defining relations. Both the problems of deciding when a two-relator group is finite as well as this latter problem posed by Magnus remain open today. 

Next, the 1931 article [Mag1931] is somewhat more difficult to summarize, as it is a rather eclectic article on many different subjects in combinatorial group theory. It is, however, explicitly inspired by the methods of [Mag1930], and the Magnus breakdown method used to prove the \textit{Freiheitssatz}. Magnus begins by formulating some essential results about amalgamated free products, now being aware of their power, and in [\S3, Mag1931] turns to apply these methods to determine the (outer) automorphism group of the \textit{figure-eight knot group} (Listing's knot group) This knot group is a one-relator group, and Dehn had in 1914 found two automorphisms $\overline{j}_0, \overline{j}_1$ of it, but was unable to prove that they generate the full outer automorphism group of the figure-eight knot group. Magnus first proves that they, together with two inner automorphisms, indeed generate the full automorphism group by using the Magnus breakdown procedure, and then determines a presentation for this group. The automorphism group turns out to be virtually free (although Magnus does not note this), and the outer automorphism group is isomorphic to the dihedral group $D_4$ with eight elements. 

In the next section [\S4, Mag1931], Magnus turns to solving the word problem in a class of one-relator groups: those with defining relation of the form $a^{\alpha_1} b^{\beta_1} a^{\alpha_2} b^{\beta_2} = 1$, where the exponents are arbitrary integers. The proof is not very far removed from his subsequent 1932 proof of the decidability of the word problem in \textit{all} one-relator groups (see below), but the inductive hypothesis is weaker, and Magnus relies on the fact that the Magnus breakdown procedure terminates quickly for the above one-relator groups. In the last section [\S5, Mag1931], Magnus considers a very different problem: determining the subgroups of the modular group $\operatorname{PSL}_2(\mathbb{Z})$, i.e.\ the group given by two generators $a, b$ subject to the defining relations $a^2 = b^3 = 1$. By using the Magnus breakdown procedure, he proves that the commutator subgroup of this group is free on two generators, and that, since the quotient of $\operatorname{PSL}_2(\mathbb{Z})$ by the commutator subgroup is cyclic of order $6$, this can be used to determine all subgroups of $\operatorname{PSL}_2(\mathbb{Z})$.

This leaves only the final of the three articles. In [Mag1932], Magnus realized the full power of the Magnus hierarchy for solving the word problem, and used it to prove that \textit{every} one-relator group has decidable word problem. The missing piece was a stronger inductive hypothesis: he needed to show that every one-relator group has decidable \textit{generalized} word problem. This problem, which today is usually called the subgroup membership problem with respect to Magnus subgroups, goes as follows. Let $G = \pres{}{A}{r=1}$, and take as input a word $w \in F_A$ and some subset $A_0 \subseteq A$. Then one is asked to decide whether or not $w$ is an element of the subgroup of $G$ generated by $A_0$. If $r$ involves all letters from $A$ (or their inverses), then any proper subset $A_0 \subset A$ generates a free subgroup by the \textit{Freiheitssatz}, and these are the amalgamated subgroups in the Magnus hierarchy. By standard results on normal forms in amalgamated products of groups, the word problem in an amalgam $G_1 \ast_H G_2$ can be reduced to the word problems in the factors $G_1, G_2$ together with the membership problem in the amalgamated subgroup $H \leq G_1, G_2$ in each factor. With a little extra work, in the case of one-relator amalgams one can also reduce the generalized word problem in this way, which thus again inductively gives the solution to the generalized word problem in all one-relator groups. Unlike the previous two articles, [Mag1932] is devoted to proving only a single result, and ends after having accomplished this; the proof is, as usual, divided into two separate cases, depending on whether some generator has exponent sum zero in $r$ or not. The article ends abruptly after completing this proof, and demonstrates how quickly one-relator group theory had matured in the short time since 1930. Finally, we remark that in the introduction of [Mag1932], Magnus claims that any non-cyclic one-relator group contains a non-abelian free subgroup, except when the group is given by two generators $a, b$ subject to the relation $aba^nb^{-1}=1$, i.e.\ a solvable Baumslag--Solitar group. No proof of this is given, and a proof would first appear in print in a short 1969 article by D.\ I.\ Moldavanskii.

\clearpage

This concludes the overview of the three remarkable articles by Magnus. From a technical point of view, I have attempted to be faithful to the original language as far as possible, choosing e.g.\ to translate \textit{Identitätsproblem} with the historical term ``identity problem'' rather than the modern term ``word problem''. Furthermore, Magnus did not yet have a consistent notation for group presentations, and often changes this throughout a given article -- I have kept all these different pieces of notation as they were, since the notation is perfectly clear in spite of its inconsistency. I have endeavoured to maintain the original design of the articles as closely as possible, including the fonts; reading the translations alongside the originals should not pose any difficulty. I have also not corrected any mathematical issues beyond obvious and small typesetting glitches like the one mentioned above; this is primarily because I only spotted a single error. This is at the very last page of [Mag1931], where Magnus computes two Möbius transformations corresponding to generators of rank $2$ free commutator subgroup in $\operatorname{PSL}_2(\mathbb{Z})$. He claims to compute the Möbius transformation for the generator $\beta_0 = aba^{-1}b^{-1}$ as $\frac{-z+1}{z-2}$, but this transformation instead corresponds to $bab^{-1}a^{-1}$, i.e.\ $\beta_0^{-1}$. The correct transformation corresponding to $\beta_0$ is $\frac{2z+1}{z+1}$. 

At this stage, I would naturally recommend any interested reader to consult the aforementioned recent survey on one-relator group theory, which covers a great deal more detail on the historical and modern context. However, I would also recommend the reader to take the time to peer through Magnus' articles, which cover a surprisingly wide range of topics with a skilled level of precision and care. The articles are all very readable to modern-day group theorists, and although some of the notation appears either clumsy or somewhat inefficient, it remains perfectly comprehensible. They are also well-written: for example, in the first article [Mag1930], I only spotted a single typo, where on p.\ 156, line 6, one reads $K_i^{(ti)}$ rather than the correct $K_i^{(\beta_i)}$, which may equally well have been caused by an overly enthusiastic typesetter. The articles also contain hints of a great deal of combinatorial group theory decades ahead of time, including the Andrews--Curtis Conjecture, and stand as a testament to the creativity of Magnus. Above all, the articles are enjoyable to read, and it has been a joy to typeset and translate them in the interest of making them more accessible to a modern audience.

\

\

\begin{flushright}
\noindent\textbf{C. F. Nyberg-Brodda}\\
\noindent{Korea Institute for Advanced Study}\\
(Seoul, Republic of Korea) \\
\today
\end{flushright}



\section*{Supplementary material}

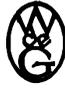

# On discontinuous groups with one defining relation. (The *Freiheitssatz*)

By *Wilhelm Magnus* in Frankfurt am Main.

## Contents

## Introduction

In 1922, Herr Dehn held a seminar (unpublished) in Leipzig on topology and general group theory. In its group-theoretic part, an approach to a general calculus for discontinuous groups is developed. The present work was created at the suggestion of Herr Dehn, and gives – with modified proofs – a large part of his results, in particular

the "Freiheitssatz" necessary for many constructions; in addition, this present work contains theorems which Herr Dehn developed in his Frankfurt lectures, as well as a number of general and special theorems as applications of the theory.

The starting point is the following definition of the *identity problem* of group theory: consider some letters $a, a^{-1}; b, b^{-1}; \ldots$. A sequence of such letters is called a *word*. The word which contains no letters is denoted by 1. Two words $W_1$ and $W_2$ are *identical*, when it is possible to transform $W_1$ to $W_2$ by insertions and deletions of the special words $aa^{-1}, a^{-1}a; bb^{-1}, b^{-1}b; \ldots$. We denote this $W_1 \equiv W_2$. In particular, $aa^{-1} \equiv 1, a^{-1}a \equiv 1, \ldots$.

The question of when it is possible to transform a given word $W$ into the empty word using insertions and deletions of special, fixed words $R_1, R_2, \ldots$ and their "reciprocals" $R_1^{-1}, R_2^{-1}, \ldots$ (for which $R_1 R_1^{-1} \equiv 1$ follows from the identical transformations) is called the identity problem for the group with the generators $a, b, \ldots$ and the defining relations $R_1 = 1, R_2 = 1, \ldots$: we say that $W = 1$ is a consequence of the relations $R_1 = 1, R_2 = 1, \ldots$.

In the present work we restrict our attention to one-relator groups. If from $R = 1$ it follows that $W = 1$, then we say that $R$ is a *root* of $W$. The identity problem asks whether or not $R$ is a root for a given word $W$. One can try to deal with the question by reversing it, and solving a much more comprehensive problem: for a given word $W$, find all of its roots. This problem is the simplest "*root problem*"; the general root problem asks to find all systems of $n$ relations $R_1 = \cdots = R_n = 1$ from which $W = 1$ follows. In fact, the solutions to the root problem include those to the identity problem; because if one wants to know, in a group with generators $a, b, \ldots$ and relation $R = 1$, whether or not the word $W$ is 1, one has – using the solution of the root problem for $W$ – only to check whether $R$ belongs to the roots of $W$ or not.

The first question that one treats when considering the root problem concerns the roots of the "simplest" word: the generators; that is, the question: what must $R$ look like in order to have $a = 1$ in the group with generators $a, b, \ldots$ and $R = 1$? The answer is the most important consequence of the "*Freiheitssatz*" (which we shall give a detailed formulation of below) and reads: we must have $R \equiv T_1 a T_1^{-1}$ or $R \equiv T_2 a^{-1} T_2^{-1}$, where $T_1$ and $T_2$ are arbitrary words. That is: $a$ only has trivial roots: the conjugates of $a^{\pm 1}$.

# First Part: The *Freiheitssatz*

## §1. Formulation

Consider generators $a_1, a_2, \ldots, a_n$ and $x$, and one relation holding between them: $R(a_1, \ldots, a_n; x) = 1$. If $R = 1$ yields some relation $W(a_1, \ldots, a_n) = 1$ holding only between the $a$, and such that this relation is not identically true (so the $x$ can be eliminated, in other words), then $R \equiv TS(a_1, \ldots, a_n)T^{-1}$, where $T$ (possibly) contains $x$, and $W = 1$ follows from $S = 1$.

That is: if $R$ contains $x$ when "cyclically written"[1], i.e.: if one cannot make the cyclically written word $R$ such that it no longer contains $x$ by identical transformations (how one decides whether or not this is possible is given in the beginning of §2), then between the $a_1, \ldots a_n$ there is no non-identical relation: $a_1, \ldots a_n$ generate in the group with the relation $R = 1$ a *free group*.

---

[1] In a cyclically written word one considers the first letter as adjacent to the last, and can define identical transformations accordingly.

## §2. Simple tools

1. We will first of all give a solution to the *identity problem for free groups*. This is done via the following theorem: if a word $W$ on the generators $a, b, c, \ldots$ is equal to 1, then it can be transformed into the empty word by only applying successive *deletions* of $aa^{-1}, a^{-1}a; \ldots$, independently of the order in which the deletions are performed. One consequence of this theorem is: if a word $W$ (written normally or cyclically) contains a certain generator $x$ even when it cannot be further "simplified" by "deletions", i.e. if one cannot delete any words $aa^{-1}, a^{-1}a, \ldots$ from it, then this generator $x$ cannot be removed from $W$ by identical transformations at all (cf. §1).

2. *Normal form for words which in a given group equal* 1. The proof of the previous theorem is carried out in the same way as that of the following theorem, which provides a kind of normal form for the word*s* $W$ over the generators $a, b, \ldots$ which as a consequence of the relations

$$(1) \qquad\qquad R_1 = 1, \quad R_2 = 1, \quad \ldots \quad R_n = 1$$

is equal to 1. In this case we then have:

$$(N) \qquad\qquad W \equiv \prod_{i=1}^{h} T_i R_{\nu_i}^{e_i} T_i^{-1}$$

where the product sequence is symbolic: the factors in the product cannot be permuted. The conjugating words $T_i$ are arbitrary words; $e_i$ has one of the values $\pm 1$; and $\nu_i$ takes one of the values $1, 2, \ldots n$. (That, conversely, $W$ is 1 in the group: $a, b, \ldots; R_1 = 1, \ldots, R_n = 1$; is trivial).

*Proof.* The theorem holds for all words $\overline{W}$ which can be transformed into the empty word using identical transformations and a *single* insertion or deletion of one of the $R_\nu^{\pm 1}$. Assume that it is proved for all words $\overline{W}$ which can be transformed into the empty word using identical transformations and at most $m$ (where $m \geq 1$) insertions or deletions of one of the $R_\nu^{\pm 1}$. Then it holds also for $W$, where $W$ can be transformed into the empty word using identical transformations and at most $m+1$ such insertions or deletions. This is proved as follows: $W$ can be transformed into a word $\overline{W}$ using identical transformations and a single insertion resp. deletion of one of the $R_\nu^{\pm 1}$. In the sequel we fix for $\overline{W}$ the identity:

$$(2) \qquad\qquad \overline{W} \equiv \prod_{k=1}^{\overline{h}} T_k R_{\nu_k}^{e_k} T_k^{-1}$$

with $e_k = \pm 1, \nu_k = 1, 2, \ldots, n$. We can now transform $\overline{W}$, by inserting or deleting a $R_\nu^{\pm 1}$ – and after the previous identical transformations – into a word which is identical to $W$. The word $R_\nu^{\pm 1}$ may be assumed to be $R_1^{+1}$. If $W$ is identical to $\overline{W}$ after an insertion of $R_1$, then we can write $W$ in two parts, such that $W \equiv W_1 W_2$ and $\overline{W} \equiv W_1 R_1 W_2 \equiv W_1 W_2 W_2^{-1} R_1 W_2$ which, by (2) gives:

$$\overline{W} \equiv \prod_{k=1}^{\overline{h}} T_k R_{\nu_k}^{e_k} T_k^{-1} \equiv W \cdot W_2^{-1} R_1 W_2 \quad \text{and} \quad W \equiv \left( \prod_{k=1}^{\overline{h}} T_k R_{\nu_k}^{e_k} T_k^{-1} \right) W_2^{-1} R_1^{-1} W_2.$$

If instead $W$ is obtained from $\overline{W}$ be a deletion of $R_1$, then we have $W \equiv \overline{W}_1 R_1 \overline{W}_2$; $\overline{W} \equiv \overline{W}_1 \overline{W}_2$, and so $W \equiv (\overline{W}_1 \overline{W}_2) \overline{W}_2^{-1} R_1 \overline{W}_2 \equiv \left( \prod_{k=1}^{\overline{h}} T_k R_{\nu_k}^{e_k} T_k^{-1} \right) \overline{W}_2^{-1} R_1 \overline{W}_2$. In both cases we have for $W$ found an expression of the form (N). We give an example: $a^2 b a^{-2} b^{-1} \equiv W$ is equal to $1$ when $ab = ba$ holds, i.e. when $aba^{-1}b^{-1} \equiv R$ is $1$. In this case we have

$$a^2 b a^{-2} b^{-1} \equiv [a(aba^{-1}b^{-1})a^{-1}](aba^{-1}b^{-1}).$$

<div style="text-align:right">□</div>

3. We will in the sequel use *invariants of identical transformations.* One such is the "*exponent sum*", which a given generator $a, b, c, \dots$ possesses in a word $W(a, b, c, \dots)$ written over them; this is defined as follows:

We can write $W \equiv a^{\alpha_1} b^{\beta_1} c^{\gamma_1} \cdots a^{\alpha_2} b^{\beta_2} c^{\gamma_2} \cdots a^{\alpha_n} \cdots$, where at least one of the $\alpha, \beta, \dots$ is non-zero. Then $\alpha = \alpha_1 + \alpha_2 + \cdots + \alpha_n$ is the exponent sum of $a$ in $W$. The invariance follows from, for example, that upon an insertion or deletion of $aa^{-1}, \dots$, the exponent sum of $a, \dots$ does not change; or from the observation that when $W \equiv W'$, i.e. $W W'^{-1} \equiv 1$, all generators must have zero exponent sum, on the basis of the following (generalizable) argument: $W W'^{-1}$ can certainly be transformed into $1$, when one allows all letters to commute, i.e. by setting $aba^{-1}b^{-1} = 1$ etc. But then $W W'^{-1} = a^{\overline{\alpha}} b^{\overline{\beta}} c^{\overline{\gamma}} \cdots$ where $\overline{\alpha}, \overline{\beta}, \dots$ resp. are the exponent sums of $a, b, \dots$ in $W W'^{-1}$, respectively. From $W W'^{-1} \equiv 1$ follows $\overline{\alpha} = \overline{\beta} = \cdots = 0$. Due to this second point of view, arguments that are based on the invariance of the exponent sum under identical transformations will be referred to as arguments "by abelianization".

4. *Transformations into new generators. Independence of the same.* Suppose $W_1$ and $W_2$ are words over $a_1, \dots, a_n$ and $b$, that $W_1(a_1, \dots, a_n, b) \equiv W_2(a_1, \dots, a_n, b)$ and that $b$ has exponent sum zero in $W_1$ (and therefore also in $W_2$). Then we can write $W_1$ and $W_2$ as words over conjugates of the $a_\nu$ $(\nu = 1, \dots n)$ by powers of $b$, i.e. as words over $b^k a_\nu b^{-k}$ $(k = 0, \pm 1, \pm 2, \dots)$. These will be denoted by $a_\nu^{(k)}$. Thus we have:

$$W_1(a_1, \dots, a_n, b) = F_1(\cdots, a_\nu^{(k)}, \dots);$$
$$W_2(a_1, \dots, a_n, b) = F_2(\cdots, a_\nu^{(k)}, \dots);$$

*The lemma then states that $F_1 \equiv F_2$ are identical as words over the $a_\nu^{(k)}$*; by abelianizing in the $a_\nu^{(k)}$ we obtain sharper statements than by abelianizing in $a_1, \dots, a_n, b$.

For a *Proof* it suffices to show that a word $W(a_1, \dots, a_n, b)$ which in the $a_1, \dots, a_n, b$ is identically $1$, is also identically $1$ in the $a_\nu^{(k)}$. The proof of this is here only given in the case that in $W$ only a single generator $a$, other than $b$, occurs; the general case follows by complete induction and applying the reasoning used in this simple case.

Suppose then that $W(a, b) \equiv 1$, and:

$$W \equiv \prod_{i=1}^{h} b^{k_i} a^{e_i} b^{-k_i} \quad \text{with } e_i = \pm 1,$$

where $k_i = 0, \pm 1, \pm 2, \dots$. Since $W \equiv 1$, two $b$-letters must cancel at least once between two $a$-letters with opposite signs; that is, for at least *one* value $i_0$ of $i$ we have $k_{i_0} = k_{i_0+1}$ and $e_{i_0} = -e_{i_0+1}$. Hence we can immediately cancel $b^{k_{i_0}} a^{e_{i_0}} b^{-k_{i_0}}$ with $b^{k_{i_0+1}} a^{e_{i_0+1}} b^{-k_{i_0+1}}$, that is: *all* deletions can be performed in the $a^{(k)}$ (where $a^{(k)} = b^k a b^{-k}$), as $b^{k_{i_0}} a^{e_{i_0}} b^{-k_{i_0}} = [a^{(k_{i_0})}]^{e_{i_0}}$ and $b^{k_{i_0+1}} a^{e_{i_0+1}} b^{-k_{i_0+1}} = [a^{(k_{i_0})}]^{-e_{i_0}}$.

## §3. Reduction of the proof of the *Freiheitssatz* to two lemmas

1. *Case distinction.* To prove the *Freiheitssatz* there are two cases to consider:

   (I) From $R(a; x) = 1$ follows a non-trivial relation for $a$, i.e. $a^n = 1$, $n \neq 0$.

   (II) From $R(a_1, \ldots, a_n; x) = 1$ follows some non-trivial relation for $a_1, \ldots, a_n$: say $W(a_1, \ldots, a_n) = 1$, and in $R$ at least two distinct $a$-letters really appear.

**Remark**. It might seem as if this case distinction is not a complete disjunction, since it could *a priori* be the case that from $R(a_1, x) = 1$ no relation for $a_1$ follows, but that there is some relation holding between the $a_1, a_2, \ldots, a_n$. With the method of the following paragraph, however, one can prove:

If some (non-trivial) relation between the $a_1, \ldots, a_n$ follows from $R(a_1, x) = 1$, then a relation for $a_1$ alone also follows.

2. *Case I*: if from $R(a, x) = 1$ it follows that $a^n = 1$, then from §2 it follows that

$$(1) \qquad a^n \equiv \prod_{i=1}^{h} T_i(a, x) R^{e_i} T_i^{-1}$$

where $e_i = \pm 1$. From this it follows that $x$ has exponent sum zero in $R$, whereas $a$ has a non-zero exponent sum in $R$.[2]

Define for all integers $k$: $x^k a x^{-k} = b_k$, so in particular $a = b_0$, and write both sides of (1) over the $b_k$. The left-hand side will be $b_0^n$. The right-hand side will be transformed as follows: it transforms $R(a, x)$ into a word $P(\ldots, b_l, \ldots)$ over the $b_k$. The index $l$ takes, in general, many different values inside $P$. The exponent sum of $x$ in $T_i$ is denoted $t_i$.

Then: $T_i = \overline{T}_i(\ldots, b_m, \ldots) \cdot x^{t_i}$ and hence

$$T_i R T_i^{-1} = \overline{T}_i(\ldots, b_m, \ldots) x^{t_i} P(\ldots, b_l, \ldots) x^{-t_i} \overline{T}^{-1}.$$

By conjugating the $b_l$ inside $P$ by $x^{t_i}$ and using the definition of $b_l$, we thus obtain:

$$(2) \qquad b_0^n \equiv \prod_{i=1}^{h} \overline{T}_i P^{e_i}(\cdots, b_{l+t_i}, \ldots) \overline{T}_i^{-1}$$

and this is by §2 also an identity in the $b_k$. This gives: the relation $b_0^n = 1$ follows from a finite number of the relations

$$(3) \qquad P(\ldots, b_{l+t}, \ldots) = 1; \quad t = 0, \pm 1, \pm 2, \ldots$$

which we will write as:

$$(3) \qquad \begin{cases} P(\ldots, b_l, \ldots) = 1, \\ P(\ldots, b_{l+1}, \ldots) = 1, \quad P(\ldots, b_{l-1}, \ldots) = 1 \\ P(\ldots, b_{l+2}, \ldots) = 1, \quad P(\ldots, b_{l-2}, \ldots) = 1 \\ \cdots\cdots\cdots\cdots\cdots\cdots\cdots\cdots\cdots\cdots\cdots\cdots \end{cases}$$

It appears that we have made the problem much more complicated.

---

[2]This is because by abelianizing it follows that when we set $\sum_{i=1}^{h} e_i = \varrho$, and $\alpha$ resp. $\xi$ are the exponent sums of $a$ resp. $x$ in $R$, then:

$$\begin{aligned} n &= \varrho \cdot \alpha, \\ 0 &= \varrho \cdot \xi, \end{aligned} \quad \text{and hence: } \varrho \neq 0, \ \alpha \neq 0, \ \xi = 0.$$

Indeed, we have that the relation $b_0^n = 1$ should follow from an infinite number of relations between infinitely many generators. But first, these generators appear in a very special manner inside the relations, and second, these relations are isomorphic: they are obtained from one another when one replaces $b_l$ by $b_{l+t}$, where $t$ is an integer (fixed for the relation), and third the number of letters in $P$ is fewer than the number of characters in $R$, namely fewer by exactly the number of $x$-letters in $R$.

The reformulation of the problem carried out here is the decisive step in proving the *Freiheitssatz*. Before we proceed, we must carry out an analogous reformulation also in Case II.

3. *Case II.* When from $R(a_1, \ldots, a_n, x) = 1$ we have $W(a_1, \ldots, a_n) = 1$, and $a$ occurs twice in $R$, then one cannot anymore assume that $x$ has exponent sum zero in $R$, because it could be the case that the exponent sums of the $a_\nu$ in $W$ are all zero. *However, one can always assume that $a_1$ in $R$ has exponent sum zero in $W$.* Indeed: either one of the $a_\nu$ that occur in $R$ has exponent sum zero. In this case, we call this $a_1$ (and proceed with $a_1$ exactly as we will do with $b_1$ (see below)). Or else, $a_1$ resp. $a_2$ has non-zero exponent sum $s_1$ resp. $s_2$ in $R$. Then we set $a_1 = b_1^{+s_2}$ and $a_2 = b_1^{-s_1} b_2$. When doing this, the words $R$ and $W$ become rewritten to $\overline{R}(b_1, b_2, a_3, \ldots, a_n, x)$ resp. $\overline{W}(b_1, b_2, a_3, \ldots, a_n)$, *in which $b_1$ occurs in $\overline{R}$ with exponent sum zero.* The identity

$$(4) \qquad W \equiv \prod_{i=1}^{h} T_i R^{e_i} T_i^{-1} \qquad\qquad (e_i = \pm 1)$$

over $a_1, \ldots, a_n, x$, which says that $R$ is a root of $W$, becomes the identity

$$(4') \qquad \overline{W} \equiv \prod_{i=1}^{h} \overline{T}_i \overline{R}^{e_i} \overline{T}_i^{-1},$$

where the $\overline{T}_i$ is obtained from the $T_i$ by replacing the $a_1$ and $a_2$ for $b_1$ and $b_2$. The identity (4') says that $\overline{R}$ is a root of $\overline{W}$. It is now important to ensure that in the transition from the $a$ to the $b$ "nothing is lost"[3], that is, one must show: when one can prove that on the basis of (4') the word $\overline{R}$, written cyclically, no longer contains $x$, then $R$, cyclically written, no longer contains $x$. This is accomplished by the following remark: if one is given a word $V(a_1, \ldots, a_n, x)$ in which no absorptions can take place between the $a_1, \ldots, a_n, x$ (that is, one cannot delete any $a_1 a_1^{-1}, \ldots$, or $x^{-1}x$), and if one makes the substitution:

$$(5) \qquad a_1 = b_1^{+s_2}, \quad a_2 = b_1^{-s_1} b_2, \quad s_2 \neq 0,$$

then $V$ becomes a word $V(b_1, b_2, a_3, \ldots, a_n, x)$ in which one may find absorptions between the $b_1$-letters, but in which (even after carrying out the absorptions of the $b_1$-letters) there are no possible absorptions between the $b_2, a_3, \ldots, a_n, x$-letters. (In particular, if $V$ is not 1, then neither is $\overline{V}$.) The proof is trivial.

The letter $b_1$ has exponent sum zero in $\overline{R}$, and hence – by (4') – also in $\overline{W}$. We conjugate both sides of (4') by $b - 1$, and set for all integers $k$:

$$c_{\nu,k} = \begin{cases} b_1^k b_2 b_1^{-k} & \text{for } \nu = 2, \\ b_1^k a_\nu b_1^{-k} & \text{for } \nu = 3, \ldots, n, \end{cases}$$

---

[3] Indeed, in general one cannot recover $b_1$ and $b_2$ expressed as words over $a_1$ and $a_2$.

and further: $x_k = b_1^k x b_1^{-k}$. Then $\overline{W}$ is transformed into a word $F(\ldots, c_{\nu,k}, \ldots)$, and $\overline{R}$ into a word $P(\ldots, c_{\nu,l}, \ldots, ; \ldots, x_m, \ldots)$, where $l, m$ take on all possible values inside $P$. It is important that $P$ does not only have fewer letters (*not* original generators) than $\overline{R}$, but also fewer than $R$. In fact, though the substitution (5) may have made the number of $b_1$-letters in $\overline{R}$ greater than the number of $a$-letters in $R$, we have that the number of $b_2, a_3, \ldots, a_n, x$-letters in $\overline{R}$ is the same as the corresponding number of $a_2, a_3, \ldots, a_n, x$-letters in $R$, and the $b_1$-letters disappear upon the conjugations by $b_1$. In the case that from the very beginning $a_1$ has exponent sum zero in $R$, we conjugate using this, set $a_1^k a_\nu a_1^{-k} = c_{\nu,k}$ for $\nu = 2, 3, \ldots, n$, and for all integers $k$ set $a_1^k x a_1^{-k} = x_k$, performing the transformation in (4) with $a_1$. We then find, just as in Case I, that $F(\ldots, c_{\nu,k}, \ldots) = 1$ follows from the relation system

(6)               $P(\ldots, c_{\nu,l+t}, \ldots; \ldots, x_{m+t}, \ldots) = 1 \quad (t = 0, \pm 1, \ldots).$

(Here $\nu, l, m$ are variable inside $P$.) This means that $F = 1$ follows from a finite number of relations in the infinite system (6).

4. *Plan of the proof.* The proof of the *Freiheitssatz* will be executed as follows: three lemmas will allow us to prove: if from a system (3) (in Case I) it follows that $b_0^n = 1$, then one generator that actually occurs in one of the relations (3) can be "eliminated" (cf. §1), and analogously in Case II: because (6) implies a relation (namely $F = 1$) which no longer contains $x$, it must be possible to eliminate $x$ from one of the relations of (6). Then we can use complete induction: the *Freiheitssatz* is trivially true in the case that $R$ contains only a single generator – especially when $R$ is just a single letter – since it cannot contain $x$.[4] Assume the theorem is proved for all words containing fewer letters than $R$. Furthermore, $R$, written cyclically, allows for no absorptions. Then it follows that none of the $P$ in (3) res. (6) allows for any absorptions.

On the other hand, the *Freiheitssatz* applies to the $P$ in (3) and (6), which yields: if a generator can be eliminated from $P$, then this word, when written cyclically, no longer contains it. (Indeed, each of the $P$ contains fewer letters than $R$). And then we are done.

5. *Formulation of three lemmas; reduction to the first two.* It is now time to formulate, apply, and prove the aforementioned lemmas (see Nr. 4). The proof will be carried out in the subsequent paragraphs. When reading the formulation of the lemmas, no heed should be paid to the previously defined terms.

**Lemma 1**. *Given four systems of generators:*

$$a_1, \ldots, a_m; \quad b_1, \ldots, b_n; \quad x_1, \ldots, x_r; \quad y_1, \ldots, y_s$$

*and two systems of (finitely many) relations between these:*

(A)               $\{P_\mu(a_1, \ldots; \quad b_1, \ldots; \quad x_1, \ldots) = 1\} \quad \mu = 1, 2, \ldots$

(B)               $\{Q_\nu(b_1, \ldots; \quad x_1, \ldots; \quad y_1, \ldots) = 1\} \quad \nu = 1, 2, \ldots$

*Then we have: if neither* (A) *alone nor* (B) *alone results in a relation between the* $a$ *and* $b$ *alone, and neither* (A) *alone nor* (B) *alone results in a relation between the* $b$ *and* $x$ *alone – then* (A) *and* (B) *taken together do not result in a relation between the* $a$ *and* $b$ *alone.*

----

[4]*Proof.* If $\prod_{i=1}^h T_i(x^m)^{e_i} T_i^{-1} \equiv W(a_1, \ldots, a_n)$, where $e_i = \pm 1$ and $m \neq 0$, then one use the observation that an identity remains true when one everywhere replaces a generator by another (which may also occur elsewhere) or by 1, and one replaces in the above identity $x$ by 1. That yields $W \equiv 1$.

**Lemma 2.** *Given three systems of generators:*

$$a_1, \ldots, a_m; \quad x_1, \ldots, x_r; \quad y_1, \ldots, y_s$$

*and two systems of (finitely many) relations:*

(Ā)     $\{\overline{P}_\mu(a_1, \ldots; \quad x_1, \ldots) = 1\}$   $\mu = 1, 2, \ldots$

(B̄)     $\{\overline{Q}_\nu(a_1, \ldots; \quad y_1, \ldots) = 1\}$   $\nu = 1, 2, \ldots$

*If neither* (Ā) *alone nor* (B̄) *yields any relation between the* $a_1, \ldots, a_m$, *then no such relation follows from* (Ā) *and* (B̄) *combined.*

**Use of the lemmas.** Remark: Lemma 2 is independent of Lemma 1; Lemma 1 will serve to reduce the missing piece of the proof of the *Freiheitssatz*, i.e. the proof of the statements made above (see Nr. 4) about the relation systems (3) resp. (6), to the proof of the following lemma:

**Lemma 3.** *Given two systems of generators:* $d_1, d_2, \ldots, d_t$ *and* $p_0, p_1, p_2, \ldots, p_H$, *and a system of* $H - K + 1$ $(K > 0)$ *relations between these with the following description:*

- *No relation, written cyclically, allows for any absorption.*
- *If* $Q_0 = 1, Q_1 = 1, \ldots, Q_{H-K} = 1$ *are the* $H - K + 1$ *relations, then for*

$$0 \leq \nu \leq H - K \quad \text{the word} \quad Q_\nu(d_1, \ldots; p_\nu, p_{\nu+1}, \ldots, p_{\nu+K})$$

  *contains, other than the* $d$-*letters (and possibly no* $d$-*letters), at most the generators* $p_\nu$ *to* $p_{\nu+K}$, *and indeed* $p_\nu$ *and* $p_{\nu+K}$ *actually occur in* $Q_\nu$.

*The relation system*

$$(\mathsf{q}) \quad \begin{cases} Q_0(d_1, \ldots; p_0, \ldots, p_K) = 1 \\ \quad Q_1(d_1, \ldots; p_1, \ldots, p_{K+1}) = 1 \\ \quad \ldots\ldots\ldots\ldots\ldots\ldots\ldots\ldots\ldots \\ \quad\quad Q_{H-K}(d_1, \ldots; p_{H-K}, \ldots, p_H) = 1 \end{cases}$$

*thus possesses very similar properties to a finite subsystem of (3) (in the case that in* (q) *no* $d$ *appears) or of (6).*

*Then we claim: if from the relation system* (q) *a relation for* $d_1, \ldots$ *and* $p_0, p_1, \ldots, p_S$ *follows, where* $0 \leq S < K$, *then necessarily in one of the relations of* (q), *say* $Q_\nu = 1$, *one of the generators* $p_\nu$ *or* $p_{\nu+K}$ *can be eliminated[5].*

The proof of this, assuming Lemma 1, is given as follows:

*Proof.* Lemma 3 is trivial when $H - K = 0$, i.e. when (q) only consists of a single relation. Assume it is proved for all relation systems (q) which consist of fewer than $H - K + 1$ relations. Then it also holds for the given relation system (q). Proof: there are two cases to consider: first: we have $H - K \leq S$. We use Lemma 1, and we identify: the generators $b_1, \ldots$ from Lemma 1 with $d_1, \ldots, d_t$ and $p_1, \ldots, p_S$; the generators $a_1, \ldots$ with $p_0$; the generators $x_1, \ldots$ with $p_{S+1}, \ldots, p_K$; the generators $y_1, \ldots$ with $p_{K+1}, \ldots, p_{H-K}$; the relation system (B) with $Q_1 = 1, Q_2 = 1, \ldots, Q_{H-K} = 1$; and finally the relation system (A) with $Q_0 = 1$. We show this assignment using

---

[5]A relation system in which the assumptions that (q) satisfies all hold, will in the sequel be called "*a system of type* (q)".

TABLE 1

|      | $a$         | $b$                                      | $x$              | $y$           |          |
|------|-------------|------------------------------------------|------------------|---------------|----------|
| (A)  | $Q_0(p_0,$  | $d_1,\ldots,p_1,\ldots\ldots\ldots,p_S,$ | $\ldots,p_K)$    |               | $=1$     |
| (B)  |             | $Q_1(p_1,d_1,\ldots,p_2,\ldots p_S,$     | $\ldots\ldots,$  | $p_{K+1})$    | $=1$     |
|      |             | $\ldots\ldots\ldots\ldots\ldots\ldots\ldots\ldots\ldots$ | $\ldots\ldots\ldots$ | $\ldots\ldots$ |          |
|      |             | $Q_{H-K}(p_{H-K},d_1,\ldots,p_S,$        | $\ldots\ldots\ldots,$ | $p_H)$       | $=1$     |

the schema[6] in Figure 1. By Lemma 1 we can either eliminate $p_0$ from $Q_0$, or eliminate one of the $p_{S+1},\ldots,p_K$ from $Q_0$, or eliminate one of the $p_{K+1},\ldots,p_H$ from $Q_1,\ldots,Q_{H-K}$.

In the case that neither of the first two cases occur, i.e. when a relation between only the $p_1,\ldots,p_K,d_1,\ldots$ follows from $Q_1 = \cdots = Q_{H-K} = 1$, then we can by using complete induction on the given assumptions, by observing that the relation system $Q_1 = 1,\ldots,Q_{H-K} = 1$ is a system of type (q) with fewer than $H - K + 1$ relations. This resolves the case of $H - K \leq S$, and leaves the case $H - K > S$. We use Lemma 1, in which we identify:

| in Lemma 1        | with | in Lemma 3                        |
|-------------------|------|-----------------------------------|
| $a_1,\ldots$      | "    | $p_0,\ldots,p_S$                  |
| $b_1,\ldots,$     | "    | $d_1,\ldots,$                     |
| $x_1,\ldots$      | "    | $p_{S+1},\ldots,p_{S+K}$          |
| $y_1,\ldots$      | "    | $p_{S+K+1},\ldots,p_H$            |
| the system (A)    | "    | $Q_0 = 1,\ldots,Q_S = 1$          |
| the system (B)    | "    | $Q_{S+1} = 1,\ldots,Q_{H-K} = 1$  |

This gives: either a relation for $p_0,\ldots,p_S,d_1,\ldots$ follows from $Q_0 = \cdots = Q_S = 1$, or one follows for $d_1,\ldots,p_{S+1},\ldots,p_{S+K}$ follows from the same, or a relation follows for $d_1,\ldots,p_{S+1},\ldots,p_{S+K}$ follows from $Q_{S+1} = \cdots = Q_{H-K} = 1$. In all these cases the relation systems $Q_0 = 1,\cdots,Q_S = 1$ resp. $Q_{S+1} = 1,\ldots,Q_{H-K} = 1$ are of "type (q)"[7] and consist of fewer than $H - K + 1$ relations, which means that we can apply complete induction. Thus, Lemma 3 is proved.                                            □

6. *Proof of the Freiheitssatz*: We apply Lemma 3 by proving: if from a finite subsystem of (3) resp. (6) the relation $b_0^n = 1$ resp. $F(\ldots,c_{\nu,k},\ldots) = 1$ follows, then one (in general a different one) of the finite subsystems of (3) resp. (6) has type (q); by Lemma 3 it then follows that in one of the relations in the system (3) resp. (6) one of the generators occurring in this relation can be eliminated, which is precisely what needs to be shown. We first deal with Case I. To deal with this, we first write the finite system of (3), from which $b_0^n = 1$ follows, in the following way:

We have $M \geq 0, N \geq 0, L \geq 0$ such that

$$(\bar{3}) \quad \left.\begin{array}{l} P_0 \equiv P(b_{-N}, b_{-N+1},\ldots,b_{-N+L}) = 1 \\ P_1 \equiv P(b_{-N+1},\ldots\ldots\ldots,b_{-N+L+1}) = 1 \\ \ldots\ldots\ldots\ldots\ldots\ldots\ldots\ldots\ldots\ldots\ldots\ldots\ldots\ldots \\ \ldots\ldots\ldots\ldots\ldots\ldots\ldots\ldots\ldots\ldots\ldots\ldots\ldots \\ P_{M-L+N} \equiv P(b_{M-L},\ldots\ldots,b_M) = 1 \end{array}\right\} \text{ implies that } b_0^n \equiv 1.$$

---

[6] where the generators in the $Q$ are partially permuted compared to their positions in (q).

[7] See Footnote 5.

We certainly have $-N < M$ and $L > 0$, as otherwise $x$ does not occur in $R$. Furthermore, we may assume $M - L \geq 0$, because if a system $(\bar{3})$ with a large value of $M$ does not yield $b_0^n \equiv 1$, then it is also *a fortiori* true for smaller values of $M$. We do not assume that all of $a_{-N}, a_{-N+1}, \ldots, a_{-N+L}$ occur in $P_0$, but we do assume that both $a_{-N}$ and $a_{-N+L}$ do actually occur in $P_0$; and in general both $a_{-N+r}$ and $a_{-N+L+r}$ actually occur in $P_r$.

The smallest of one of the numbers $-N$ and $M$ is non-zero. It is no loss of generality to assume that this is $M$. Were now $N = 0$, then the system $(\bar{3})$ would be of type (q), so we can apply Lemma 3[8]. We can not, however, assume that $N = 0$, and so we make use of Lemma 2, in which we identify:

| in Lemma 2 | with | in the current case |
|---|---|---|
| $a_1, \ldots$ | " | $\ldots, b_0, \ldots, b_{L-1}$ |
| $x_1, \ldots$ | " | $b_L, \ldots, b_M$ |
| $y_1, \ldots$ | " | $b_{-N}, \ldots, b_{-1}$ |
| the system $(\bar{A})$ | " | $P_N = 1, \ldots, P_{M-L+N} = 1$ |
| the system $(\bar{B})$ | " | $P_0 = 1, \ldots, P_{N-1} = 1$ |

This is because if from $(\bar{3})$ a relation for $b_0$ alone follows, then *a fortiori* a relation for $b_0, \ldots, b_{L-1}$ also follows; and Lemma 2 then implies that either from

$$P_0 = 1, \ldots, P_{N-1} = 1 \quad \text{or from} \quad P_N = 1, \ldots, P_{M-L+N} = 1$$

a relation for $b_0, \ldots, b_{L-1}$ must follow.

Now we use Lemma 3, by identifying, for example, in the first case:

| in Lemma 3 | with | in the current case |
|---|---|---|
| $p_0, \ldots, p_S$ | " | $\ldots, b_0, \ldots, b_{L-1}$ |
| $p_{S+1}, \ldots, p_K$ | " | $b_{-1}$ |
| $p_{K+1}, \ldots, p_H$ | " | $b_{-2}, \ldots, b_{-N}$ |
| and $Q_\nu = 1$ | " | $P_{N-\nu-1} = 1$ $\qquad (\nu = 0, \ldots, N-1)$. |

and additionally the special case of Lemma 3, in which the $d_1, \ldots$ do not appear. We find that one of the generators actually occurring in one of the relations $P = 1$ can be eliminated, and the same applies when from $P_N = 1, \ldots, P_{M-L+N} = 1$ a relation for $b_0, \ldots, b_{L-1}$ follows, which is proved in exactly the same way as above. Case II is settled in a very similar manner. We have here $N \leq M$ and $L \geq 0$ such that:

$$(\bar{6}) \left. \begin{array}{l} P_0 \equiv P(\ldots, c_{\nu, l+N}, \ldots; x_{+N}, \ldots, x_{+N+L}) = 1 \\ \cdots\cdots\cdots\cdots\cdots\cdots\cdots\cdots\cdots\cdots\cdots\cdots\cdots \\ P_{M-N} \equiv P(\ldots, c_{\nu, l+M}, \ldots; x_M, \ldots, x_{M+L}) = 1 \end{array} \right\} \text{implies } F(\ldots, c_{\nu, k}, \ldots) \equiv 1.$$

In $P_\nu$, both $x_{N+\nu}$ and $x_{N+L+\nu}$ actually appear. We want to show that already in a single one of the relations in $(\bar{6})$ there is an actually occurring generator $x$ which can

---

[8]Indeed, if $b_0^n \equiv 1$ follows from $P(b_0, \ldots, b_L) = 1, \ldots, P(b_{M-L}, \ldots, b_M) = 1$, then we identify:

| in Lemma 3 | with | in the current case |
|---|---|---|
| $p_0, \ldots, p_S, \ldots$ | " | $b_0$ |
| $p_{S+1}, \ldots, p_K$ | " | $b_1, \ldots, b_L$ |
| $p_{K+1}, \ldots, p_H \cdots$ | " | $b_{L+1}, \ldots, b_M$ |
| $Q_\nu = 1$ | " | $P(b_\nu, \ldots, b_{\nu+L}) = 1$ |

and use the special case of Lemma 3, in which none of the relations (q) contains a generator $d_1, \ldots$ (the $d_1, \ldots$ do not appear whatsoever).

be eliminated. In the case that $N = M$, i.e. when $(\bar{\bar{6}})$ consists only of a single relation, then this is trivially true. We can also dispose of the case when $L = 0$; because if in each relation $(\bar{\bar{6}})$ only a single $x$ appears, then we can use Lemma 2 in the following way: from $(\bar{\bar{6}})$: $P(\ldots, c_{\nu,l+N+t}, \ldots; x_{N+t}) = 1$ (where $t = 0, \ldots, M - N$) a relation for the $c_{\nu,K}$ alone should follow. We identify:

| in Lemma 2 | with | in the current case |
|---|---|---|
| $a_1, \ldots$ | " | $\ldots, c_{\nu,K}, \ldots$ |
| $x_1, \ldots,$ | " | $x_N$ |
| $y_1, \ldots,$ | " | $x_{N+1}, \ldots, x_M$ |

and obtain: either we can eliminate $x_N$ from $P_0 = 1$, or else from the system $(\bar{\bar{6}}\text{a})$: $P_1 = 1, \ldots, P_{M-N} = 1$ a relation for the $\ldots, c_{\nu,K}, \ldots$ alone follows. The system $(\bar{\bar{6}}\text{a})$ is constructed such that it contains one less relation than $(\bar{\bar{6}})$. If $(\bar{\bar{6}}\text{a})$ contains more than one relation, then we can by repeated application of Lemma 2 prove, that one of the relations from $(\bar{\bar{6}})$ contains a generator $x$ which can be eliminated.

We now consider the case when $L > 0$. If from $(\bar{6})$ a relation for the $\ldots, c_{\nu,K}, \ldots$ follows, then *a fortiori* such a relation must also follow for the $\ldots, c_{\nu,K}, \ldots$ and $x_N, x_{N+1}, \ldots, x_{N+L-1}$ together. We now use Lemma 3, identifying as follows:

| in Lemma 3 | with | in the current case | |
|---|---|---|---|
| $d_1, d_2, \ldots$ | " | $\ldots, c_{\nu,K}, \ldots$ | |
| $p_0, \ldots, p_S$ | " | $x_N, \ldots, x_{N+L-1}$ | |
| $p_{S+1}, \ldots, p_K$ | " | $x_{N+L}$ | |
| $p_{K+1}, \ldots, p_H$ | " | $x_{N+L-1}, \ldots, x_M$ | |
| the relation $Q_\nu = 1$ | " | $P_{N+\nu} = 1$ | $(\nu = 0, \ldots, M - N)$. |

We obtain as a consequence, that also in this case in one of the relations of $(\bar{6})$ there is an actually occurring generator $x$ which can be eliminated. This completes the reduction of the proof of the *Freiheitssatz* to the proofs of Lemma 1 and Lemma 2.

## §4. Proof of the two lemmas

1. **Lemma 1**. In the proof of Lemma 1, we will replace each of the four generating systems – for the sake of brevity – by a single generator. The course of the proof will not be affected by this. We thus have $a, b, x, y$ and the system:

$$\text{(A)} \quad \{P_\mu(a, b, x) = 1\} \qquad \text{and} \qquad \text{(B)} \quad \{Q_\nu(b, x, y) = 1\}.$$

If now from (A) and (B) there follows a relation $R(a, b) = 1$ only between the $a$ and $b$, then there is an identity of the form:

$$\text{(1)} \qquad R(a, b) \equiv \prod_{i=1}^{h} T_i K_i^{(\beta_i)} T_i^{-1},$$

where $\beta_i$ takes one of the values $1, 2$, such that for $\beta_i = 1$ we have that $K_i^{(1)}(a, b, x)$ is a word over $a, b, x$ such that $K_i^{(1)}(a, b, x) = 1$ follows from (A), and when $\beta_i = 2$ we have that $K_i^{(2)}$ is a word over $b, x, y$ equal to $1$ using the relations (B). For example, we can choose to take as $K_i^{(\beta_i)}$ the $P_\mu^{\pm 1}$ and $Q_\nu^{\pm 1}$, but this more general representation has its advantages. Before showing this, however, we give some terminology:

In the right-hand side of (1) the word $F_i \equiv T_i K_i^{(\beta_i)} T_i^{-1}$ is called the $i$th "*factor*", $T_i$ is a "*conjugator*", and $K_i^{(\beta_i)}$ is a "*kernel*". The kernel is thus defined as what remains when one performs all possible absorptions in the cyclically written factors.

*We will postulate the following:*

(I) *No factor allows any absorptions.* If necessary, this can always be achieved by replacing one kernel with that obtained from it by performing all possible cyclic deletions on it.

(II) *The right-hand side of* (1) *is a representation of $R$ which has the minimal number of factors, considered over all possible such representations.*

(III) *Of all the possible representations in* (II) *for the right-hand side of* (1), *we have chosen one with the minimal number of $a$- and $y$-letters.*

From the assumptions in the statement of the lemma (§3, Nr. 5) we have:

Each kernel contains either $a$-letters or $y$-letters, but not both. And: there are kernels in which $y$-letters really occur. The $a$-letters and $y$-letters are called the *characteristic letters* of the kernels.

From the postulates (I) to (III) and the aforementioned assumptions we will now derive a contradiction. First, some small and helpful remarks:

($\alpha$) Between two neighbouring factors with $y$-containing kernels, we cannot remove the $a$-letters in the conjugators – otherwise, we could combine the two factors into a single one, which would decrease the number of factors in contradiction with (II).

($\beta$) The characteristic letters of a kernel can never be absorbed with a neighbouring kernel – because if both kernels have the same characteristic letters, then one can in this case again remove the number of factors.

($\gamma$) No factor can have more than half of its $a$-letters and $y$-letters be absorbed by the *conjugators* of a neighbouring factor; (this case occurs if and only if more than half of the characteristic letters of the kernel of the conjugators of the neighbouring factor are absorbed), because in this case we can reduce the number of $a$-letters and $y$-letters – contradicting (III)[9].

We now write $R(a,b) \equiv F_1 F_2 \cdots F_h$. On the right-hand side all $y$-letters must be eliminated. Absorption can only happen when two factors are merged together. It is obviously not possible, after all absorptions between the kernels of each factor have been carried out, that only characteristic letters remain, since then no further absorptions would be possible, and yet some factors would contain $y$-letters in their kernel, which would thus not be eliminated.

On the other hand: if one performs absorptions as far as possible at the point where to factors touch (without worrying about the other places), then because of ($\beta$) and $\gamma$) we must have that in no kernel can more than half of its characteristic letters be absorbed, and that, if any factor as a result of such a process will have all of its characteristic letters absorbed, then these must be absorbed with the *conjugators* of the neighbouring factors. Thus we have:

---

[9] Let $F_2 \equiv \varphi \overline{\varphi}$, where $\varphi$ contains more than half of the characteristic letters of the kernel $K_2^{(\beta_2)}$ (and thus more than $\overline{\varphi}$), and let $F_1 \equiv T_1 K_1^{(\beta_1)} T_1^{-1} \equiv \varphi \tau_1 K_1^{(\beta_1)} \tau_1^{-1} \varphi^{-1}$ (where on the right-hand side no absorptions occur). Then we have:

$$F_1 F_2 \equiv T_1 K_1^{(\beta_1)} T_1^{-1} F_2 \equiv \varphi \tau_1 K_1^{(\beta_1)} \tau_1^{-1} \varphi^{-1} \varphi \overline{\varphi} \equiv \varphi \overline{\varphi} (\overline{\varphi}^{-1} \tau_1) K_1^{(\beta_1)} (\overline{\varphi}^{-1} \tau_1)^{-1} \equiv F_2 \overline{T}_2 K_2^{(\beta_2)} \overline{T}_2^{-1}$$

where $\overline{T}_2$ contains fewer $a$- and $y$-letters than $T_2$.

There are factors whose $a$-letters and $y$-letters are absorbed in exactly equal parts by the conjugators of their neighbouring factors. Such factors will be called "***permutable***" on the basis of the following proposition.

If in $F_1 F_2 F_3$ the factor $F_2$ is permutable, then $F_1 F_2 F_3 \equiv F_2 \overline{F}_1 F_3 \equiv F_1 \overline{F}_3 F_2$ where $\overline{F}_1$ and $\overline{F}_3$ are conjugates of $F_1$ resp. $F_3$, containing equally many $a$- and $y$-letters as $F_1$ resp. $F_3$.[10] Thus: if a factor is permutable, then its neighbours are not.[11]

We now only have to prove the following: we can, without violating the assumptions made, transform the right-hand side of (1) such that it does not contain any permutable factors.

To do this, we will assume that the representation (1) not only satisfies the postulates (I) to (III), but also that the representation:

(IV) contains a *minimal number of factors*,

and, among all possible representations satisfying this, choosing one for which:

(V) *the first permutable factor has the smallest possible index.*

Let now $F_r$ in $F_1 \cdots F_{r-2} F_{r-1} F_r F_{r+1} \cdots F_h$ be the first permutable factor. We permute $F_r$ with $F_{r-1}$ and obtain: $R \equiv F_1 \cdots F_{r-2} F_r \overline{F}_{r-1} F_{r+1} \cdots F_h$ (where $\overline{F}_{r-1}$ is a conjugate of $F_{r-1}$, which really contains both $a$-letters and $y$-letters). It is possible that $F_r$ is now permutable again. If then another factor had not become permutable, then (V) would have been violated. Therefore, since nothing else has changed and $\overline{F}_{r-1}$ cannot be permutable, $F_{r+1}$ must be permutable.[12] But the same must be the case if $F_r$ is not permutable, as otherwise (IV) would be violated.

Thus, when we permute $F_{r+1}$ with $\overline{F}_{r-1}$, we have:

$$R \equiv F_1 \cdots F_{r-2} F_r F_{r+1} \overline{\overline{F}}_{r-1} F_{r+2} \cdots F_h.$$

Since $F_{r+1}$ is no longer permutable[13] we can proceed as above to deduce the existence and permutability of $F_{r+2}$. But $h$ is finite. Therefore, we have proved a contradiction, and therefore also the non-existence of permutable factors. $\qquad \square$

We have thus proved Lemma 1. The proof of Lemma 2 is carried out by the same method, and we will therefore only provide a sketch of the proof.

2. **Lemma 2**. Again, we will return to representing the generating systems by a single letter. We thus have

$$(\overline{\mathrm{A}}) \quad \{P_\mu(a, x) = 1\} \qquad \text{and} \qquad (\overline{\mathrm{B}}) \quad \{Q_\nu(a, y) = 1\}.$$

If from this a relation $R(a) = 1$ for $a$ follows (which is not identically true), then we have the following identity:

$$(2) \qquad R(a) \equiv \prod_{i=1}^{h} T_i K_i^{(\beta_i)} T_i^{-1}$$

---

[10] The proof is like the proof of $(\gamma)$.

[11] It cannot be the case that both $F_1$ absorbs half of the $a$-letters and $y$-letters of the conjugator of $F_2$ and that $F_2$ absorbs half of the $a$-letters and $y$-letters of the conjugator of $F_1$. Indeed, in the first case, the conjugator of $F_2$ must contain more letters than $F_1$, since it completely absorbs the latter and, moreover, also absorbs letters from the kernel of $F_1$. In the second case, analogous reasoning gives that $F_1$ must contain more letters than $F_2$.

[12] Hence, $F_{r+1}$ exists.

[13] Because we assumed that in $F_1 \cdots F_{r-1} F_r F_{r+1} \cdots$ the factor $F_r$ is permutable.

where $\beta_i = 1, 2$, and $K_i^{(1)}(a, x)$ resp. $K_i^{(2)}(a, y)$ is equal to 1 on the basis of $(\overline{A})$ resp. $(\overline{B})$. We make the following postulates for the representation of $R$ in the right-hand side of (2):

(I′)  No factor allows any absorptions.

(II′)  In the right-hand side of (2) there is a minimal number of factors.

(III′)  In the right-hand side of (2) there is a minimal number of $x$- and $y$-letters.

From the assumptions in the statement of the lemma it now follows: each kernel $K_i^{(\beta_i)}$ contains, depending on whether $\beta_i = 1$ or $= 2$, either $x$-letters or $y$-letters, but not both. The $x$-letters resp. $y$-letters, which appear in some kernel, are called the characteristic letters of the corresponding factor.

Based on these assumptions, in the right-hand side of (2) there can only be absorptions if two factors are merged, and it must not be the case that, after carrying out all such absorptions, characteristic letters remain, since then in the right-hand side of (2) not all $x$-letters or $y$-letters could be eliminated. We then proceed as in the proof of Lemma 1; one has to only replace the words $a$- and $y$-letters by $x$- and $y$-letters, respectively.                                                                                      □

**§5. Generalization of Lemma 1. Related lemmas. Independence theorems.**

The following results are associated to Lemmas 1 and 2; however, we will not use them until the *second part*.

1.  **Generalization of Lemma 1.**  Let $a, b, x, y$ be four systems of generators, between which there are two systems of relations:

(A)  $\{P_\mu(a, b, x) = 1\}$       and   (B)  $\{Q_\nu(b, x, y) = 1\}$.

and if from (A) a system of relations follows between $a$ and $b$, but no relation holds for $b$ alone, and if no relation for $b$ and $x$ alone follows from either (A) or (B) alone – then all relations for $a$ and $b$ alone which follow from (A) and (B) follows from (A) alone.

*Proof.* The proof of this is already contained[14] in the proof of Lemma 1; because in

$$(1) \qquad\qquad R(a, b) \equiv \prod_{i=1}^{h} T_i K_i^{(\beta_i)} T_i^{-1}$$

($\beta_i = 1, 2$; $K_i^{(1)}(a, b, x) = 1$ on the basis of (A) and $K_i^{(2)}(b, x, y) = 1$ on the basis of (B)) on the right-hand side there can be no factor whose kernel contains $y$-letters, when the postulates (I) to (V) are made (§4, Nr. 1); but by assumption all $K_i^{(2)}$ contain $y$-letters, and therefore no $K_i^{(2)}$ can appear whatsoever.                       □

2.  **Lemma 4.**  *Consider four systems of generators, which we will represent by a single generator each. We thus have* $a, b, x, y$ *and two systems of relations between them:*

(A₁)  $\{P_\mu(a, b, x) = 1\}$       *and*   (B₁)  $\{Q_\nu(b, x) = 1\}$.

----

[14]The actual reason for this is as follows: when proving Lemma 1, one first proves – regardless of whether from (A) a relation for $a$ and $b$ (in which $a$ actually appears, as otherwise it is a relation for $b$ and $x$) follows – that, based on the other assumptions, the relations following from (B) are not needed. Since, by assumption, from (A) no $R = 1$ follows, we are done.

*Suppose that from neither* (A₁) *nor* (B₁) *alone a relation for* $x$ *follows, and suppose that all relations between* $a$ *and* $x$ *which follow from* (A₁) *are defined by the (possibly infinite) system of relations* (C₁) $\{S_\varrho(a,x) = 1\}$.

*Then the following holds: all relations* $R(a,b) = 1$ *between* $a$ *and* $b$ *which follow from* (A₁) *and* (B₁)*, follow already from* (C₁) *and* (B₁)*.*

*Proof.* We have the identity

$$(1_1) \qquad\qquad R(a,b) \equiv \prod_{i=1}^{h} T_i K_i^{(\beta_i)} T_i^{-1}; \quad \beta_i = 1, 2$$

where $K_i^{(1)}(a,x,y)$ is 1 on the basis of (A₁), and $K_i^{(2)}(b,x) = 1$ on the basis of (B₁). Each $K_i^{(\beta_i)}$ contains, by assumption, $a$-letters or $b$-letters – but not both – or only $y$-letters and no $b$-letters. The $a$-letters and $y$-letters resp. the $b$-letters, which appear in the kernel $K_i^{(\beta_i)}$ of a factor $T_i K_i^{(\beta_i)} T_i^{-1}$, are called characteristic letters of the kernels resp. the factors. Just as in Lemma 1, we add the following postulates for the right-hand side of the representation (1₁):

(I′)   No factor allows any absorptions.

(II′)  In the right-hand side of (1₁) there is the smallest number of factors.

(III′) Of all representations of $R$ with a minimal number of factors, the right-hand side of (1₁) has a minimal number of $a$-, $y$-, and $b$-letters.

If we call a factor permutable when it can absorb at least half of its $a$-, $y$-, and $b$-letters with the conjugator of its neighbouring factor, then we postulate further:

(IV′) In the right-hand side of (1₁) there is a minimal number of permutable factors;

(V′)  and the first such permutable factor has the smallest possible index.

One can now show – word for word as in Lemma 1 –: no absorption is possible except where two factors touch. If one performs all possible absorptions at one such place, then no factor can have more than half of its $a$-, $b$-, and $y$-letters be absorbed. Since there are no permutable factors, it follows: if one performs all possible absorptions in (1₁), then there will still remain characteristic letters in the kernel. Since there can be no $y$-letters, it easily follows that: the right-hand side of (1₁) contains, when the postulates (I′) to (V′) hold, no $y$-letters[15]. In particular, $K_i^{(\beta_i)}$ therefore cannot contain any such letters, and since $R = 1$ follows from $K_i^{(\beta_i)} = 1$, the lemma is proved.                                                                                  □

3. The following **Lemma 5** requires a new idea for its proof. If we are given three systems of generators – which, in the sequel, will be denoted by a single generator – say, $a, b, t$, and two systems of relations between them:

$$(\alpha) \qquad \{P_\mu(a,t) = 1\} \qquad \text{and} \qquad (\beta) \qquad \{Q_\nu(b,t) = 1\}.$$

Suppose that from neither of these systems alone does a relation for $t$ alone, or a relation for $a$ or $b$ alone, follow; but that both together result in a relation $R(a,b) = 1$.

---

[15] *Proof:* Indeed, by the above reasoning, the $y$-letters in each factor must be absorbed by the conjugators of the neighbouring factors; thus, if some factor contains a $y$-letter, then this is also true of the conjugators of all other factors. But $T_1$ and $T_h^{-1}$ do not contain any $y$-letters, as these cannot be cancelled.

We have the identity

$$(1_2) \qquad R(a,b) \equiv \prod_{i=1}^{h} T_i K_i^{(\beta_i)} T_i^{-1}$$

with $\beta_i = 1, 2$ and $K_i^{(1)}(a,t) = 1$ on the basis of $(\alpha)$, and $K_i^{(2)}(b,t) = 1$ on the basis of $(\beta)$. Each $K_i^{(\beta_i)}$ contains, by assumption, either $a$-letters or $b$-letters, but not both. If we call the $a$- resp. $b$-letters as the *characteristic letters* of $K_i^{(1)}$ resp. $K_i^{(2)}$, then we have the following on the basis of the often used reasoning in the proof of Lemma 1:

(I) *The right-hand side of* $(1_2)$ *has its factors* $T_i K_i^{(\beta_i)} T_i^{-1}$ *chosen such that when all possible absorptions between each factor is performed, then there remain characteristic letters in the kernels.* Furthermore, we postulate:

(II) If the right-hand side of $(1_2)$ is chosen such that (I) holds, then at least one of the expressions $K_i^{(\beta_i)}$ is obtained from a word $A(a)\Theta(t)$ resp. $B(b)\Theta(t)$ by conjugation $(A, B, \Theta$ are, respectively, words over $a, b, t$, which are not identically 1), or, stated otherwise: *from $(\alpha)$ or $(\beta)$ a relation $A^{-1}(a) = \Theta(t)$ resp. $B^{-1}(b) = \Theta(t)$ follows.*

And finally it follows directly from (I) that:

(III) *If* (I) *holds, then the number of kernels* $K_i^{(1)}$ *resp.* $K_i^{(2)}$ *in the right-hand side* $(1_2)$ *is less than or equal to the number of $a$- resp. $b$-letters in $R(a,b)$, when this is written such that no absorptions are possible.*

Number (III) is used in the case that $R$ consists of only $a$- and $b$-letters (§7). All that remains to be proved is (II). We do this in the following way: because of the assumption of the validity of (I), we know that upon performing all absorptions in the right-hand side of $(1_2)$, no more $t$-letters can appear there, and one may also make the assumption that in the right-hand side of $(1_2)$ the number of factors is the least possible, since this assumption can always be made by the very existence of a representation of the right-hand side of $(1_2)$ satisfying (I). We now define: if one performs all possible absorptions at a place where two factors meet (where both factors contain $t$-letters in their kernels!), and there remains from the kernel of the second only a word $A_1(a)\Theta(t)$ resp. $B_1(A)\Theta_1(t)$ (where $\Theta_1 \equiv 1$ is possible, but by (I) this is not possible for $A_1$ or $B_1$), then we say that the second factor is "*destructively absorbed*" into the first. We now prove: when no kernel is a conjugate of a word $A(a)\Theta(t)$ resp. $B(b)\Theta(t)$, and when the $i$th factor is not destructively absorbed into the $(i-1)$st, then the $(i+1)$st factor is not absorbed into the $i$th factor. With this we will be done, because by assumption the first factor does not have a kernel of the form $A_1(a)\Theta(t)$ resp. $B_1(b)\Theta_1(t)$, and since the last factor cannot be destructively absorbed, we must in this case have that there must remain $t$-letters on the right-hand side of $(1_2)$ even after performing all possible absorption, which contradicts our assumptions. Hence – to prove the above statement – let $F_{i-1}, F_i, F_{i+1}$ be three consecutive factors in the right-hand side of $(1_2)$. The kernel of $F_i$ really contains $a$. If one performs all possible absorptions on the boundary between $F_{i-1}$ and $F_i$, then there remains in $F_i$ a piece $S_i$ of the following form:

$$S_i \equiv A_0(a)\Theta_1(t)A_1(a)\Theta_2(t)A_2(a,t)T_i^{-1},$$

where $T_i$ is the conjugator of $F_i$, $A_0\Theta_1 A_1\Theta_2 A_2$ belongs to the kernel of $K_i^{(\beta_i)}$, and $\Theta_1$ and $A_1$ are not $\equiv 1$. In $\Theta_1 A_1 \Theta_2 A_2 T_i^{-1}(T_{i+1} K_{i+1}^{\beta_{i+1}} T_{i+1}^{-1})$ both $\Theta_2$ and $A_1$ must

absorb with $T_{i+1}$[15a], since $K_{i+1}^{(\beta_{i+1})}$ cannot have any of its characteristic letters absorb with $K_i^{(\beta_i)}$, and one can thus in this case reduce the number of factors. If $\Theta_1$ is also fully absorbed with $T_{i+1}$, then we are done[16]. In the case that $\Theta_1$ entirely or partially is absorbed with $K_{i+1}^{(\beta_{i+1})}$, we proceed as follows: since $K_{i+1}^{(\beta_{i+1})}$ is not a conjugate of $A\Theta$ resp. $B\Theta$, it can be written as:

$$K_{i+1}^{(\beta_{i+1})} \equiv \overline{\Theta}_1 C_1 \overline{\Theta}_2 C_2 H,$$

where $\overline{\Theta}_1$ and $\overline{\Theta}_2$ are words over $t$, $C_1$ and $C_2$ are over $a$ or $b$ (corresponding to whether $\beta_{i+1} = 1$ or $= 2$), and $H$ (which is possibly $\equiv 1$) is written over $t$ and $a$ resp. $b$. The word $K_{i+1}^{(\beta_{i+1})}$ must begin with $t$-letters, because otherwise $\Theta_1$ cannot be absorbed. But even if $\overline{\Theta}_1$ is completely absorbed with $\Theta_1$, then, as we can see, $K_{i+1}^{(\beta_{i+1})}$ cannot be destructively absorbed.                                          □

# Second Part:
# Applications of the *Freiheitssatz*

Even though the *Freiheitssatz* seems almost trivial, it is nevertheless an extraordinarily powerful tool, as will be seen in the following sections.

### §6. Equivalent relations. Primitive elements.
### A strong form of the *Freiheitssatz*.

1. *Preliminary remark: the "strong form".* To use the *Freiheitssatz*, it is first useful to express it in the following form, which we will make constant use of in the sequel.

Given generators: $b_0, b_1, \ldots, b_M$ and $c$, and a system of relations between them:

$$(1) \qquad \left.\begin{aligned} &S_0(c; b_0, \ldots, b_K) = 1 \\ &\cdots\cdots\cdots\cdots\cdots \\ &S_{M-K}(c; b_{M-K}, \ldots, b_M) = 1 \end{aligned}\right\} K \leq M,$$

where the $S$, written cyclically, admit no absorptions, and each $S_i$ really contains the generators $b_i$ and $b_{i+K}$. Then we have: if $0 \leq H \leq L \leq M$ are integers, and if from (1) there follows a relation for $c$ and $b_H, \ldots, b_L$ alone, then this relation follows from *the* relations in (1) which only contain $c$ and $b_H, \ldots, b_L$.

The proof of this is as follows: if the system (1) consists of only a single relation, then this is just the statement of the *Freiheitssatz*. Using the *generalization* of Lemma 1 (§5) and complete induction, the claim made follows in general by using exactly the same proof idea by which Lemma 3 (§3) is derived from Lemma 1 (§3). Note also that the *Freiheitssatz* in the above form remains true even if one replaces the generator $c$ by a system of generators $c_1, c_2, \ldots$.                                          □

2. *Equivalent relations.* It is now not difficult to prove the following theorem: *if* $R_1(a_1, a_2, \ldots, a_n) = 1$ *and* $R_2(a_1, \ldots, a_n) = 1$ *are two relations over the generators* $a_1, \ldots, a_n$ *such that* $R_1 = 1$ *follows from* $R_2 = 1$, *and vice versa, then we have:* $R_1$ *is conjugate to* $R_2$. We say that $R_1 = 1$ and $R_2 = 1$ are equivalent.

---

[15a]Indeed, $\Theta_1$ must also be absorbed somehow.

[16]As can be seen from the form of $K_{i+1}^{(\beta_{i+1})}$ that we are about to present.

The proof of the above statement is by complete induction in the following manner: if $R_1$ or $R_2$ only contain a single generator, then the theorem is an immediate consequence of the *Freiheitssatz*. More generally – slightly modified – the words $R_1$ and $R_2$ are, cyclically rewritten, identical, and we conclude further: when the theorem is true for all pairs of equivalent relations $P_1 = 1, P_2 = 1$, for which $P_1$ resp. $P_2$, which when cyclically written[17] contain fewer letters than $R_1$ resp. $R_2$, cyclically written[17], then it also holds for the pair $R_1 = 1, R_2 = 1$.

*Proof.* We have two identities:

$$\text{(2)} \qquad R_1 \equiv \prod_{i=1}^{h} T_i R_2^{e_i} T_i^{-1} \quad (e_i = \pm 1)$$

$$\text{(}\overline{2}\text{)} \qquad R_2 \equiv \prod_{i=1}^{\overline{h}} \overline{T}_i R_1^{\overline{e}_i} \overline{T}_i^{-1} \quad (\overline{e}_i = \pm 1).$$

By abelianizing it follows that $R_1$ and $R_2$ has the same exponent sum in all generators. Since $R_1$ and $R_2$, written cyclically, admit no absorptions, $R_1$ contains precisely the same generators that really appear in $R_2$, and vice versa. (*Freiheitssatz.*) Suppose $a_1$ appears in both $R_1$ and $R_2$. We may assume that $a_1$ has exponent sum zero $R_1$ (and hence also in $R_2$). This is because: if in $R_1$ only $a_1$ appears, then we are done (see the above). If in $R_1$ some other generators appear, then if one of them has exponent sum zero, then we call this $a_1$. Otherwise, in $R_1$ there are two generators, say $a_1$ and $a_2$, which have non-zero exponent sums $s_1$ resp. $s_2$. If we perform the substitution $a_1 = b_1^{-s_2}$, $a_2 = b_1^{s_1} b_2$, the number of $a_2, a_3, \ldots$-letters in $R_1$ will be equal to the number of $b_2, a_3, \ldots$-letters after the introduction of the letters $b_1$ and $b_2$, and $b_1$ has in $R_1$ exponent sum zero[18]. We now denote $b_1$ resp. $b_2$ by $a_1$ and $a_2$. We set $a_1^K a_\nu a_1^{-K} = a_{\nu,K}$ for $\nu = 2, \ldots, n$ and all integers $K$. This changes the words $R_1$ resp. $R_2$ into words $P_1(\ldots, a_{\nu,K}, \ldots)$ resp. $P_2(\ldots, a_{\nu,K}, \ldots)$, which contain fewer letters than $R_1$ resp. $R_2$. For all integer $\lambda$, we set:

$$P_1(\ldots, a_{\nu,K+\lambda}, \ldots) = P_{1,\lambda}$$
$$P_2(\ldots, a_{\nu,K+\lambda}, \ldots) = P_{2,\lambda}.$$

If we let $t_i$ resp. $\tau_i$ be the exponent sum of $a_i$ in $T_i$ resp. $\overline{T}_i$, then the identities (2) and ($\overline{2}$) become:

$$\text{(3)} \qquad P_{1,0} \equiv \prod_{i=1}^{h} T_i P_{2,t_i}^{e_i} T_i^{-1}$$

$$\text{(}\overline{3}\text{)} \qquad P_{2,0} \equiv \prod_{i=1}^{\overline{h}} \overline{T}_i P_{1,\tau_i}^{\overline{e}_i} \overline{T}_i^{-1},$$

where the $T_i$ and $\overline{T}_i$ are words in the $a_{\nu,K}$. Using the strong form of the *Freiheitssatz* it follows[19] that there is some $\lambda_0$ such that $P_{1,0} = 1$ and $P_{2,\lambda_0} = 1$ (resp. $P_{2,0} = 1$ and $P_{1,-\lambda_0} = 1$) are equivalent relations. Thus we are done, since we assumed by

---

[17] Here it is assumed that the considered words, when cyclically written, admit no absorptions.

[18] The fact that this substitution in $R_1$ could introduce more $b_1$-letters than there were $a_1$-letters does not matter, since we will conjugate by $b_1$ later.

[19] There must be a $P_{2,\lambda_0}$ which (by (3)) which only contains generators that occur in $P_{1,0}$, and by ($\overline{3}$) there must be a $P_{1,\lambda_1}$ which contains only generators occurring in $P_{2,0}$. From this it follows that *all* generators

induction that the theorem holds for $P_{1,0} = 1, P_{2,\lambda_0} = 1$, and so it immediately also follows for $R_1 = 1, R_2 = 1$. $\qquad\qquad\square$

3. *Roots of* $aba^{-1}b^{-1}$. As a further application of the *Freiheitssatz*, we will prove a familiar theorem[20] on the automorphisms of the free group on two generators. *If $a$ and $b$ are generators of a free group, then two words $\alpha$ and $\beta$ on $a$ and $b$ form a pair of connected primitive elements (i.e. they can be considered as new generators), when identically in $a, b$ we have:*

$$(4) \qquad\qquad \alpha\beta\alpha^{-1}\beta^{-1} \equiv Taba^{-1}bT^{-1}.$$

The proof of this result is done using the following theorem: *a root of* $aba^{-1}b^{-1}$ *is either a primitive element or a conjugate of* $aba^{-1}b^{-1}$. The proof of this goes as follows: if $R = 1$ is a root of $aba^{-1}b^{-1}$, then we have two possible cases: (1) the letter $b$ has exponent sum zero in $R$. In this case, we set $b^k a b^{-k} = a_k$. By the often used conclusion we have: if $R$ is rewritten over the $a_k$, as $P(\ldots, a_k, \ldots)$, then $a_0 a_1^{-1} = 1$ follows from a finite number of the relations of the system: $P(\ldots, a_k, \ldots) = 1, P(\ldots, a_{k+1}, \ldots) = 1, \ldots$. The strong form of the *Freiheitssatz* now gives: if $P$, written cyclically, does not allow for any absorptions, then either $a_0 a_1^{-1} = 1$ follows from a relation $P(a_0, a_1) = 1$, or from two relations $P(a_0) = 1, P(a_1) = 1$. The second case gives $P(a_0) \equiv a_0^{\pm 1}$ and $R \equiv Ta^{\pm 1}T^{-1}$, and thus: $R$ is a conjugate of $a^{\pm 1}$.

The first case thus gives, when we set $a_0 a_1^{-1} = b_0, a_1 = b_1$: that $b_0 = 1$ follows from $Q(b_0, b_1) = 1$ (where $Q(b_0, b_1) = P(a_0, a_1)$). From the *Freiheitssatz* we have $Q \equiv T_2 b_0^{\pm 1} T_2^{-1}$, and hence:

$$P \equiv T_1(a_0 a_1^{-1})^{\pm 1}T_1^{-1}, \quad R \equiv Taba^{-1}b^{-1}T^{-1}.$$

The case when both $a$ and $b$ have a non-zero exponent sum in $R$ can be reduced to this first case. Indeed, if $s_1$ resp. $s_2$ denotes the exponent sum of $a$ resp. $b$ in $R$, and we let $d$ be the greatest common denominator of $s_1$ and $s_2$, then we can from the coprime numbers $\sigma_1 = \frac{s_1}{d}$ and $\sigma_2 = \frac{s_2}{d}$ find a pair of connected primitive elements $\gamma$ and $\delta$, such that $\gamma$ has exponent sum $\sigma_1$ in $a$ resp. $\sigma_2$ in $b$.[21] Since $\gamma\delta\gamma^{-1}\delta^{-1} = 1$ is a root of $aba^{-1}b^{-1} = 1$, and vice versa, it holds that: if one expresses $R$ in $\gamma$ and

---

$\delta$, i.e. if $R(a, b) = P(\gamma, \delta)$, then $P(\gamma, \delta) = 1$ is a root of $\gamma\delta\gamma^{-1}\delta^{-1} = 1$. Since $\delta$ has exponent sum zero in $P$, is hence from the first part of the proof $P$ a conjugate of $\gamma$, and hence a primitive element.

Now it is easy to prove that (4) holds if and only if $\alpha$ and $\gamma$ are connected primitive elements. First: if $\alpha$ and $\beta$ are connected primitive elements, then $\alpha\beta\alpha^{-1}\beta^{-1} = 1$ implies that the group is abelian; and hence $aba^{-1}b^{-1} = 1$ holds, and since $\alpha\beta\alpha^{-1}\beta^{-1}$ has exponent sum zero in $a$ and $b$, it is a conjugate of $aba^{-1}b^{-1}$. On the other hand: if $\overline{\alpha}$ and $\overline{\beta}$ are two words such that

$$(\overline{4}) \qquad \overline{\alpha}\overline{\beta}\overline{\alpha}^{-1}\overline{\beta}^{-1} \equiv T(aba^{-1}b^{-1})T^{-1}$$

then either $\overline{\alpha}$ or $\overline{\beta}$ is a primitive element, since both are roots of $aba^{-1}b^{-1} = 1$, and both cannot be a conjugate of $aba^{-1}b^{-1}$.[22] Thus, $\overline{\alpha}$ is a primitive element, and will henceforth be denoted $\alpha$. It remains to show that $\overline{\beta}$ is a primitive element connected to it. We choose a primitive element $\beta$ connected to $\alpha$. Then $\overline{\beta}$ is a word in $\alpha$ and $\beta$, and the identity $(\overline{4})$ becomes the following identity in $\alpha$ and $\beta$:

$$(5) \qquad \alpha\overline{\beta}\alpha^{-1}\overline{\beta}^{-1} \equiv T\alpha\beta\alpha^{-1}\beta^{-1}T^{-1}.$$

Let $s$ be the exponent sum of $\alpha$ in $\beta$, and $t$ that of $\alpha$ in $T$. If we set $\alpha^K\beta\alpha^{-K} = \beta_K$, then the identity (5) becomes an identity in the $\beta_K$ as:

$$(6) \qquad \overline{\overline{\beta}}(\dots, \beta_{K+1}, \dots)\overline{\overline{\beta}}^{-1}(\dots, \beta_K, \dots) \equiv \overline{T}\beta_{t+1}\beta_t^{-1}\overline{T}^{-1},$$

where $\overline{\overline{\beta}}(\dots, \beta_K, \dots)\alpha^s = \overline{\beta}$.

As one can immediately see, one cannot on the left-hand side of (6) – when one writes it cyclically – delete the first letter with the last; because if $\overline{\overline{\beta}}(\dots, \beta_{K+1}, \dots)$ begins with $\beta_L$, then $\overline{\overline{\beta}}^{-1}(\dots, \beta_K, \dots)$ ends with $\beta_{L-1}^{-1}$. Hence we either have $\overline{T} \equiv 1$ or $\overline{T} \equiv \beta_t^{+1}$ or $\overline{T} \equiv \beta_{t+1}^{-1}$, whence it follows that: $\overline{\beta} \equiv \alpha^{s_1}\beta^{\pm 1}\alpha^{s_2}$, where $s_1 + s_2 = s$. Since $\alpha$ and $\beta$ are connected primitive elements, the same is also true of $\alpha$ and $\overline{\beta}$.   □

## §7. Examples of finding all roots of a given word

1. *Preliminary remarks:* The roots of $aba^{-1}b^{-1}$ could be determined relatively easily, since commutativity of the generators is a characteristic (of a particular type that is independent of the presentation) property of the group. The word $aba^{-1}b^{-1}$ has infinitely many pairwise non-conjugate roots; this is essentially due to the fact that $aba^{-1}b^{-1}$ has exponent sum zero in $a$ and $b$, so that the exponent sums of $a$ and $b$ in the roots of $aba^{-1}b^{-1}$ have no essential restrictions.

By contrast, in the following examples we will treat roots of a "simple" word, in which $a$ or $b$ has non-zero exponent sum. Our treatment can only be carried out in full in very special cases, e.g. the question of the roots of $a^2b^p$ (where $p$ is a prime number), or $a^2b^{2p}$.

The tools used in the sequel can be used in other cases, too, but are insufficient to deal with the general case.

2. *The approach:* We are seeking the roots of $\overline{a}^n\overline{b}^m$. We set $\overline{a} = ab^{-m}$ and $\overline{b} = b^{-n}$. (We may assume that $m \neq 0, n \neq 0$, as otherwise everything is trivial). Each root

---

[22] One sees this immediately when one sets $b^K ab^{-K} = a_K$ (for $K = 0, \pm 1, \dots$) and $(\overline{4})$, transforms into an identity in the $a_K$ and makes it abelian.

of $\bar{a}^n \bar{b}^m$ corresponds to one of $[ab^{-m}]^n b^{-mn}$, although not vice versa[23]. The word $W \equiv (ab^{-m})^n b^{-mn}$ has exponent sum zero in $b$, although not in $a$; the same is hence also true of every root $R$ of $W$. In the sequel, we will only consider those roots $R$ of $W$ which admit no absorption when written cyclically; conjugates of such words should not be listed specifically; and similarly if two roots are inverses of one another, then we will only list one.

The word $R(a, b)$ becomes, when we set $a_K = b^K a b^{-K}$ for all integers $K$, a word $Q(\dots, a_K, \dots)$. *We then have that*

$$W = a_0 a_m a_{2m} \dots a_{m(n-1)} = 1$$

*follows from finitely many relations $Q_\lambda \equiv Q(\dots, a_{K+\lambda}, \dots) = 1$ (with $\lambda = 0, \pm 1, \dots$). By the strong form of the Freiheitssatz we have only to consider those $Q_\lambda$ in which the $a_K$ with $0 \le K \le m(n-1)$ appear.* If we assume (without loss of generality) that in $Q_0$ at most $a_0, a_1, \dots, a_s$ appear (and $a_0$ and $a_s$ really do appear), then the relation

(1)                 $a_0 a_m \cdots a_{m(n-1)} = 1$    follows from

$$
\begin{array}{l}
Q(a_0, a_1, \dots, a_s) = 1 \\
Q(a_1, a_2, \dots, a_{s+1}) = 1 \\
\cdots\cdots\cdots\cdots\cdots\cdots\cdots \\
\cdots\cdots\cdots\cdots\cdots\cdots\cdots \\
Q(a_{m(n-1)-s}, \dots, a_{m(n-1)}) = 1
\end{array}
\left.\right\}
\quad
\begin{array}{l}
m(n-1) + 1 - s \quad \text{relations} \\
\text{over } m(n-1) + 1 \quad \text{generators.}
\end{array}
$$

(2)

The existence of an identity

(3)                 $$a_0 a_m \cdots a_{m(n-1)} \equiv \prod_{i=1}^{h} T_i Q_{l_i}^{e_{l,i}} T_i^{-1}$$

where $e_{l,i} = \pm 1$ and $l_i$ is a number in the sequence $0, 1, \dots, m(n-1)-s$, has already been used; and to make statements about the exponent sums of $a_0, \dots, a_s$ in $Q_0$, we set: $\sum_i e_{l,i} = e_l$ (where $e_l$ is, so to speak, the exponent sum of $Q_l$ in the right-hand side of (3)); if we now denote the exponent sum of $a_K$ in $Q_0$ (where $K = 0, 1, \dots, s$) by $d_k$ (which is then equal to the exponent sum of $a_{K+l}$ in $Q_l$) and set $c_i = \begin{cases} 1 & \text{for } m \mid i \\ 0 & \text{otherwise} \end{cases}$ (where $c_i$ is the exponent sum of $a_i$ (for $0 \le i \le m(n-1)$) in $a_0 a_m \cdots a_{m(n-1)}$), then from (3) by abelianizing it follows that:

(4)                 $$c_i = \sum_{0 \le i-l \le s} e_l d_{i-l}$$

for $i = 0, 1, \dots, m(n-1)-s$. The $c_i$ are known. *We want to find all integer solutions $e_l, d_k$ to (4)*; we can do this easily using the following remark: the equations (4) are necessary and sufficient conditions for the following to be identically true in the variable $z$:

(5)                 $$\sum_{i=0}^{m(n-1)} c_i z^i = \left\{ \sum_{l=0}^{m(n-1)-s} e_l z^l \right\} \left\{ \sum_{k=0}^{s} d_k z^k \right\}.$$

The left-hand side is the polynomial

$$1 + z^m + z^{2m} + \cdots + z^{m(n-1)} = \frac{z^{mn} - 1}{z^m - 1}.$$

We thus obtain restrictions both for the number $s$ as well as the exponent sums of $a_0, \ldots, a_s$ in $Q_0$ by the condition: $\sum_{k=0}^{s} d_k z^k$ *is a divisor of* $\frac{z^{mn}-1}{z^m-1}$. Since the $d_k$ are integers, there are only finitely many integral solutions to (4). Unfortunately, it is not possible to prove that for every solution of (4) there are only a finite number of systems of relations (2) from which (1) follows. On the other hand, it may well be that for a given solution to (4) (i.e. a system (2) that is "possible in the abelianization") there is no system of relations (2) from which (1) follows.

The main tool of this section will be the following result:

3. **Main Lemma: If from the system** (2) **(in Nr. 2) the relation** (1) **follows, then** $a_0$ **and** $a_s$ **only occur once in** $Q_{(a_0, \ldots, a_s)}$. (Note that $Q$, written cyclically, admits no absorptions). To prove this, the following remark suffices: if the relation (1) follows from (2), then from the system $(\overline{2})$ on the generators $a_0, a_1, \ldots, a_{mn}$:

$$(\overline{2}) \qquad \begin{cases} Q(a_0, \ldots, a_s) = 1 \\ \cdots\cdots\cdots\cdots\cdots \\ \cdots\cdots\cdots\cdots\cdots \\ Q(a_{mn-s}, \ldots, a_{mn}) = 1 \end{cases}$$

the relation

$$(\overline{1}) \qquad\qquad a_0 a_{mn}^{-1} = 1$$

follows. This is because from (2) follows not only $a_0 a_m \cdots a_{m(n-1)} = 1$ but also $a_m a_{2m} \cdots a_{mn} = 1$, and from this we hence have $(\overline{1})$. Now we have to show: if the relation $(\overline{1})$ follows from $(\overline{2})$, then $a_0$ and $a_s$ appear only once in $Q_0 \equiv Q(a_0, \ldots, a_s)$. To prove this, we use the lemmas from §5.

First of all, from Lemma 4, in the case that $s > 0$ (if $s = 0$ all is trivial), it follows: all relations for $a_0$ and $a_{mn}$ which follow from $(\overline{2})$, are obtained from $Q_{mn-s} \equiv Q(a_{mn-s}, \ldots, a_{mn}) = 1$ and the relations for $a_0$ and $a_{mn-s}, \ldots, a_{mn-1}$ obtainable from $Q_0 \equiv Q(a_0, \ldots, a_s) = 1, \ldots, Q_{mn-s-1} \equiv Q(a_{mn-s-1}, \ldots, a_{mn-1}) = 1$, and can be denoted by $S_r(a_0, a_{mn-s}, \ldots, a_{mn-1}) = 1$ (where $r = 0, 1, 2, \ldots$).[24] Each relation $S_r = 1$ actually contains the letter $a_0$.

---

[24]*Proof:* We identify

| in Lemma 4, §5 | with | in the current case |
|---|---|---|
| $a$ | " | $a_0$ |
| $b$ | " | $a_{mn}$ |
| $x$ | " | $a_{mn-s}, \ldots, a_{mn-1}$ |
| $y$ | " | $a_1, \ldots, a_{mn-s-1}$ |
| the system (A$_1$) | " | $Q_0 = 1, \ldots, Q_{mn-s-1} = 1$ |
| the system (B$_1$) | " | $Q_{mn-s} = 1.$ |

That the assumptions of Lemma 4 hold follows easily from the strong form of the *Freiheitssatz*.

We now use Lemma 5 from §5, by identifying:

| in Lemma 5, §5 | with | in the current case |
|---|---|---|
| $a$ | ” | $a_0$ |
| $b$ | ” | $a_{mn}$ |
| $R(a, b)$ | ” | $a_0 a_{mn}^{-1}$ |
| $l$ | ” | $a_{mn-s}, \dots, a_{mn-1}$ |
| the system $(\alpha)$ | ” | the system $S_r = 1$ |
| the system $(\beta)$ | ” | the relation $Q_{mn-s} = 1$. |

Then by (III) in Lemma 5 we have: there is some relation in the system $S_r = 1$, say, $S_{r_0} = 1$, such that for some suitable choice of signs:

$$(a_0 a_{mn}^{-1})^{\pm 1} \equiv T_1 S_{r_0}^{\pm 1}(a_0, a_{mn-s}, \dots, a_{mn-1}) T_1^{-1} T_2 \overline{Q}_{mn-s}^{\pm 1} T_2^{-1},$$

where $\overline{Q}_{mn-s}(a_{mn-s}, \dots, a_{mn}) = 1$ follows from $Q_{mn-s} = 1$. By abelianizing it follows that $a_0$ in $S_{r_0}$, and $a_{mn}$ in $\overline{Q}_{mn-s}$, must have exponent sum $+1$ or $-1$. By using (II) of Lemma 5, it now follows: either $a_0$ appears in $S_{r_0}$ only once, or $a_{mn}$ appears in $\overline{Q}_{mn-s}$ only once. By a direct consideration of absorptions, it is easily seen that from one the other follows. Hence, $a_{mn}$ appears only once in $\overline{Q}_{mn-s}$; hence $\overline{Q}_{mn-s}$ is a primitive element. Since $Q_{mn-s}$ is a root of $\overline{Q}_{mn-s}$, it is consequently a conjugate of $\overline{Q}_{mn-s}$. From this, it follows that $a_{mn}$ only appears once in $Q_{mn-s}$. In exactly the same way, we see that $a_0$ only appears once in $Q_0$. Because of the "isomorphy" of the relations (2), it hence follows: in each relation $Q_t = 1$ (where $t = 0, \dots, mn - s$) the letters $a_t$ and $a_{t+s}$ appear exactly once. □

4. *Applications.* To use the aforementioned remark, we will use the following notation: if one is searching for the relation systems (2) from which (1) follow, then the exponent sums $d_k$ of the $a_k$ in $Q_0$ will satisfy the conditions (4); we express this by saying that: $Q_0$ is $= a_0^{d_0} a_1^{d_1} \cdots a_s^{d_s}$ in the abelianization.

4.a) *First example.* To find the roots of $\overline{a}^2 \overline{b}^p$ (where $p$ is a prime), we proceed as follows: we set $\overline{a} = ab^{+p}$, $\overline{b} = b^{-2}$. The roots of $ab^p ab^p b^{-2p} = ab^p ab^{-p}$ are found, once we set $b^k ab^{-k} = a_k$, and we seek the system

$$\left. \begin{array}{c} Q(a_0, \dots, a_s) = 1 \\ \dots\dots\dots\dots\dots \\ Q(a_{p-s}, \dots, a_p) = 1 \end{array} \right\}$$

from which $a_0 a_p = 1$ follows. According to the divisors with integer coefficients of $z^p + 1$, the word $Q(a_0, \dots, a_s)$ is in the abelianization $= a_0^{\pm 1}$ or $= (a_0 a_1)^{\pm 1}$ or $= (a_0 a_1^{-1} a_2 \cdots a_{p-1}^{\pm 1}$ or $= (a_0 a_p)^{\pm 1}$.

Since $a_0$ and $a_s$ only occur once in $Q(a_0, \dots, a_s)$, we thus have that either $Q_0 \equiv a_0^{\pm 1}$ or $Q_0 \equiv (a_0 a_1)^{\pm 1}$ resp. $(a_1 a_0)^{\pm 1}$, or else $Q_0^{\pm 1}$ is a conjugate of

$$W_0 \equiv a_0 H_1(a_1, \dots, a_{p-2}) a_{p-1} H_2(a_1, \dots, a_{p-2}).$$

Hence, from $W_0 = 1$ and from

$$W_1 \equiv a_1 H_1(a_2, \dots, a_{p-1}) a_p H_2(a_2, \dots, a_{p-1}) \equiv 1$$

the relation $a_0 a_p = 1$ follows, and so we must have

$$H_1(a_1, \dots) a_{p-1} H_2(a_1, \dots) \equiv H_1^{-1}(a_2, \dots) a_1^{-1} H_2^{-1}(a_2, \dots).$$

But we can now easily, by considering possible absorptions, establish that this is not possible. In this case, there is hence no possible $Q$ associated to a solution to (4). If instead $Q_0$ is in the abelianization $= (a_0 a_p)^{\pm 1}$, then we consider $a_0 a_p$ as a new generator (for example together with $a_1$ to $a_p$), and find $Q_0 \equiv (a_0 a_p)^{\pm 1}$ resp. $\equiv (a_p a_0)^{\pm 1}$.

Hence: the roots of $a b^p a b^{-p}$ are given as follows: $a = 1$; or $abab^{-1} = 1$; or $ab^p ab^{-p} = 1$. From $1 = abab^{-1}$ and $ab^p ab^{-p} = 1$, for $p \neq 2$ we have the roots $(\overline{a}\overline{b}^{\frac{p-1}{2}})^2 \overline{b} = 1$ resp. $\overline{a}^2 \overline{b}^p = 1$ of $\overline{a}^2 \overline{b}^p$; however, to $a = 1$ belong no roots of $\overline{a}^2 \overline{b}^p$. For $p = 2$ we have from $abab^{-1} = 1$ no roots of $\overline{a}^2 \overline{b}^p$, but from $a = 1$ resp. $ab^2 ab^{-2} = 1$ we have the roots $\overline{a}\overline{b} = 1$ resp. $\overline{a}^2 \overline{b}^2 = 1$.

4.b) *Similar examples.* By very similar methods we can find all pairwise non-conjugate roots of $a^2 b^{2^k}$, $a^p b^{p^k}$ and similar words.

5. *Simple, unsolved problems.* Even the simple question of finding the roots of $ab^6 a^{-1} b^{-6}$ which have exponent sum zero in $b$, is not simple to answer: indeed, we should have

$$a_0 a_6^{-1} = 1 \text{ following from } \begin{cases} Q(a_0, \dots, a_4) = 1 \\ Q(a_1, \dots, a_5) = 1 \\ Q(a_2, \dots, a_6) = 1. \end{cases}$$

(And $Q_0$ is $= a_0 a_2 a_4$ in the abelianization), and so we first of all – in contrast to the previous example – from *one* solution to (4) get *two* non-conjugate and non-inverse solutions $Q_0 \equiv a_0 a_2 a_4$ and $Q_0 \equiv a_2 a_0 a_4$, and second, using the previous method we have not been able to prove that there is not an infinite number of other such solutions. We must here, for example, be able to prove that $a_2$ appears exactly once in $Q_0$.

# Appendix:

## §8. Remarks on groups with two defining relations

Let $a, b, c, \dots$ be generators, and $R_1(a, b, c, \dots) = 1$, $R_2(a, b, c, \dots) = 1$ be the defining relations of a group. We can – as in the case of a one-relator group – as questions about the words $W(a, b, c, \dots)$ which are equal to 1 on the basis of $R_1 = 1, R_2 = 1$ (the identity problem), or ask to find all relation pairs $R_1 = 1, R_2 = 1$ on the basis of which a given word $W$ is equal to 1 (the root problem). The solutions of these problems only provide something new if they cannot be reduced to the question: to find those (individual) relations $R = 1$ on the basis of which $W$ becomes equal to 1. (This is the root problem associated to one-relator groups.) Such a reduction of the aforementioned problem is possible, for example, when $W = 1$ follows from the *pair* of relations (1) $\{R_1 = 1, R_2 = 1\}$ by virtue of the fact that $W = 1$ already follows from $R_1 = 1$ (or $R_2 = 1$) alone; in general, such a reduction is possible when the pair of relations $R_1 = 1, R_2 = 1$ is, in a sense to be specified, "equivalent" to a pair of relations (1′) $\{R_1' = 1, R_2' = 1\}$, such that $W = 1$ follows from either $R_1' = 1$ or $R_2' = 1$ alone. The relation pairs (1) and (1′) will be said to be equivalent, if the existence of (1) implies the existence of (1′), and vice versa.

By analogy with the terminology of the theory of integral algebraic numbers, a pair of relations (1), when it is not equivalent to a pair (1′) from which $W = 1$ follows from $R'_1 = 1$ alone, will be called an "*ideal root*" of $W$. Before we introduce the questions related to finding ideal roots, it may be appropriate to make a remark about the question of when two pairs of relations are equivalent, since this question, by the above discussion, is important if one wants to decide whether a given pair of relations (1) is an *ideal* root of $W$. The corresponding question which arose in our treatment of one-relator groups (to decide when two relations are equivalent) was solved in §6 of the second part solved using the *Freiheitssatz*. The "natural" generalization of the result thereby obtained (equivalent relations are obtained from one another by conjugacy) would thus be the following theorem: equivalent pairs of relations are obtained from one another by the following process: if $R_1 = 1$ and $R_2 = 1$ are the original relations, then we form – by analogy with the formation of connected primitive elements from the original generators[25] – the expressions:

$$R'_1 \equiv T_1 R_1^{\pm 1} T_1^{-1} T_2 R_2^{\pm 1} T_2^{-1}$$
$$R'_2 \equiv T R_2 T^{-1} \quad (\text{resp. } R'_2 \equiv T R_1 T^{-1})$$

The pair of relations $R'_1 = 1, R'_2 = 1$ is then said to be equivalent with the pair $R_1 = 1, R_2 = 1$, and the presumed theorem would say that, all pairs $R'_1 = 1, R'_2 = 1$ obtained from $R_1 = 1, R_2 = 1$ by repeating this procedure, gives all pairs of relations equivalent to $R_1 = 1, R_2 = 1$. Currently, there is no known way to prove this. The connection between the question of equivalent pairs of relations and the question of the ideal roots of a word $W$ emerges when we ask for the ideal roots of the "easiest" word, the generators, i.e. when we set $W \equiv a$. Indeed, here we have the conjecture that $a$ does not possess any ideal roots, i.e. that every pair of relations $R_1 = 1, R_2 = 1$ from which $a = 1$ follows, is equivalent to a pair: $a = 1, R = 1$, one of whose relations is $a = 1$. This conjecture is thus a generalization of the theorem proved in the first part, which says that $P = T a^{\pm 1} T^{-1}$ whenever $a = 1$ follows from $P = 1$. Here it is notable, that some powers of $a$ possess ideal roots; for $a^3$ we have the following example: from (2) $\begin{cases} b^2 = 1 \\ bab^{-1} = a^{+2} \end{cases}$ follows $a^3 = 1$; indeed, we have

$$b^2 a b^{-2} = b(bab^{-1})b^{-1} = ba^{+2}b^{-1} = (bab^{-1})^{+2} = a^4 = b^2 ab^{-2} = a.$$

The group defined by (2) is the symmetric group $S_{3!}$ of permutations on three symbols.[26] It should now be shown that there is no presentation of $S_{3!}$ equivalent to (2) can be obtained from a pair of relations (2′) $\{a^3 = 1, Q(a,b) = 1\}$. A presentation from $\{a = 1, \overline{Q}(a,b) = 1\}$ is out of the question, since $a$ is not equal to 1). If we denote by $\alpha$ resp. $\beta$ the exponent sum of $a$ resp. $b$ in $Q$, then we must have, if (2′) is to be equivalent to (2), by a frequently used way have to solve following equations in integers $\lambda_1, \ldots$:

$$\lambda_1 + \alpha \lambda_2 = 0 \qquad 3\mu_1 + \alpha \mu_2 = -1$$
$$\beta \lambda_2 = 2 \qquad \beta \mu_2 = 0.$$

---

By $\lambda_2 \beta = 2$ we have $\lambda_2 \neq 0$ and $\alpha \equiv 0 \pmod 3$, and this is a contradiction to $3\mu_1 + \alpha\mu_2 = -1$. In this example, it should be noted that (2) defines a finite, non-cyclic group. Since in finite groups a power of each generator is always equal to 1, it is logical to look for ideal roots of powers of generators among pairs of relations that yield finite, non-cyclic groups. In fact, for example, the following pairs of relations for any integer $n$ and for $m = 2, 3, 4, 5$: (3) $a^2 = b^3$; $(ab)^m b^{3n} = 1$ always yields a finite non-cyclic group; the previous presentation (3) arises from the investigations of Herr Dehn on three-dimensional manifolds that arise from knots. Indeed, depending on whether $m = 2, 3, 4$ or 5, the element $a^2 = b^3$ generates a normal cyclic subgroup of order $|3n + 5|$ resp. $|4n + 10|$ resp. $|6n + 20|$ resp. $|12n + 50|$. The factor group of this normal subgroup is then respectively the group $S_{3!}$, the tetrahedral, octahedral, resp. icosahedral group. In fact, the pair of relations (3) for the different words over $n$ and $m$ yields ideal roots for different powers of $b$ or $a$; however, using the pair of relations (3), nor anywhere else, I have not been able to find ideal roots of $a^2$, so that at least the possibility remains that other than $a$ itself, $a^n$ has no ideal roots for some exponent $n > 1$.

Finally, there is one more problem, which we could equally well have stated at the beginning of this section: the question of when a group with generators $a, b, c, \ldots$ and relations $R_1 = 1, R_2 = 1$ is *essentially two-relator*, i.e. whether the group is not isomorphic to any other group which can be defined by one (or no) relations. This occurs, for example, when $R_1$ is a root of $R_2$; – and to decide whether this is the case or not, is in general unsolved.

# Investigations on some infinite discontinuous groups.

by
Wilhelm Magnus in Frankfurt am Main.

### Contents

## §1.
## Introduction and statement of the results.

In what follows we will apply methods which were introduced by M. Dehn in his treatment of fundamental problems in the theory of infinite discontinuous groups. First, we will in §3 prove a theorem on the automorphism group of Listing's knot, namely that *all* automorphisms of this group arise as automorphisms of the kind given by Dehn in his work "On the two trefoil knots"[1]. The proof will be carried out by assigning the automorphisms to transformations of a binary quadratic form into itself[2] via substitutions with determinant $\pm 1$ and integer coefficients.

Second, we will in §4 solve the identity problem for all groups with two generators $a, b$ and the defining relation $a^{\alpha_1} b^{\beta_1} a^{\alpha_2} b^{\beta_2} = 1$; that is, we will give a procedure which decides, in a finite number of steps, whether a given expression in the generators is equal to the identity element or not. We will, as part of this, make essential use of the theorem by O. Schreier[3] on the "existence of free products with amalgamated subgroups" of two groups. Thus far, the identity problem has only been solved (other than in special cases) when the defining relation is a "binomial" $a^\alpha b^\beta$ equal to the identity.[4]

Third, we will in §5 give a short proof for a theorem (based on unpublished results by M. Dehn) on the subgroups of the modular group; which we can easily state here: each subgroup of the modular group contains a free normal subgroup, the quotient group of which is a cyclic group of order $1, 2, 3$, or $6$. The proof of this is based on proving in particular: the commutator group of the modular group is a free group

---

[1] M. Dehn, *On the two trefoil knots*, Math. Annalen **75** (1914), p. 412.

[2] Up to a factor of $\pm 1$.

[3] O. Schreier, *The subgroups of free groups.* Abhandlungen aus dem Mathematischen Seminar der Hamburgischen Universität **5**, Leipzig 1927, p. 161 ff.

[4] See H. Gieseking, *Analytic investigation on topological groups*, Ph. D. thesis, Münster 1912.

on two generators (from which it follows, incidentally, that all free groups can be realized, in an easily specified manner, as subgroups of the modular group).

The methods required for the proofs are explained in detailed in an article on "discontinuous groups with one defining relation"[5], but nevertheless we will briefly go through these tools in §2 (without proofs). The main result of the cited article, being the so-called *Freiheitssatz*, is perhaps mostly dispensable in the following, since this result, being also rather difficult to prove, is not used here in its full generality, but in any case the methods devised for its proof will be used extensively here.

Sections 3, 4, and 5 can be read independently of one another.

# §2.
# The method.

**Theorem 1**. *Let $G$ be any group given via generators $a, b, c, \dots$ and some relations holding between them:*

$$R_i(a, b, c, \dots) = 1 \qquad\qquad (i = 1, 2, \dots)$$

*Let $W(a, b, c, \dots)$ be an expression in the generators (a "word") which is equal to 1 on the basis of the defining relations $R_i = 1$. Then $W$ can be transformed by using insertions and deletions of $aa^{-1}, bb^{-1}, \dots$ into an expression of the form*

$$T_1 R_{i_1}^{\varepsilon_1} T_1^{-1} T_2 R_{i_2}^{\varepsilon_2} T_2^{-1} \cdots T_k R_{i_k}^{\varepsilon_k} T_k^{-1}$$

*where the $\varepsilon_1, \dots, \varepsilon_k$ are $\pm 1$, and $i_1$ to $i_k$ are arbitrary natural numbers. The $T$ are arbitrary expressions in the generators. We write this in abbreviated form as:*

$$W \equiv \prod_{\lambda=1}^{k} T_\lambda R_{i_\lambda}^{\varepsilon_\lambda} T_\lambda^{-1} \qquad\qquad (\equiv \text{ is called "identical".})$$

**Theorem 2**. *If in the group $G\{a, b, c, \dots; R_i(a, b, c, \dots) = 1\}$ the word $W(a, b, c, \dots)$ is equal to 1, then $W$ is (furthermore) equal to 1 in the group $\overline{G}$ with the generators $a, b, c, \dots$ and the relations $R_i = 1$ of $G$ together with some new relations $\overline{R}_k = 1$ (the group $\overline{G}$ is a quotient group of $G$).*

In particular, in what follows the addition of $ab = ba, ac = ca, bc = cb, \dots$ as relations, which we call the *abelianization* of the group in question, will play an important role. When doing this, we will use the concept of the *exponent sum*; if $W$ is a word in the generators $a, b, c, \dots$ and $W$ is written as $a^{\alpha_1} b^{\beta_1} c^{\gamma_1} \cdots a^{\alpha_2} b^{\beta_2} c^{\gamma_2} \cdots a^{\alpha_n} b^{\beta_n} c^{\gamma_n}$, then $\alpha = \alpha_1 + \alpha_2 + \cdots \alpha_n$ is called the *exponent sum* of $a$ in $W$ (of course, some or all of the exponents $\alpha_1, \dots, \alpha_n$ can be zero). Exponent sums are invariants under identical transformations. It then follows immediately by abelianizing: if $W$ is identically 1 (that is, it is equal to 1 in the free group generated by $a, b, c, \dots$), then $W$ has exponent sum zero in all generators $a, b, c, \dots$. This results in the following:

**1. Application of Theorem 2**. Let $G$ be a group with two generators $a, b$ and one defining relation $R(a, b) = 1$ (abbreviated $G\{a, b, R(a, b) = 1\}$). Suppose $R$ has exponent sum $\varrho_1$ resp. $\varrho_2$ in the generator $a$ resp. $b$. Suppose the word $W(a, b)$ over $a$ and $b$ has exponent sum zero $\alpha$ resp. $\beta$ in $a$ resp. $b$, and that it is equal to 1 on the

---

basis of $R = 1$. By Theorem 1 we identically have

$$W \equiv \prod_{\lambda=1}^{k} T_\lambda R_{i_\lambda}^{\varepsilon_\lambda} T_\lambda^{-1}. \qquad\qquad (\varepsilon_\lambda = \pm 1)$$

By abelianizing, it follows that

$$\alpha = \varrho_1 \cdot \sum_{\lambda=1}^{k} \varepsilon_\lambda, \quad \beta = \varrho_2 \cdot \sum_{\lambda=1}^{k} \varepsilon_\lambda$$

and consequently, in the case that $\alpha = 0$ and $\varrho_1 \neq 0$, that

$$\sum_{\lambda=1}^{k} \varepsilon_\lambda = 0.$$

In this case, $W = 1$ follows from the relations $R \cdot a = a \cdot R$ and $R \cdot b = b \cdot R$. Indeed,

$$W \equiv (T_1 R^{\varepsilon_1} T_1^{-1} R^{-\varepsilon_1}) R^{\varepsilon_1} (T_2 R^{\varepsilon_2} T_2^{-1} R^{-\varepsilon_2}) R^{-\varepsilon_1} \cdot R^{\varepsilon_1 + \varepsilon_2} (T_3 \cdots R^{-\varepsilon_3}) \cdots$$
$$\cdots R^{\varepsilon_1 + \varepsilon_2 + \cdots + \varepsilon_{k-1}} (T_k R^{\varepsilon_k} T_k^{-1} R^{-\varepsilon_k}) \cdot R^{-\varepsilon_1 - \cdots - \varepsilon_{k-1}} \cdot R^{\varepsilon_1 + \cdots + \varepsilon_k}.$$

By assumption, however, $\varepsilon_1 + \cdots + \varepsilon_k = 0$, so $R^{\varepsilon_1 + \cdots + \varepsilon_k} = 1$, and from $Ra = aR$ and $Rb = bR$ it follows that $T_\lambda R^{\varepsilon_\lambda} = R^{\varepsilon_\lambda} T_\lambda$ (for $\lambda = 1, \ldots, k$), and hence also

$$R^{\varepsilon_1 + \cdots + \varepsilon_{\lambda-1}} \cdot (T_\lambda R^{\varepsilon_\lambda} T_\lambda^{-1} R^{-\varepsilon_\lambda}) \cdot R^{-(\varepsilon_1 + \cdots + \varepsilon_{\lambda-1})} = 1$$
$$\text{for } \lambda = 2, 3, \ldots, k.$$

Generalizing this to groups with more generators and relations is easily done; for later (§5) we will make a remark: if in the modular group, which is defined by the generators $a, b$ and the relations $a^2 = b^3 = 1$, a word $W(a, b)$, which has exponent sum zero in $a$ and in $b$, is equal to 1, then $W$ is also equal to 1 on the basis of the relations $a^2 b a^{-2} b^{-1} = 1$ and $b^3 a b^{-3} a^{-1} = 1$.

For §4 the following:

## 2. Application of Theorem 2 is important. If in a group

$$G(a, b, c, \ldots; R_i(a, b, c, \ldots) = 1)$$

all $R_i$ have exponent sum zero in $a$, then a word $W(a, b, c, \ldots)$ which has a non-zero exponent sum zero in $a$ cannot be equal to 1 in $G$. To solve the identity problem for such a group $G$ it is therefore enough to solve the identity problem in the subgroup $H_a$ of $G$ consisting of all words which have exponent sum zero in $a$. The following two theorems are concerned with subgroups of such a type $H_a$:

**Theorem 3**. *Let $G$ be the free group on the generators $a, b, c, \ldots$. We set*

$$a^k b a^{-k} = b_k, \quad a^k c a^{-k} = c_k, \ldots \qquad \text{for all } k = 0, \pm 1, \ldots.$$

*Then we have: the (normal) subgroup $H_a$ of $G$, consisting of all words which have exponent sum zero in $a$, is generated by the $b_k, c_k, \ldots$ $(k = 0, \pm 1, \pm 2, \ldots)$, and furthermore is free on these generators.*

(In particular, the free group on two generators contains a normal subgroup on infinitely many generators). Combining Theorem 1 and Theorem 3 now yields:

**Theorem 4.** [6] *Suppose that in $G\{a, b, c, \dots; R_i(a, b, c, \dots) = 1\}$ all $R_i$ have exponent sum zero in $a$. Consider the subgroup $H_a$ consisting of all words which have exponent sum zero in $a$. Its generators are, by Theorem 3, given by*

$$a^k b a^{-k} = b_k, \quad a^k c a^{-k} = c_k, \quad \dots$$

*By assumption, we can write the $R_i$ over the $b_k, c_k, \dots$. Let the result of this be*

$$R_i(a, b, c, \dots) = \overline{R}_i(b_{k_1}, b_{k_2}, \dots, b_{k_r}; c_{l_1}, \dots, c_{l_s}; \dots).$$

*Then all relations of $H_a$ are given by*

$$\overline{R}_i(b_{k_1+\lambda}, \dots, b_{k_r+\lambda}; c_{l_1+\lambda}, \dots, c_{l_s+\lambda}; \dots) = 1,$$

*where $\lambda$ runs across all (positive and negative) integers (including zero).*

Further theorems, and above all a theorem by O. Schreier and the *Freiheitssatz*, will be formulated in §4, in which they will also all be used.

# §3.
# The automorphisms of the group of Listing's knot.

The group of Listing's knot (Fig. 1) is defined[7] by the generators $a_1, a_2, a_3, a_4, a_5$ and the relations

$$a_3 a_4^{-1} a_1 = a_1 a_2^{-1} a_3 = a_1 a_4^{-1} a_5 a_2^{-1} = a_3 a_2^{-1} a_5 a_4^{-1} = 1.$$

The automorphisms of this group include, in addition to the inner automorphisms, the following:

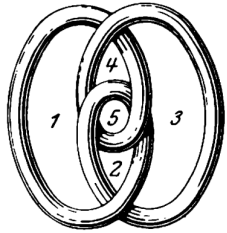

Fig. 1.
Listing's knot

$$\left.\begin{aligned} a_1' &= a_3^{-1} \\ a_2' &= a_2^{-1} \\ a_3' &= a_1^{-1} \\ a_4' &= a_4^{-1} \\ a_5' &= a_5^{-1} \end{aligned}\right\} \text{ called } \overline{j}_0, \text{ and } \left\{\begin{aligned} a_1' &= a_2 a_5^{-1} \\ a_2' &= a_3 a_5^{-1} \\ a_3' &= a_4 a_5^{-1} \\ a_4' &= a_1 a_5^{-1} \\ a_5' &= a_5^{-1} \end{aligned}\right.$$

which we call $\overline{j}_1$. Then $\overline{j}_0^2, \overline{j}_1^4$ and $(\overline{j}_0 \overline{j}_1)^2$ are inner automorphisms; the normal subgroup, of the group of all automorphisms, generated by the inner automorphisms thus admits the dihedral group of order eight as a quotient group. Indeed, we have the following:

*All automorphisms of $G$ are generated by $\overline{j}_0, \overline{j}_1$ and inner automorphisms.*

---

We will alter $G$ by first re-expressing all generators as words over $a_1$ and $a_3$, and introduce new generators $u, v$ by

$$a_1 a_3 a_1^{-2} = u, \quad a_1 = vu,$$
$$a_1^3 a_3^{-1} a_1^{-1} = v, \quad a_3 = (vu)^{-1} u (vu)^2.$$

This presents $G$ as

$$G\{u, v; u^3 = v^{-1} u v^2 u v^{-1}\},$$

and the automorphisms $\bar{j}_0 \bar{j}_1^2$ resp. $\bar{j}_1$ are given[8] by

$$j_0: \quad \{u' = u^{-1}; \; v' = v^{-1}\}$$

and

$$j_1: \quad \{u' = v^{-1} u^{-1} v u; \; v' = u^{-1} v^{-1} u^{-2} v^{-1} u v\}.$$

We now consider the subgroup $H$ of $G$, which consists of all words which have exponent sum zero in $v$. The subgroup $H$ is normal in $G$, and, as we shall shortly prove, is even characteristic in $G$ (i.e. $H$ is mapped to itself by every automorphism of $G$). Indeed, this follows from a property of our group that we will make frequent use of in what follows:

*If $J$ is an automorphism of $G$ defined by $u' = Q(u, v)$ and $v' = S(u, v)$, then $Q$ has exponent sum zero in $v$, and $S$ has exponent sum $\pm 1$ in $v$.*

*Proof.* Since the commutator subgroup of $G$ is a characteristic subgroup, $J$ is also an automorphism of this commutator subgroup and hence also by its quotient group in $G$, which we call $\overline{G}$. Then $\overline{G}$ is defined by $u, v; uvu^{-1}v^{-1} = 1$, $u^3 = v^{-1}uv^2uv^{-1}$ or $u, v; u = 1$, which thus possesses only the automorphisms $u' = 1$, $v' = v^{\pm 1}$. Since $\overline{G}$ is obtained from $G$ by the addition of the relation $u = 1$, we must hence have $Q(1, v) = 1$ and $S(1, v) = v^{\pm 1}$. □

Consequently, words in $u$ and $v$ which have exponent sum zero in $v$ are again mapped to such words upon replacing $u$ and $v$ by $Q(u, v)$ resp. $S(u, v)$. Hence $H$ is mapped into itself by $J$.

If we now denote the automorphisms $J$ for which $v$ in $S$ has exponent sum $+1$ by $J^+$, and the others by $J^-$, then the $J^+$ clearly form a normal subgroup of index 2 in the group of all automorphisms $J$. (We cannot obtain the entire group, as $j_0$ and $j_1$ are automorphisms of the form $J^-$). In the following, we must prove that all automorphisms of the form $J^+$ can be generated by $j_1^2, j_0 j_1$, and inner automorphisms. To do this, we first prove using the group $H$ defined above: *every automorphism $J^+$ of $G$, which "induces"[9] an inner automorphism of $H$, is itself an inner automorphism of $G$.*

*Proof.* We investigate $H$. If we set

$$v^k u v^{-k} = u_k \quad \text{for} \quad k = 0, \pm 1, \pm 2, \dots,$$

then by Theorem 4, §2 we have that $H$ is given by

$$H\{u_k; u_k^3 = u_{k-1} u_{k+1}\},$$

---

and since we can here eliminate all relations and all generators except two, say $u_0$ and $u_1$, it follows that $H$ is the free group on the generators $u_0$ and $u_1$. To consider the action of $J^+$ on $H$, we set $S(u, v) \equiv S'(u, v) \cdot v$; then $S'$ has, just as $Q$, exponent sum zero in $Q$; thus $Q$ and $S'$ belong to $H$, and can thus be expressed in terms of $u_0$ and $u_1$: $Q = \overline{Q}(u_0, u_1)$, $S' = \overline{S}(u_0, u_1)$, and $J^+$ induces in $H$ an automorphism $K^+$, which we write as $u_0' = \overline{Q}(u_0, u_1)$, $u_1' = \overline{S}(u_0, u_1)\overline{Q}(u_1, u_0^{-1}u_1^{+3})\overline{S}^{-1}$ (since $u_0^{-1}u_1^{+3} = u_2$), and when $K^+$ is an inner automorphism – e.g. the identity – then $Q \equiv u$, and $\overline{S} \equiv u_1^\ell$ (where $\ell$ is an integer), and hence $S = v \cdot u^\ell$, and since $u'^3 = v'^{-1}u'v'^2u'v'^{-1}$ or $u^3 = u^{-\ell}v^{-1}uvu^\ell vuv^{-1}$ or $u_0^3 = u_0^{-\ell}u_{-1}u_0^\ell u_1$ or, since by $u_{-1} = u_0^3 u_1^{-1}$ we must have $u_0^{-\ell}u_0^3 u_1^{-1}u_0^\ell u_1$, we have $\ell = 0$.                               $\square$

Consider now quotient group $\overline{H}$ corresponding to the commutator subgroup of $H$; then $\overline{H}$ is an abelian group on the generators $u_0, u_1$, being defined by

$$\overline{H}\{u_k; u_k^3 = u_{k-1}u_{k+1},\ u_0u_1 = u_1u_0 \quad (k = 0, \pm 1, \dots)\}.$$

Then $\overline{H}$ is also characteristic for $G$; each automorphism $J^+$ induces an automorphism $K^+$ of $\overline{H}$. We then have[10]:

*If $J^+$ induces the identity automorphism on $\overline{H}$, then it induces an inner automorphism on $H$, because if $J^+$ induces $K^+$ on $H$, and $K^+$ induces the identity automorphism on $\overline{H}$, then $K^+$ is an inner automorphism of $H$.*

We now say: if an automorphism $\overline{K}^+$ of $\overline{H}$ is induced by an automorphism $J^+$ of $G$, then we will call this: $\overline{K}^+$ can conversely be "*extended*" to an automorphism $J^+$ of $G$. If we can now prove the following: *each automorphism $\overline{K}^+$ of $\overline{H}$, which can be extended to an automorphism $J^+$ of $G$, is induced by such an automorphism $J^+$ which is generated by $j_0, j_1, j_1^2$, and inner automorphisms*, then we will be done.

An arbitrary automorphism $\overline{K}^+$ of $\overline{H}$ can be written as

$$u_0' = u_0^{s_0}u_1^{s_1}, \quad u_1' = u_0^{\sigma_0}u_1^{\sigma_1}; \quad \begin{vmatrix} s_0 & s_1 \\ \sigma_0 & \sigma_1 \end{vmatrix} = \pm 1.$$

If this is induced by an automorphism $J^+$ of $G$, then it follows from $u_0' = u_0^{s_0}u_1^{s_1}$ and the commutativity of the $u_k$, that $u_1' = u_1^{s_0}u_2^{s_1} = u_0^{s_0}u_1^{3s_1}u_0^{-s_1}$ (since $u_2 = u_0^{-1}u_1^3$), and thus $\sigma_0 = -s_1, \sigma_1 = s_0 + 3s_1$, which yields[11] by substitution into the determinant:

$$(1) \qquad\qquad s_0^2 + 3s_0s_1 + s_1^2 = \pm 1.$$

An arbitrary automorphism $\overline{K}^+$ of $\overline{H}$, which is induced from an automorphism $J^+$ of $G$, is uniquely determined by the associated values of $s_0$ and $s_1$, since $\sigma_0 = -s_1$ and $\sigma_1 = s_0 + 3s_1$ are uniquely determined by $s_0$ and $s_1$. On the basis of our previous investigation it hence follows that each automorphism $J^+$ of the entire group $G$ is determined uniquely, up to inner automorphism, by the exponents $s_0$ and $s_1$ which define the automorphism $\overline{K}^+$ of $\overline{H}$ induced by $J^+$.

Thus it is now sufficient to show:

I. *For each solution $s_0, s_1$ of* (1), *there is an automorphism $J^+$ of $G$, which induces an automorphism of $\overline{H}$ defined by $u_0' = u_0^{s_0}u_1^{s_1}$ and $u_1' = u_0^{\sigma_0}u_1^{\sigma_1}$.*

---

[10]See J. Nielsen, *The automorphisms of general infinite groups with two generators*, Math. Annalen **78** (1918), p. 393.

[11]The reduction of our task to the Diophantine equation (1) – and thus to the presently following number-theoretic considerations – is naturally dependent on our manner of presenting our group by the appropriately chosen generators $u$ and $v$.

II. *The automorphisms $J^+$ from* I *are generated by inner automorphisms and the automorphisms $j_1^2$ and $j_0 j_1$.*

To prove this, we first need a number-theoretic remark on the solutions to (1). Namely, we have[12]:

All integer matrices

$$\begin{pmatrix} s_0 & s_1 \\ \sigma_0 & \sigma_1 \end{pmatrix}$$

with determinant $\pm 1$, which satisfy equation (1) for $s_0, s_1$ (and hence also for $\sigma_0, \sigma_1$), are given by the matrices

$$\pm \begin{pmatrix} -1 & 1 \\ -1 & 2 \end{pmatrix}^{\ell}, \qquad\qquad (\ell = 0, \pm 1, \pm 2, \dots)$$

which are the matrices of all linear binary substitutions with integer coefficients, which transform the binary form $x^2 + 3xy + y^2$ on two variables $x, y$ into itself, up to sign. These matrices obviously form an abelian group with the generators

$$\begin{pmatrix} -1 & 0 \\ 0 & -1 \end{pmatrix} \quad \text{and} \quad \begin{pmatrix} -1 & 1 \\ -1 & 2 \end{pmatrix}.$$

Since now for any two automorphisms $\overline{K_s^+}$ and $\overline{K_t^+}$ of $\overline{H}$, the exponents of which are given by the matrices

$$\begin{pmatrix} s_0 & s_1 \\ \sigma_0 & \sigma_1 \end{pmatrix} \quad \text{and} \quad \begin{pmatrix} t_0 & t_1 \\ \tau_0 & \tau_1 \end{pmatrix},$$

their product $\overline{K_s^+} \cdot \overline{K_t^+}$ has exponents given by the matrix

$$\begin{pmatrix} s_0 & s_1 \\ \sigma_0 & \sigma_1 \end{pmatrix} \cdot \begin{pmatrix} t_0 & t_1 \\ \tau_0 & \tau_1 \end{pmatrix},$$

it follows that proving statements I and II above reduces to proving the following[13]:

I'. For the solutions $s_0 = -1, s_1 = 0$ and $s_0 = -1, s_1 = 1$ of (1) – which express the matrices $\begin{pmatrix} -1 & 0 \\ 0 & -1 \end{pmatrix}$ and $\begin{pmatrix} -1 & 1 \\ -1 & 2 \end{pmatrix}$ – automorphisms $J^+$ of $G$, which induce automorphisms $\overline{K^+}$ of $\overline{H}$, given by $u_0' = u_0^{-1}, u_1' = u_1^{-1}$ and $u_0' = u_0^{-1} u_1, u_1' = u_0^{-1} u_1^2$, respectively.[13]

II'. The automorphisms $J^+$ from I' are generated by inner automorphisms and the automorphisms $j_0 j_1$ and $j_1^2$.

But now we are done, since – up to inner automorphisms – $j_0 j_1$ resp. $j_1^2$ are given by

$$u' = u^{-1} v u v^{-1}; \quad v' = v u^{-1} v^{-1} u v u^2 \quad \text{resp.} \quad u' = u^{-1}; \quad v' = u^2 v u,$$

and these are automorphisms $J^+$, which induce automorphisms $\overline{K^+}$ of $\overline{H}$, given by $u_0' = u_0^{-1}, u_1' = u_1^{-1}$ resp. $u_0' = u_0^{-1} u_1, u_1' = u_0^{-1} u_1^2$.

We will now give an abstract representation of the automorphism group of $G$: if we denote by $j_u$ resp. $j_v$ the inner automorphisms

$$u' = u, \ v' = uvu^{-1} \quad \text{resp.} \quad u' = vuv^{-1}, v' = v,$$

---

[12]See e.g. F. Klein, *Lectures on the theory of elliptic modular functions*, published by R. Fricke, Vol. **2** (Leipzig 1892), p. 161 f. There, only substitutions with determinant $+1$ are treated; this causes a slight deviation from what is stated above.

[13]This is because if these two automorphisms $J^+$ induce automorphisms in $\overline{H}$ given by $s_0 = -1, s_1 = 0$ resp. $s_0 = -1, s_1 = 1$, then by the aforementioned remark we can, for each solution $s_0, s_1$ of (1) construct an automorphism $J^+$ which in $\overline{H}$ induces an automorphism in $\overline{K^+}$ given by $u_0' = u_0^{s_0} u_1^{s_1}, u_1' = \cdots$.

then between these we only have the relation

$$(2) \qquad j_u^3 = j_v^{-1} j_u j_v^2 j_u j_v^{-1} \qquad \text{(corresponding to } u^3 = v^{-1} u v^2 u v^{-1})$$

since the centre of $G$ consists only of the identity element, as can be proved with little effort. The automorphism group of $G$ is hence, as is easily shown, given by the generators $j_0, j_1, j_u, j_v$ and the relation (2) together with

$$(3) \qquad \begin{cases} j_0^2 = 1; \quad (j_0 j_u)^2 = 1; \quad (j_0 j_v)^2 = 1, \\ \qquad (j_0 j_1)^2 = j_u j_v; \quad j_1^4 = 1, \\ j_u j_1 = j_1 \cdot j_u^{-1} j_u^{-1} j_v j_u; \quad j_v j_1 = j_1 \cdot j_u^{-1} j_u^{-1} j_u^{-2} j_v^{-1} j_u j_v. \end{cases}$$

The group of transformations of the quadratic form $x^2 + 3xy + y^2$ into itself by unimodular substitutions arises by setting $j_u$ to be equal to 1 in the subgroup generated by the $j_1^2, j_0 j_1, j_u, j_v$.

Finally, the following is noteworthy: since $H$ is characteristic in $G$, i.e. every automorphism of $G$ maps $H$ into itself, and since $H$ is a free group on the generators $u = u_0$ and $vuv^{-1} = u_1$, the element $u_0 u_1 u_0^{-1} u_1^{-1}$ is necessarily mapped by all automorphisms of $G$ into a conjugate of itself[14]; this fact, used in an essential way above, corresponds geometrically to the fact that $u_0 u_1 u_0^{-1} u_1^{-1}$ is a "longitudinal curve"[15] of the knotted tube belonging to Listing's knot, which is the boundary of the exterior space, and that all such longitudinal curves are conjugate to one another.

## §4.
## The identity problem for the groups $a^{\alpha_1} b^{\beta_1} a^{\alpha_2} b^{\beta_2} = 1$

In this section, we will above all make use of a theorem of O. Schreier[16], which we will presently state; furthermore, we will prove two Lemmas, which allow us to solve the identity problem in certain groups when it is solved for certain other groups. In addition, the *Freiheitssatz*[17] will be used, an appropriate version of which will therefore be formulated here; its use could probably be avoided, although it will be quicker if we make use of it.

We thus begin with the theorem of O. Schreier, which we (somewhat differently) state as follows: let $G$ be a group with two systems of generators $x_1; x_2; \ldots$ and $y_1; y_2; \ldots$ (either system can contain infinitely many generators). Suppose the generators $x_1, x_2, \ldots$ form a group $G_x$ with the defining relations

$$R_i(x_1, x_2, \ldots) = 1 \qquad\qquad (i = 0, 1, 2, \ldots)$$

and similarly the $y_1, y_2, \ldots$ form a group $G_y$ with the defining relations

$$S_k(y_1, y_2, \ldots) = 1. \qquad\qquad (k = 1, 2, \ldots)$$

Suppose further that $G_x$ and $G_y$ have two isomorphic subgroups $H_x$ resp. $H_y$; if we denote by $h_{x,l}$ resp. $h_{y,l}$ $(l = 1, 2, \ldots)$ the elements of $H_x$ resp. $H_y$, the notation expressing a specified isomorphism between $H_x$ and $H_y$, then we suppose that the

---

[14]See J. Nielsen, *loc. cit.* p. 393, equation (11).

[15]For these topological concepts, see M. Dehn, *loc. cit.* p. 412. The longitudinal curves mentioned there are boundaries of the exterior space.

[16]O. Schreier, *loc. cit.*, see Remark 3.

[17]See §1.

only defining relations of $G$ – other than $R_i = 1$ and $S_k = 1$ – will be the relations $h_{x,l} = h_{y,l}$, where $h_{x,l}$ and $h_{y,l}$ ranges across all of $H_x$ resp. $H_y$.

*$G$ is called the free product of $G_x$ and $G_y$ with amalgamated subgroups $H_x$ and $H_y$.[18] We have the following: if one can solve the identity problem in $G_x$ and $G_y$, and if one can decide in $G_x$ and $G_y$ when an expression in the $x_1, x_2, \ldots$ resp. the $y_1, y_2, \ldots$ belongs to $H_x$ resp. $H_y$ – for infinite groups $H_x, H_y$ this is in general more than the solution of the identity problems for $G_x$ and $G_y$ – then one can also solve the identity problem in $G$.[18]*

Indeed: if we write an arbitrary element $W$ of $G$ in the form

$$(1) \qquad W \equiv g_{x,1} g_{y,1} g_{x,2} g_{y,2} \cdots g_{x,n} g_{y,n},$$

where the $g_{x,\nu}$ resp. $g_{y,\nu}$ $(\nu = 1, \ldots, n)$ are elements of $G_x$ resp. $G_y$, i.e. are words over the $x_1, \ldots$ and $y_1, \ldots$, then when $W = 1$ we must have that at least one of the $g_{x,\nu}$ belongs to $H_x$ or at least one of the $g_{y,\nu}$ belongs to $H_y$; if we then replace $g_{x,\nu}$ resp. $g_{y,\nu}$ by its corresponding expression in $H_y$ resp. $H_x$, then we obtain a new representation of the form (1) for $W$, the right-hand side of which has fewer than $2n$ of the expressions $g_{x,\nu}, g_{y,\nu}$, and we can apply the same process to this new representation. The proof that this procedure is correct follows directly from Schreier's formulation of his theorem. $\qquad \square$

Our first application of the theorem will be:

**Lemma 1.** *If in the group $G\colon G(a, b; R(a, b) = 1)$[19] one can solve not only the identity problem, but can also decide (in a finite number of steps) when a word $W(a, b)$ over $a$ and $b$ is equal to a power of $a$, then one can also solve the identity problem in the group*

$$\overline{G}(a, b, t; a = t^n, \ R(a, b) = 1).$$

*Proof.* [20] From Theorem 2, §2 it easily follows that it suffices to solve the identity problem in the subgroup $\overline{H}_t$ of $\overline{G}$ consisting of all the words in which $t$ has exponent sum zero; because if $W'(a, b, t)$ has non-zero exponent sum $\tau$ in $t$, then $W' = 1$ implies: $\tau = \lambda \cdot n$ ($\lambda$ an integer) and $W = W' \cdot a^{+\lambda} t^{-\tau}$ lies in $\overline{H}_\tau$ and is also equal to 1. If one sets $t^k a t^{-k} = a_k$ and $t^k b t^{-k} = b_k$, then by §2, Theorem 4 and an easy generalization of the first usage of §2, Theorem 2, $\overline{H}_t$ is given by

$$\overline{H}_t \{ a_k, b_k; R(a_k, b_k) = 1, \ a_k^{-1} b_{k+n} a_k = b_k, \ a_k = a_{k+1}; \quad k = 0, \pm 1, \pm 2, \ldots \}$$

or, by eliminating the superfluous generators by means of $a_k^{-1} b_{k+n} a_k = b_k, a_k = a_{k+1}$, we have that $\overline{H}_t$ is given by

$$\overline{H}_t \{ a_\nu, b_\nu; R(a_\nu, b_\nu) = 1; \ a_0 = a_1 = \cdots = a_{n-1}; quad \nu = 0, \ldots, n-1 \}.$$

Obviously, we now see that $\overline{H}_t$ is a free product with amalgamated subgroups of groups $G_\nu$ $(\nu = 0, \ldots, n-1)$, all isomorphic to $G$, and defined by

$$G_\nu \{ a_\nu, b_\nu; R(a_\nu, b_\nu) = 1 \}$$

where the subgroups of $G_\nu$ generated by the $a_\nu$ are the amalgamated subgroups. Since $G_\nu$ is isomorphic to $G$, and since by assumption one can decide in $G$ whether a word is equal to a power of $a$, the proof now follows by applying Schreier's theorem. $\qquad \square$

---

[18] It is easy to construct a generalization of this notion to free products with amalgamated subgroups of more than two groups, as well as a corresponding generalization of the associated theorem on the solution of the identity problem.

[19] $a, b$ are generators, and $R(a, b) = 1$ is the defining relation of $G$.

[20] This can be proved in a more immediate and natural way by constructing Dehn's *Gruppenbild* [= Cayley graph] associated to $\overline{G}$.

In the same manner as Lemma 1, we can also prove:

**Lemma 2.** *If one can solve the identity problem in $\overline{G}\{a, b, t; a = t^n, \; R(a, b) = 1\}$, then one can also solve it in $G(a, b; R(a, b) = 1)$, and furthermore a word $W(a, b)$ in $\overline{G}$ is equal to $1$ if and only if it is equal to $1$ in $G$ (thus $R(a, b) = 1$ is the defining relation of the subgroup of $G$ generated by $a$ and $b$).*

*Proof.* We use the same notation as in Lemma 1. If we consider it as a word in $\overline{G}$, the word $W(a, b)$ always lies in $\overline{H}_t$; we can present $\overline{H}_t$ in the same manner as in the proof of Lemma 1 by

$$\overline{H}_t\{a_0; b_0, b_1, \ldots, b_{n-1}; R(a_0, b_\nu) = 1; \quad \nu = 0, 1, \ldots, n-1\}.$$

We now only have to prove that $W(a_0, b_0) = 1$ in $\overline{H}_t$ follows from $R(a_0, b_0) = 1$ alone. This is proved as follows: by §2, Theorem 1, we have, when $W(a_0, b_0) = 1$ in $\overline{H}_t$, an identity

$$(2) \qquad W(a_0, b_0) \equiv \prod_{\lambda=1}^{k} T_\lambda R^{\varepsilon_\lambda}(a_0, b_{s_\lambda}) T_\lambda^{-1},$$

where $\varepsilon_\lambda = \pm 1$ and the $s_\lambda$ are numbers in the sequence $0, 1, \ldots n-1$. Since an identity remains true when one everywhere replaces a letter by another, it follows that if we everywhere in the identity (2) replace $b_\nu$ by $b_0$ for $\nu = 1$ to $n-1$, then the left-hand side of (2) does not change, and (2) becomes

$$(2') \qquad W(a_0, b_0) \equiv \prod_{\lambda=1}^{k} T_\lambda' R^{\varepsilon_\lambda}(a_0, b_0) T_\lambda'^{-1},$$

where the $T_\lambda'$ are obtained from the $T_\lambda$ by everywhere replacing $b_\nu$ ($\nu = 1, \ldots, n-1$) by $b_0$. Hence $W(a_0, b_0) = 1$ follows from $R(a_0, b_0) = 1$ alone. □

The last of the lemmas needed in the sequel will be the *Freiheitssatz*, which we now formulate.

**Freiheitssatz.** *Suppose we have two systems of generators of a group $G$: $a, b, c, \ldots$ is one, and $x_k$ ($k = 0, \pm 1, \pm 2, \ldots$) is the second. If $m$ is a fixed integer $m \geq 0$, then suppose we have the following defining relations between these generators: for $\lambda = 0, \pm 1, \ldots$, we have*

$$R_\lambda\{a, b, c, \ldots; x_\lambda, x_{\lambda+1}, \ldots, x_{\lambda+m}\} = 1.$$

*Assume that no $R_\lambda$, which is not identically $1$, is identically equal to a word*

$$T_\lambda S_\lambda(a, b, c, \ldots; x_\lambda, \ldots, x_{\lambda+m-1}) T_\lambda^{-1}$$

*or*

$$T_\lambda S_\lambda(a, b, c, \ldots; x_{\lambda+1}, \ldots, x_{\lambda+m}) T_\lambda^{-1}$$

*i.e. where $S_\lambda$ does not contain $x_\lambda$ or $x_{\lambda+m}$. (On the other hand, $x_{\lambda+1}$ to $x_{\lambda+m-1}$ need not occur in $R_\lambda$, and neither does any of the $a, b, c, \ldots$). Here $T_\lambda$ is any word.*

*Then: all relations for the subgroup of $G$ generated by the $a, b, c, \ldots; x_0, x_1, \ldots, x_M$ (where $M$ is an integer $\geq 0$) are given by*

$$R_\mu\{a, b, c, \ldots; x_\mu, x_{\mu+1}, \ldots, x_{\mu+m}\} = 1$$

*for $0 \leq \mu \leq M - m$. If $M < m$, then $a, b, c, \ldots; x_0, x_1, \ldots, x_M$ generate a free group.*

We can now begin with the actual task at hand: we are tasked with finding a procedure, which in finitely many steps decides whether or not an arbitrary word

$W(a, b)$ on the generators $a, b$ of the group $G$ with the defining relation $R(a, b) \equiv a^{\alpha_1} b^{\beta_1} a^{\alpha_2} b^{\beta_2} = 1$ is equal to 1 on the basis of this relation.

The $\alpha_1, \beta_1, \alpha_2, \beta_2$ are integers $\geq 0$. There are three cases to investigate, depending on whether both of the numbers $\alpha_1 + \alpha_2 = \alpha$ and $\beta_1 + \beta_2 = \beta$ are equal to zero, if neither of them are, or if exactly one of them vanishes.

*First Case:* $\alpha_1 + \alpha_2 = 0;$ $\quad \beta_1 + \beta_2 = 0.$

In groups $G(a, b, R(a, b) = 1)$, in which $a$ and $b$ have exponent sum zero in $R$, one can always decide whether a word $W$ is a power of a generator if one can solve the identity problem – indeed, if $W = a^k$, then $k$ must be equal to the exponent sum of $a$ in $W$, and hence is uniquely determined by $W$ – so in this case we can apply Lemma 1 twice and reduce to the case of $G(a, b, aba^{-1}b^{-1} = 1)$, which is trivial.

*Second Case:* $\alpha_1 + \alpha_2 \neq 0;$ $\quad \beta_1 + \beta_2 \neq 0.$

We set $a = s \cdot t^{-\beta}, b = t^{\alpha}$ (with $\alpha = \alpha_1 + \alpha_2$ and $\beta = \beta_1 + \beta_2$), and

$$G(a, b; a^{\alpha_1} b^{\beta_1} a^{\alpha_2} b^{\beta_2} = 1)$$

is thus changed into a group

$$\overline{G}\{s, t; (st^{-\beta})^{\alpha_1} t^{\alpha \cdot \beta_1} (st^{-\beta})^{\alpha_2} t^{\alpha \beta_2} = 1\}.$$

By using Lemma 2[21] it is easy see that it suffices to solve the identity problem for $\overline{G}$ and – since the defining relation has exponent sum zero in $t$ – it even suffices to solve the identity problem in the subgroup $\overline{H}_t$ of $\overline{G}$ which consists of all words of $\overline{G}$ which have exponent sum zero in $t$. For a proof, see §2, second application of Theorem 2. From Theorem 4 of §2, it furthermore follows that $\overline{H}_t$ is defined using the $s_k$ ($s_k = t^k s t^{-k}; k = 0, \pm 1, \pm 2, \dots$) and the relations:

$$\overline{P}_\lambda \equiv s_\lambda s_{-\beta + \lambda} \cdots s_{-\beta(\alpha_1 - 1) + \lambda} \cdot s^\varepsilon_{\mu + \lambda} s^\varepsilon_{\mu + \lambda - \beta \cdot \varepsilon} \cdots s^\varepsilon_{\mu + \lambda - \beta(\alpha_2 - \varepsilon)} = 1,$$

where $\alpha_1 > 0$ is assumed, $\varepsilon$ is the sign of $\alpha_2$, and $\lambda$ ranges across all integers $0, \pm 1, \pm 2, \dots$, while $\mu = \alpha_2 \beta_1 - \alpha_1 \beta_2 + \beta$ in the case that $\varepsilon = -1$, and is otherwise equal to $\alpha_2 \beta_1 - \alpha_1 \beta_2$, in the case that $\varepsilon = +1$. It is now easy to see that – except for in the easily disposed-of cases of $\alpha_1 = \alpha_2, \alpha_2 \beta_1 - \alpha_1 \beta_2 = 0$, i.e. $\beta_1 = \beta_2$, that either the $s_k$ with the largest or the $s_k$ with the smallest $k$ appears only once in $\overline{P}_\lambda$. Indeed, the indices of the $s_k$ in $\overline{P}_\lambda$ from $s_\lambda$ to $s_{-\beta(\alpha_1-)+\lambda}$ and from $s^\varepsilon_{\mu+\lambda}$ to $s^\varepsilon_{\mu+\lambda-\beta(\alpha_2-\varepsilon)}$ are altered in a monotone way in $\beta$, and hence for the largest resp. smallest indices of the $s_k$ in $\overline{P}_\lambda$ only the numbers $\lambda; \mu + \lambda; -\beta(\alpha_1 - 1) + \lambda; \mu + \lambda - \beta(\alpha_2 - \varepsilon)$ can appear, and these must be pairwise equal, if not at least one only appears once. In the case that $\varepsilon = -1$, we must have, since the indices from $\lambda$ to $-\beta(\alpha_1 - 1) + \lambda$ and $\mu + \lambda$ to $\mu + \lambda - \beta(\alpha_2 - \varepsilon)$ change in the opposite directions, we have $\mu = -\beta(\alpha_1 - 1)$ and $\mu + \lambda - \beta(\alpha_2 - \varepsilon) = \lambda$; this easily gives $\alpha_1 + \alpha_2 = 0$, which was assumed not to be the case. Analogously, $\varepsilon = +1$ only gives the possibility $\mu = 0, \alpha_1 = \alpha_2$. This brings us to an exceptional case; in this case, $G$ is defined by

$$G\{a, b; (a^{\alpha_1} b^{\beta_1})^2 = 1.$$

By using Lemma 1 we can reduce this to the case of the identity problem in the group $G'\{a, b; (ab^{\beta_1})^2 = 1\}$, which is easy to solve, since we can introduce $ab^{\beta_1}$ as a new generator.

---

[21] and by introducing a new primitive element; we first set $b = t^\alpha$ and use, instead of $a$ and $t$, new generators $s = at^\beta$ and $t$.

Other than this exceptional case, we have proved that in each relation $\overline{P}_\lambda = 1$ the $s_k$ with the largest or the $s_k$ with the smallest index only appears exactly once. Since every word $\overline{H}_t$ only contains finitely many $s_k$, it is sufficient to solve the identity problem in *the* subgroup of $\overline{H}_t$ which is generated by finitely many $s_k$: $s_N, s_{N+1}, \ldots, s_M$ ($N \leq M$), in order to solve it in $\overline{H}_t$ itself. By the *Freiheitssatz* the subgroup of $\widetilde{H}_t$ generated by the $s_N$ to $s_M$ has as defining relations finitely many of the relations $\overline{P}_\lambda = 1$, and since one can by the above use these relations to eliminate the superfluous generators, each such subgroup must be a free group, in which one can solve the identity problem. Thus we are left with the:

*Third Case:* $\alpha_1 + \alpha_2 \neq 0$; $\quad \beta_1 + \beta_2 = 0$.

By Lemma 1 this case can be reduced to, in the case when $\beta_1 = 1, \beta_2 = -1$, and $\alpha_1$ is coprime with $\alpha_2$: $(\alpha_1, \alpha_2) = 1$, solving in the group

$$G\{a, b; a^{\alpha_1} b a^{\alpha_2} b^{-1} = 1; \quad (\alpha_1, \alpha_2) = 1\}$$

not only the identity problem, but also deciding when an element $W(a, b)$ of $G$ is equal to a power of $a$, and to deciding in

$$G'\{a, b; a^{d\alpha_1} b a^{d\alpha_2} b^{-1} = 1\}$$

($d$ an integer, $|d| > 1$) when an element is equal to a power of $b$, in the case that the identity problem is solved in $G'$. This second case is solved on the basis of Theorem 2, §2 and a conclusion already used in the "first case" $\alpha = \beta = 0$. We will solve the first case by using the following course of reasoning.

*Plan for the proof.* Suppose the group $G(a, b; a^{\alpha_1} b a^{\alpha_2} b^{-1} = 1)$ with $\alpha_1 + \alpha_2 \neq 0$, $(\alpha_1, \alpha_2) = 1$ is given. If one wishes to solve only the identity problem in $G$ (we should also be able to decide when a word of $G$ is a power of $a$ or not), then it suffices to be able to solve the identity problem in the subgroup $H_b$ of $G$, which consists of all words of $G$ which have exponent sum zero in $b$ (§2, Theorem 2, 2nd application). Now, $H_b$ appears as a group with infinitely many generators; however, since every word in $H_b$ only uses finitely many generators, it suffices to solve the identity problem in some subgroups $\overline{H}_m$ of $H_b$, in which the $\overline{H}_m$ has only $m + 1$ generators, which, as it will turn out, can all be expressed in terms of only two of them. The $\overline{H}_m$ will in this way appear as groups with only two generators $a_0$ and $a_m$. To solve the identity problem in $\overline{H}_m$, it is furthermore sufficient to solve it in the subgroup of $\overline{H}_m$ consisting of all words which have exponent sum zero in $a_m$. In this subgroup, it is easy to solve the identity problem through repeated usage of the theorem of Schreier.

Thus we will have solved the identity problem for $G$; however, the task to solve the problem of deciding when a word $W(a, b)$ in $G$ is equal to a power of $a$ remains. We will now show that, in fact, this problem can be solved when one can solve the identity problem in the groups $\overline{H}_m$. This will complete our proof.

We thus begin with solving the identity problem for the subgroup $H_b$ of

$$G(a, b; a^{\alpha_1} b a^{\alpha_2} b^{-1} = 1; \quad \alpha_1 + \alpha_2 \neq 0; (\alpha_1, \alpha_2) = 1)$$

which consists of all words of $G$ which have exponent sum zero in $b$. By §2, Theorem 4 we can present $H_b$ (when we set $b^k a b^{-k} = a_k$ for $k = 0, \pm 1, \pm 2, \ldots$) by

$$H_b\{a_k; a_k^{\alpha_1} = a_{k+1}^{-\alpha_2}; \quad k = 0, \pm 1, \pm 2, \ldots\}.$$

Since every word in $H_b$ only contains finitely any generators, it suffices to solve the identity problem for all subgroups of $H_b$ generated by finitely many of the $a_k$, and without loss of generality it thus suffices to solve it in the subgroup $\overline{H}_m$ generated by $a_0, a_1, \ldots, a_m$ ($m$ an integer $\geq 0$). By the *Freiheitssatz* this is defined by

$$H_m\{a_0, a_1, \ldots, a_m; a_\mu^{\alpha_1} = a_{\mu+1}^{-\alpha_2}; \quad \mu = 0, 1, \ldots, m-1\}.$$

Now we will show that since $(\alpha_1, \alpha_2) = 1$, the generators $a_1, a_2, \ldots, a_{m-1}$ can be expressed in terms of $a_0$ and $a_m$. From $a_0^{\alpha_1} = a_1^{-\alpha_2}, a_1^{\alpha_1} = a_2^{-\alpha_2}$ it follows that $a_0^{\alpha_1^2} = a_1^{-\alpha_2 \cdot \alpha_1} = a_1^{\alpha_1 \cdot (-\alpha_2)} = a_2^{(-\alpha_2)^2}$, and more generally

$$(3) \qquad a_0^{\alpha_1^{m-\mu}} = a_{m-\mu}^{(-\alpha_2)^{m-\mu}}; \quad a_{m-\mu}^{\alpha_1^\mu} = a_m^{(-\alpha_2)^\mu} \qquad (\mu = 0, 1, \ldots, m-1)$$

Since we have $(\alpha_1, \alpha_2) = 1$, we can for each $\mu = 1, 2, \ldots, m-1$ find a pair of integers[22] $\lambda_1, \lambda_2$ so that

$$(4) \qquad\qquad \lambda_1 \alpha_1^\mu + \lambda_2 (-\alpha_2)^{m-\mu} = 1$$

and consequently also

$$(5) \qquad\qquad a_{m-\mu} = a_0^{\lambda_2 \alpha_1^{m-\mu}} \cdot a_m^{\lambda_1 (-\alpha_2)^\mu}.$$

Thus $\overline{H}_m$ is generated by $a_0$ and $a_m$ alone. We will now determine a defining set of (essential) relations between $a_0$ and $a_m$. From (3) it follows for $\mu = 0$ that

$$(6) \qquad\qquad a_0^{\alpha_1^m} = a_m^{(-\alpha_2)^m},$$

and furthermore, if $\rightleftarrows$ denotes "commutes with", then

$$(7) \qquad\qquad a_0^{\alpha_1^{m-\mu}} \rightleftarrows a_m^{(-\alpha_2)^\mu} \qquad\qquad (\mu = 1, 2, \ldots, m-1)$$

The relations (6) and (7) form a complete system of defining relations for $\overline{H}_m$; indeed, since we have $\overline{H}_m$ as

$$\overline{H}_m(a_0, a_1, \ldots, a_m; a_\mu^{\alpha_1} = a_{\mu+1}^{\alpha_2}; \quad \mu = 0, 1, \ldots, m-1),$$

and if we express the $a_\mu$ using $a_0$ and $a_m$ by (5), then $a_\mu^{\alpha_1} = a_{\mu+1}^{\alpha_2}$, using (7) and as a simple consequence of (6), becomes of the form $a_0^{\tau a_1^m} = a_m^{\tau(-\alpha_2)^m}$, where $\tau$ is an integer.

Thus $\overline{H}_m$ is now defined by

$$\overline{H}_m\{a_0, a_m; a_0^{\alpha_1^m} = a_m^{(-\alpha_2)^m}; a_0^{\alpha_1^{m-\mu}} \rightleftarrows a_m^{(-\alpha_2)^\mu}; \quad \mu = 1, \ldots, m-1\}.$$

To now solve the identity problem in $\overline{H}_m$, we note that it is sufficient to solve it for the subgroup $\overline{H}_m^*$ of $\overline{H}_m$ consisting of all words in $\overline{H}_m$ which have exponent sum zero in $a_M$. Indeed, in every word $\overline{W}(a_0, a_m)$ from $\overline{H}_m$ which is equal to 1, the exponent sum of the letter $a_m$ must be divisible by $(-\alpha_2)^m$, say $M \cdot (-\alpha_2)^m$; hence the word

$$\overline{W}^* \equiv \overline{W} \cdot a_0^{M \cdot \alpha_1^m} a_m^{-M \cdot (-\alpha_2)^m}$$

lies in $\overline{H}_m^*$, and is also equal to 1.

To now solve the identity problem in $\overline{H}_m^*$, we consider the first application of §2, Theorem 2, which says that every word $\overline{W}^*(a_0 a_m)$ belonging to $\overline{H}_m^*$ is or is not equal to 1 only on the basis of the relations

$$(8) \qquad\qquad a_0^{-\alpha_1^m} a_m^{(-\alpha_2)^m} \rightleftarrows a_0; \quad a_0^{-\alpha_1^m} a_m^{(-\alpha_2)^m} \rightleftarrows a_m$$

---

[22] Of course dependent on $\mu$.

($\rightleftarrows$ means "commutes with") and

(9) $$a_0^{\alpha_1^{m-\mu}} \rightleftarrows a_m^{(-\alpha_2)^\mu}. \qquad (\mu = 1, 2, \ldots, m-1)$$

If we now set, for $k = 0, \pm 1, \pm 2, \ldots$

$$a_m^k a_0 a_m^{-k} = d_k,$$

then by §2, Theorem 4 (since we can write (8) and (9) over the $d_k$) $\overline{H}_m^*$ is given by

$$\overline{H}_m^*\{d_k; \; d_{k+(-\alpha)^\mu}^{\alpha_1^{m-\mu}} = d_k^{\alpha_1^{m-\mu}}; \quad k = 0, \pm 1, \ldots; \; \mu = 0, 1, \ldots, m\}.$$

In the remainder, we will make use of the following:

**Lemma 3.** *In groups of the following type, which we will denote by a $*$:*

$$\text{Generators:} \quad e_i \; (i = 1, 2, \ldots),$$
$$\text{Relations:} \quad e_i^{n_{ik}} = e_k^{n_{ik}}$$

*(where the $n_{ik}$ are arbitrary integers $\gtreqless 0$; and $k$ runs across the same indices as $i$), one can decide when a word is equal to a power of a generator when one can solve the identity problem in the group.*

*Proof.* If $W(e_{i_1}, e_{i_2}, \ldots) = e_k^l$, and if $\varepsilon_1, \varepsilon_2, \ldots$ is the exponent sums of $e_{i_1}, e_{i_2}, \ldots$ in $W$, then $l = \varepsilon_1 + \varepsilon_2 + \cdots$, and hence is uniquely determined by $W$. Indeed, if one considers the quotient group of our group obtained by adding the relations $e_1 = e_2 = e_3 = \cdots = e$, then $W(e_{i_1}, e_{i_2}, \ldots) = e_k^l$ becomes $e^{\varepsilon_1 + \varepsilon_2 + \cdots} = e^l$, and since no relation for $e$ alone follows, we must have $\varepsilon_1 + \varepsilon_2 + \cdots = l$. $\qquad\square$

By using Schreier's theorem, the following is now a consequence: if one can solve the identity problem in all groups of $*$-type, then we can always solve it in free products with amalgamated subgroups of such groups, when the amalgamated subgroups are such subgroups generated by a single generator or a power of a generator (these subgroups are hence all isomorphic with the free group on one generator).

We now return to $\overline{H}_m^*$. We can assume $-\alpha_2 > 0$, and by elimination of superfluous generators $\overline{H}_m^*$ is defined by

$$\overline{H}_m^*\{d_k; d_{k+(-\alpha_2)^\mu}^{\alpha_1^{m-\mu}}; \; k = 0, 1, \ldots, ((-\alpha_2)^m - 1); \mu = 0, \ldots, m-1\},$$

where the relations are naturally defined such that as indices of $d$ none other than $0, 1, \ldots ((-\alpha_2)^m - 1)$ appear.

To show that we can solve the identity problem in $\overline{H}_m^*$, it will be useful to begin with the following *simple example*: if $m = 2$, and $-\alpha_2 = +2$, then $\overline{H}_m^*$ is defined by the generators $d_0, d_1, d_2, d_3$ and the relations

$$d_0^{\alpha_1^2} = d_1^{\alpha_1^2} = d_2^{\alpha_1^2} = d_3^{\alpha_1^2};$$
$$d_0^{\alpha_1} = d_2^{\alpha_1}; \quad d_1^{\alpha_1} = d_3^{\alpha_1}.$$

We can see that our group is a free product with amalgamated subgroups of the following groups $\overline{H}_2^*$ and $\overline{H}_2^*$:

$$\overline{H}_2^*\{d_0, d_2; \; d_0^{\alpha_1} = d_2^{\alpha_1}\}$$

and

$$\overline{H}_2^*\{d_1, d_3; \; d_1^{\alpha_1} = d_3^{\alpha_1}\}$$

where the amalgamated subgroups are the subgroups generated by the elements $d_0^{\alpha_1^2}$ $(= d_2^{\alpha_2^2})$ resp. $d_1^{\alpha_1^2}$ $(= d_3^{\alpha_1^2})$. The group $\overline{H}^*$ is, for its own part, a free product with amalgamated subgroups of the groups generated by $d_0$ and $d_2$, where the amalgamated subgroups are the subgroups generated by $d_0^{\alpha_1}$ resp. $d_2^{\alpha_1}$. The analogous statement is true of $H_2^*$. By a second application of Lemma 3, we thus have a solution to the identity problem in this simple example.

We now solve the identity problem for $\overline{H}_m^*$ in general in an analogous way. Indeed, $\overline{H}_m^*$ is the free product with amalgamated subgroups of the following $|\alpha_2|$ groups $\overline{H}_{m,\sigma_1}^*$ ($\sigma_1 = 0, 1, \ldots, |\alpha_2| - 1$) of $^*$-type: if $\tau_1$ ranges across the numbers $0, 1, \ldots, ((-\alpha_2)^{m-1} - 1)$ and $\sigma_1$ the numbers $0, 1, \ldots, ((-\alpha_2) - 1)$, then $\overline{H}_{m,\sigma_1}^*$ is defined[23] by:

$$\overline{H}_{m,\sigma_1}^*\{d_{\sigma_1 + \tau_1 \cdot (-\alpha_2)};\ d_{\sigma_1 + (\tau_1 + (-\alpha_2)^\mu)(-\alpha_2)}^{\alpha_1^{m-\mu-1}} = d_{\sigma_1 + \tau_1 \cdot (-\alpha_2)}^{\alpha_1^{m-\mu-1}}\},$$

where $\mu = 0, 1, \ldots, m - 2$, and as amalgamated subgroups we take the subgroup of $\overline{H}_{m,\sigma_1}^*$ generated by

$$d_{\sigma_1}^{\alpha_1^m}\ (= d_{\sigma_1 + (-\alpha_2)}^{\alpha_1^m} = = d_{\sigma_1 + 2 \cdot (-\alpha_2)}^{\alpha_1^m} = \cdots = d_{\sigma_1 + ((-\alpha_2)^{m-1} - 1)(-\alpha_2)}^{\alpha_1^m}).$$

The groups $\overline{H}_{m,\sigma_1}^*$ are themselves free products of the following $-\alpha_2$ groups $\overline{H}_{m,\sigma_1,\sigma_2}^*$ ($\sigma_2 = 0, 1, \ldots, ((-\alpha_2) - 1)$) of $^*$-type with amalgamated subgroups: if $\tau_2$ ranges across the numbers $0, 1, \ldots, ((-\alpha_1)^{m-2} - 1)$, then $\overline{H}_{m,\sigma_1,\sigma_2}^*$ is given by

$$\overline{H}_{m,\sigma_1,\sigma_2}^*\{d_{\sigma_1 + \sigma_2 \cdot (\alpha_2) + \tau_2 \cdot (-\alpha_2)^2};\ d_{\sigma_1 + \sigma_2 \cdot (-\alpha_2) + (\tau_2 + (-\alpha_2)^\mu)(-\alpha_2)^2}^{\alpha_1^{m-\mu-2}} = d_{\sigma_1 + \sigma_2 \cdot (-\alpha_2) + \tau_2 \cdot (-\alpha_2)^2}^{\alpha_1^{m-\mu-2}}\},$$

where $\mu$ ranges across the numbers $0, 1, \ldots, m - 3$, and the amalgamated subgroups are the subgroups of $\overline{H}_{m,\sigma_1,\sigma_2}^*$ generated by[24]

$$d_{\sigma_1 + \sigma_2(-\alpha_2)}^{\alpha_1^{m-1}}\ (= d_{\sigma_1 + \sigma_2(-\alpha_2) + 1 \cdot (-\alpha_2)^2}^{\alpha_1^{m-1}} = \cdots = d_{\sigma_1 + \sigma_2(-\alpha_2) + ((-\alpha_2)^{m-2} - 1) \cdot (-\alpha_2)^2}^{\alpha_1^{m-1}}).$$

One can now easily see how this continues, and that after $m$ applications of Lemma 3, the identity problem in $\overline{H}_m^*$ reduces to that in the free group on one generator.

If $-\alpha_2 < 0$, then the above formulae hold when replacing $\alpha_2$ by $-\alpha_2$.

As indicated earlier (p. 43), we now have to provide a proof that we can decide in the original group $G(a, b;\ a^{\alpha_1} b a^{\alpha_2} b^{-1} = 1)$ whether a word $W(a, b)$ is equal to a power of $a$, given that we can solve the identity problem in $\overline{H}_m^*$. We thus ask:

What does it mean to decide when a word $W(a, b) = a^l$ in $G$ is equal to a power of $a$, for $H_b$? Since $W$ belongs by §2, Theorem 2 to $H_b$, there is some suitable power $b^\lambda$ of $b$ with $\lambda \geq 0$, so that $b^\lambda W b^{-\lambda}$ belongs to $\overline{H}_m^*$ (for a sufficiently large $m$). Such a $\lambda$ is easy to construct, and we choose $m$ large enough such that $0 \leq \lambda \leq m$. If we set $b^\lambda W(a, b) b^{-\lambda} = \overline{W}(a_0, a_m)$, then we must decide in $\overline{H}_m$ whether $\overline{W}(a_0, a_m) = a_\lambda^l$. If we set $\lambda = m - \mu$ ($0 \leq \mu \leq m$), then we must (by (5)) have that $\overline{W}(a_0, a_m) = a_0^{l \cdot \lambda_2 \alpha_1^{m-\mu}} a_m^{l \lambda_1 (-\alpha_2)^\mu}$ follows from $a_0^{\alpha_1^m} = a_m^{(-\alpha_2)^m}$ and (7); if $\gamma_0$ resp.

---

[23] In the case that $m = 1$, no relations appear here; and we are therefore done.
[24] In the case that $m = 2$, we are therefore done.

$\gamma_m$ are the exponent sums of $a_0$ resp. $a_m$ in $\overline{W}$, then we must by a generalization[25] of §2, Theorem 2, first application, have an integer $\nu$ such that

$$\nu \cdot \alpha_1^m = \gamma_0 - l \cdot \lambda_2 \alpha_1^{m-\mu}$$
$$-\nu(-\alpha_2)^m = \gamma_m - l \cdot \lambda_1 (-\alpha_2)^\mu.$$

By using $\lambda_1 \alpha_1^\mu + \lambda_2(-\alpha_2)^{m-\mu} = 1$, we can thus determine $\mu$ and $l$ uniquely; and hence $l$ is uniquely determined[26] by $W$ and hence it suffices to solve the identity problem for $\overline{H}_m$ in order to decide in $G$ whether or not a word is equal to a power of $a$.

# §5.
# The subgroups of the modular group

The modular group $G$ is defined by the two generators[27] and the relations

$$a^2 = b^3 = 1.$$

We will first show that the commutator subgroup $C$ of the modular group is a free group on two generators. The commutator subgroup $C$ of the modular group consists of all words $W(a, b)$, which have exponent sum zero in $a$ and $b$.[28] By §2, Theorem 2, first application, such a word $W$ in the modular group $G$ is equal to $1$ on the basis of the relations

(1) $$a^2 b a^{-2} b^{-1} = 1; \quad b^3 a b^{-3} a^{-1} = 1.$$

Stated otherwise, we have: the commutator subgroup $C$ of the modular group $G$ is (single-step) isomorphic to the commutator subgroup of the group $\overline{G}$, which is given by the generators $a, b$ and the defining relations (1). The commutator subgroup of $\overline{G}$ will hence also be denoted by $C$.

To find generators and defining relations of $C$, we first investigate the subgroup $H_a$ of $\overline{G}$, which consists of the words in $\overline{G}$ which have exponent sum zero in $a$. Indeed, $C$ is a subgroup of $H_a$.

By §2, Theorem 4, $H_a$ is generated by

$$b_k = a^k b a^{-k} \qquad\qquad (k = 0, \pm 1, \pm 2, \dots)$$

and all relations between the $b_k$ follow from the relations

(2) $$b_{k+2} b_k^{-1} = 1; \quad b_k^3 b_{k+1}^{-3} = 1. \qquad\qquad (k = 0, \pm 1, \dots)$$

All $b_k$ can thus be written in terms of $b_0$ and $b-1$, and between these we now have the relation

(3) $$b_1^3 b_0^{-3} = 1.$$

---

[25] In §2, Theorem 2, first application, we only had that $W = 1$ followed from a single relation $R = 1$; *here* the relation $a_0^{\alpha_1^m} = a_m^{(-\alpha_2)^m}$ appears together with the additional relations (7), which, however, vanish upon abelianization since all generators have exponent sum zero.

[26] Indeed, $l = \gamma_0 \left( \frac{-\alpha_2}{\alpha_1} \right)^{m-\mu} + \gamma_m \left( \frac{\alpha_1}{-\alpha_2} \right)^\mu$.

[27] The generators $a$ and $b$ correspond, respectively, to the linear substitutions $z' = -\frac{1}{z}$ resp. $z' = -\frac{1}{z+1}$ of a single variable $z$.

[28] Of course, there are also words in $C$ which do not have exponent sum zero in $a$ and $b$; however, these words can, by using the defining relations of $G$, be transformed into words which have exponent sum in $a$ and $b$, e.g. $(ab)^6 = (ab)^6 a^{-6} b^{-6}$ belongs to $C$.

Thus, if a word $\overline{W}(b_0, b_1)$ in $H_a$, when one uses $b_0 = b$, $b_1 = aba^{-1}$ to transform it into a word $\overline{W}(b, aba^{-1})$ in $\overline{G}$, also has exponent sum zero in $b$ (and not only in $a$), then we must obviously have that in $\overline{W}(b_0, b_1)$ the sum of the exponent sums of $b_0$ and $b_1$ must be equal to zero. Then $C$, considered as a subgroup of $H_a$, consists of the words from $H_a$ for which the sum of the exponent sums of $b_0$ and $b_1$ is zero.

To now define $C$ by means of generators and relations, we consider next the subgroup $F'$ of the *free* group $F(b_0, b_1)$ generated by $b_0$ and $b_1$, which consists of all words for which the sum of the exponent sums of $b_0$ and $b_1$ is equal to zero. Then it is clear that $F'$ is the free group on the generators

$$\beta_i = b_0^i b_1 b_0^{-1} b_0^{-i}. \qquad\qquad (i = 0, \pm 1, \pm 2, \dots)$$

Indeed, it is first of all clear that the $\beta_i$ are free, since if one sets

$$b_0 = x, \; b_1 = y \cdot x, \quad \text{then} \quad \beta_i = x^i y x^{-i}$$

and the $x^i y x^{-i}$ are by §2, Theorem 3 free. Furthermore, it is not hard to see that the $\beta_i$ actually generate $F'$, since a word

$$\overline{W} \equiv b_0^{\gamma_1} b_1^{\delta_1} b_0^{\gamma_2} b_1^{\delta_2} \cdots b_0^{\gamma_m} b_1^{\delta_m}$$

from $F'$, for which thus $\gamma_1 + \gamma_2 + \cdots + \gamma_m + \delta_1 + \delta_2 + \cdots + \delta_m = 0$, can also be written as

$$\overline{W} \equiv b_0^{\gamma_1} (b_1^{\delta_1} b_0^{-\delta_1}) b_0^{-\gamma_1} \cdot b_0^{\gamma_1 + \delta_1 + \gamma_2} \cdot (b_1^{\delta_2} b_0^{-\delta_2}) b_0^{-\gamma_1 - \delta_1 - \gamma_2} \cdots$$
$$\cdots b_0^{\gamma_1 + \delta_1 + \cdots + \gamma_m} (b_1^{\delta_m} b_0^{-\delta_m}) b_0^{-\gamma_1 - \delta_1 - \cdots - \gamma_m} \cdot b_0^{\gamma_1 + \delta_1 + \cdots + \gamma_m + \delta_m},$$

and since by assumption $b_0^{\gamma_1 + \delta_1 + \cdots + \gamma_m + \delta_m} = 1$, we now only have to prove that each word $b_0^\gamma (b_1^\delta b_0^{-\delta}) b_0^{-\gamma}$ can be written in terms of the $\beta_i$, and this is clearly possible, since $b_0^\gamma b_1^\delta b_0^{-\delta} b_0^{-\gamma} = \beta_\gamma \cdot \beta_{\gamma+1} \cdot \beta_{\gamma+2} \cdots \beta_{\gamma + \delta - 1}$.

We now return to the question of defining $C$ as a subgroup of $H_a(b_0, b_1; b_1^3 b_0^{-3} = 1)$ via generators and relations. As we just proved, $C$ is generated by the $\beta_i = b_0^i b_1 b_0^{-1} b_0^{-i}$, and the question is now which relations for the $\beta_i$ follow from $b_1^3 b_0^{-3} = 1$. To investigate this, we know first of all that for every word $\overline{W}(b_0, b_1)$ from $H_a$ which is equal to $1$ on the basis of $b_1^3 b_0^{-3} = 1$, can by §2, Theorem 1 gives rise to an identity:

$$\overline{W} \equiv \prod_{\lambda=1}^{k} T_\lambda (b_1^3 b_0^{-3})^{\varepsilon_\lambda} T_\lambda^{-1},$$

where $\varepsilon_\lambda = \pm 1$. Clearly, $\overline{W}$ belongs to $C$, since in $T_\lambda b_1^3 b_0^{-3} T_\lambda^{-1}$ the sum of the exponent sums of $b_0$ and $b_1$ is zero. If the sum of the exponent sum zero of $b_0$ and $b_1$ in $T_\lambda$ is equal to $t_\lambda$, then $T_\lambda \equiv T_\lambda' b_0^{t_\lambda}$, where $T_\lambda'$ belongs to $C$ and can hence be written as a product of the $\beta_i$. From this it follows: all relations between the $\beta_i$, which follows from $b_1^3 b_0^{-3} = 1$, are obtained from the relations $b_0^{t_\lambda} (b_1^3 b_0^{-3})^{\varepsilon_\lambda} b_0^{-t_\lambda} = 1$ written over the $\beta_i$, which gives the relations

$$(4) \qquad\qquad \beta_i \beta_{i+1} \beta_{i+2} = 1 \qquad\qquad (i = 0, \pm 1, \pm 2, \dots)$$

and because $\overline{W}$ can be written as

$$\prod_{\lambda=1}^{k} T_\lambda' (\beta_{t_\lambda} \beta_{t_\lambda+1} \beta_{t_\lambda+2})^{\varepsilon_\lambda} T_\lambda'^{-1}$$

it follows that $\overline{W}$, considered as a word over the $\beta_i$, is equal to 1 on the basis of (4).[29]

Since in the relations (4) each generator only appears once, one can thus express all generators in terms of two of them, say, $\beta_{-1}$ and $\beta_0$. In this way, the relations (4) are "resolved" (i.e. by identities in $\beta_{-1}$ and $\beta_0$) and so $\beta_{-1}$ *and* $\beta_0$ *are free generators of* $C$. As elements of the modular group $G$ they are

$$\beta_{-1} = b^{-1}aba^{-1}; \quad \beta_0 = aba^{-1}b^{-1}.$$

They correspond to the substitutions

$$z' = \frac{2z-1}{-z+1}, \quad \text{resp.} \quad z' = \frac{-z+1}{z-2}$$

of a single variable $z$ (in the case that $a$ and $b$ are represented by the substitutions $z' = \frac{-1}{z}$ and $z' = \frac{-1}{z+1}$).

Let now $\Gamma$ be an arbitrary subgroup of the modular group $G$. Denote the intersection of $\Gamma$ and $C$ by $\Delta$; then $\Delta$ is a normal subgroup of $\Gamma$, and, as a subgroup of the free group $C$, is itself free by a theorem of O. Schreier[30] (where we regard the identity element as being the free group on zero generators).

The quotient group $\Gamma/\Delta$ is one-step isomorphic to a subgroup of the quotient group $G/C$ of $C$ in $G$; but $G/C$ is defined by the generators $a, b$ and the relations

$$a^2 = b^3 = aba^{-1}b^{-1} = 1,$$

since one obtains the quotient group of a normal subgroup by adding, to the relations of the entire group, a new relation making each element of the normal subgroup 1.

Thus $G/C$ is a cyclic group of order 6, and $\Gamma/\Delta$ has order 1, 2, 3, or 6.

---

(Received 25 October 1930)

---

[29]The same reasoning just used can also be used to derive Theorem 4 of §2 from Theorems 1 and 3 of §2.

[30]O. Schreier, *The subgroups of free groups*, Abhandlungen aus dem Matematischen Seminar der Hamburgischen Universität **5**, p. 161, Leipzig 1927.

# The identity problem for groups with one defining relation

By
W. Magnus in Göttingen.

## Introduction.

Let a group be given by certain (finitely or infinitely many) generating elements $a_1, a_2, a_3, \ldots$ and certain "defining relations" holding between these generators:

$$R_k(a_1, a_2, a_3, \ldots) = 1. \qquad (k = 1, 2, 3, \ldots)$$

Each of the expressions (or "*words*", as we will call them) formable from the generators $a_1, a_2, a_3, \ldots$ and their inverses $a_1^{-1}, a_2^{-1}, a_3^{-1}, \ldots$ thus represents an element of the group; however, not in a unique manner: furthermore, each element can be represented by words in infinitely many different ways. The identity or word problem is then the task of finding a procedure which for two given arbitrary words $W_1$ and $W_2$ decides, in finitely many steps, whether they represent the same group element, or, which is equivalent, deciding whether a given word represents the identity or not.

The identity problem is first of all directly relevant to topology[1]; but second, it is also of importance in investigations of infinite groups; if one e.g. wants to determine generators and defining relations for a subgroup $H$ of a group $G$ given by generators and defining relations, then one must specify[2] a system of representatives for the cosets of $H$ in $G$ and a procedure by which to find, for every element of $G$, the corresponding representative. In the case that $H$ is a normal subgroup, then this implies the identity problem for $G/H$.

In the present article we will solve the identity problem in groups with only *one* defining relation. Naturally, we will also gain some insight into the structure of groups with one defining relation, such as a theorem, formulated below, about the "free" groups contained in them as subgroups. Free groups are those with no defining relations[3]; two words $W_1$ and $W_2$, which represent the same element of a free group, can be transformed into one another by applications of the rules $a_1 a_1^{-1} = 1, a_1^{-1} a_1 = 1$; $a_2 a_2^{-1} = 1, a_2^{-1} a_2 = 1, \ldots$; in this case, we say: $W_1$ is identical to $W_2$, and write $W_1 \equiv W_2$.

---

[1] See M. Dehn, *On infinite discontinuous groups*. Math. Annalen **71** (1912), p. 116.

[2] See K. Reidemeister, *Knots and groups*. Treatises of the mathematical seminar of the University of Hamburg **5** (1927), p. 7; and O. Schreier, *The subgroups of free groups*, ibid. p. 161.

[3] Free groups are usually defined as those which can be represented by some system of generators between which no relation holds. In the sequel, however, we will usually speak of a group being free when no relation holds between the specified generators. For example, $a, b, c$ generate a free group on two generators when the relation $(abc)^2 ab = 1$ holds between them; but there is still no general procedure to decide whether an arbitrary given group is isomorphic to a free group or not.

We have the following: among all pairwise identical words $W$, there is exactly one, $W_0$, which cannot be written with fewer letters[4]; we find this word by deleting the words $a_1 a_1^{-1}, a_1^{-1} a_1; a_2 a_2^{-1}, a_2^{-1} a_2, \ldots$ as many times as possible from $W$; the same applies if one writes all words $W$ "*cyclically*", i.e. that the first letter is considered to be adjacent to the last, and one performs the identical transformations accordingly. We can thus in a free group solve not only the identity problem, but also the conjugacy problem, i.e. whether or not for two arbitrary given words $W_1$ and $W_2$ one can find a third, $T$, such that

$$W_1 \equiv T W_2 T^{-1}.$$

Furthermore, one can also solve the "extended" word problem in free groups, which we will describe below.

These facts will be used frequently in the sequel, since firstly some tasks are reduced e.g. to the conjugacy problem in a free group, and secondly since groups with one defining relation are closely linked with free groups in general; we will see this in its structure, using the construction, due to Schreier, of the "free product with amalgamated subgroups" of two groups, as well as e.g. the theorem that groups with one defining relation "in general" will contain a free subgroup on two (and hence also on infinitely many) generators. The only exceptions are the following cases:

(1) The group has only a single generator,
(2) The group has two generators $a$ and $b$, and is isomorphic to a group with a single defining relation

$$a b a^n b^{-1} = 1 \qquad\qquad (n = 0, \pm 1, \pm 2, \ldots)$$

In these cases, the commutator subgroup is abelian.

The proof of this will essentially[5] be conducted using the following methods.

When solving the identity problem it turns out to be useful to immediately create and solve a more general problem, which in the sequel will be called the "*extended identity problem*". This is stated as follows:

Let $a_{i_1}, a_{i_2}, a_{i_3}, \ldots$ be an arbitrary, possibly empty, subset of the generators $a_i$ $(i = 1, 2, 3, \ldots)$ of the group $G$ under consideration. Let

$$W(\ldots, a_i, \ldots)$$

be an arbitrary word on the $a_i$. Then we seek a procedure which in a finite number of steps decides whether or not, as a consequence of the relation $R = 1$ in $G$, the word can be transformed into a word

$$\overline{W}(\ldots, a_{i_\lambda}, \ldots)$$

consisting only of the generators $a_{i_1}, a_{i_2}, \ldots$ (in the case that the set is empty, then of course $\overline{W} \equiv 1$). It should always be required (even if it is not explicitly stated) that, in the case that $W$ is equal to at least one word $\overline{W}$, one can actually specify at least

---

[4]If we have a word written as

$$a_{i_1}^{\varepsilon_1} a_{i_2}^{\varepsilon_2} a_{i_3}^{\varepsilon_3} \cdots a_{i_k}^{\varepsilon_k},$$

where $\varepsilon_1 = \pm 1, \varepsilon_2 = \pm 1$, and where the numbers $i_1, i_2, \ldots, i_k$ are any numbers in the sequence $1, 2, 3, \ldots$, then $k$ is called the number of *letters* in the word. On the other hand, the number of *distinct* $a_i$ appearing in the word is called the number of *generators* in the word.

[5]In addition to the theorems in the following paragraphs, we will also use some investigations into the "root problem". See Magnus, *On discontinuous groups with one defining relation*, §7, in Journal für die reine und angewandte Mathematik **163** (1930).

one such word $\overline{W}$. The extended identity problem is thus a question which entirely depends on the specific representation of the group.

The solution of the extended identity problem will be carried out in several steps. In §1, we will first show, using a theorem of O. Schreier[6] on the "existence of the free product with amalgamated subgroups" of two groups and the so-called *Freiheitssatz*[7], certain simplifications of the questions under consideration; including the solution of the extended identity problem in a special case. In §2 we will show that, under some constraining assumptions, we can reduce the extended identity problem in the original group $G$ to the solution of the same problem in a one-relator group $H_0$, in which the defining relation contains fewer letters than that of $G$. Essential use will here be made of the aforementioned theorem of Schreier and the *Freiheitssatz*. In §3, we will show that the constraining assumptions can be fulfilled in a group $\overline{G}$ into which the given group $G$ can be embedded.

## §1.
## Simplifications and useful tools.

Let $a_1, a_2, a_3, \ldots$ be the generators of our group $G$. The number of generators need not be finite, although it is of course always the case that only a finite number of them will appear in the defining relation $R = 1$ of $G$. We will then say that a given generator, say $a_1$, will "***actually***" appear in the relation $R = 1$ when $R$, written cyclically, cannot be transformed via identical transformations into a word which no longer contains $a_1$, or, stated otherwise, if there is no word $T$ such that $T R T^{-1}$ is identical with a word which does not contain $a_1$. Then the *Freiheitssatz*[7] says: *if some generator, say $a_1$, appears in the relation $R = 1$ of $G$, then there is no relation among the set of all generators excluding $a_1$; that is, this set generates a free group.*

This provides a first simplification of the extended identity problem: if one wants to decide, for a word $W$ over the $a_i$ $(i = 1, 2, 3, \ldots)$, whether or not it is equal to a word $\overline{W}$ over the generators $a_{i_1}, a_{i_2}, a_{i_3}, \ldots$, and if $a_1$ is not among the generators $a_{i_1}, a_{i_2}, a_{i_3}, \ldots$, then it is enough to find a procedure for deciding whether or not $W$ is equal to a word $W'$ over the generators $a_2, a_3, \ldots$, and if so transform $W$ into such a word $W'$; then $W'$ is equal to a word $\overline{W}$ if and only if $W' \equiv \overline{W}$, and whether or not this is the case is easy to decide[8].

A further simplification comes from a theorem of O. Schreier[9], which for our purposes we will state as follows:

Let a group $G$ be given with two systems of generators $a_\mu$ $(\mu = 1, 2, 3, \ldots)$ and $b_\nu$ $(\nu = 1, 2, 3, \ldots)$, between which there are relations

$$R_i(\ldots, a_\mu, \ldots) = 1 \quad \text{and} \quad S_k(\ldots, b_\nu, \ldots) = 1. \quad (i, k = 1, 2, \ldots)$$

The group generated by the $a_\mu$ resp. the $b_\nu$ with the defining relations $R_i = 1$ resp. $S_k = 1$ will be denote $\mathfrak{A}$ resp. $\mathfrak{B}$. The subsets $a_{\mu_1}, a_{\mu_2}, a_{\mu_3}, \ldots$ resp. $b_{\nu_1}, b_{\nu_2}, b_{\nu_3}, \ldots$ generate subgroups $\mathfrak{C}_1$ resp. $\mathfrak{C}_2$ of $\mathfrak{A}$ resp. $\mathfrak{B}$. Furthermore, we assume $\mathfrak{C}_1$ and $\mathfrak{C}_2$ are holoedrically isomorphic, and moreover in such a way that the assignment of the $a_{\mu_\lambda}$

---

to $b_{\mu_\lambda}$ ($\lambda = 1, 2, 3, \dots$) forms an isomorphic image of $\mathfrak{C}_1$ onto $\mathfrak{C}_2$. In addition to the relations $R_i = 1, S_k = 1$, our group will also have the defining relations

$$a_{\mu_\lambda} = b_{\nu_\lambda}. \qquad\qquad (\lambda = 1, 2, 3, \dots)$$

Then $G$ is called the *free product of $\mathfrak{A}$ and $\mathfrak{B}$ with amalgamated subgroups $\mathfrak{C}_1$ and $\mathfrak{C}_2$*; in the case that $\mathfrak{C}_1$ and $\mathfrak{C}_2$ only consist of the identity element (i.e. when the set of $a_{\mu_\lambda}$ is empty), then $G$ is simply called the free product of $\mathfrak{A}$ and $\mathfrak{B}$. We set $a\mu_\lambda = b_{\mu_\lambda} = c_\lambda$ and denote by $A_1, A_2, \dots$ resp. $B_1, B_2, \dots$ esp. $C, C_1, \dots$ arbitrary words over $a_\mu$ resp. $b_\nu$ resp. $c_\lambda$. Then we have:

Every word $W$ in $G$ can be written in the form

$$C; \quad A_1 B_1 A_2 B_2 \cdots A_K B_K$$

where $A_1$ and $B_K$ may be equal to 1, but among the other words $A$ and $B$ none is equal to a word $C$ – except for when $W$ is equal to $C$. Such a representation will be called reduced. If $W \neq C$, and

$$W = \overline{A}_1 \overline{B}_1 \overline{A}_2 \overline{B}_2 \cdots \overline{A}_L \overline{B}_L$$

is another reduced representation of $W$, then we have $K = L$ and

$$\overline{C}_1 A_1 = \overline{A}_1 C_1, \quad \overline{C}_2 A_2 = \overline{A}_2 C_2, \quad \dots \quad \overline{C}_L A_L = \overline{A}_L C_L$$

$$\overline{C}'_1 B_1 = \overline{B}_1 C'_1, \quad \overline{C}'_2 B_2 = \overline{B}_2 C'_2, \quad \dots \quad \overline{C}'_L B_L = \overline{B}_L C'_L$$

(where it should be noted that the $c_\lambda$ can be both considered as generators $a_\mu$ as well as $b_\nu$). – From this, we have the following:

**Lemma**. *If $W$ is any word in the group $G$, and if one can decide in $\mathfrak{A}$ and $\mathfrak{B}$ whether a word $A$ or $B$ can be written as a word over the $c_\lambda$, then one can decide whether or not $W$ can be transformed into a word $A$ (or a word $B$).*

*Proof.* Indeed, if we have

$$W \equiv A_1 B_1 A_2 B_2 \cdots A_M B_M$$

is a not necessarily reduced representation of $W$ over the generators $a_\mu$ and $b_\nu$, then we consider – in the case that $M > 1$ – whether any of the words $A$ or $B$ can be transformed into one over the $c_\lambda$. If this is possible for $A_m$, then – since we have $R_i = 1$ – it follows that $A_m = C_m$, and hence equal to a word over $b_\nu$, and we can set

$$B_{m-1} A_m B_m = B'_m$$

where $B'_m$ is a new word over the $b_\nu$. We now have a new representation of $W$ as:

$$W = A_1 B_1 A_2 B_2 \cdots A_{m-1} B'_m A_{m+1} \cdots A_M B_M$$

and can perform the same reduction procedure. After at most $M$ steps, we obtain a reduced representation for $W$, and because of the aforementioned result it is now easy to decide whether it is equal to a word $A$ or $B$. $\qquad\square$

From this lemma, we will first make the following application for our original group $G$ with the generators $a_1, a_2, a_3, \dots$ and the defining relation $R = 1$: if the generators $a_1, a_2, a_3, \dots$ actually appear in $R$, and no other, then $G$ is a free product of the free group generated by $a_{n+1}, a_{n+2}, \dots$ by the group $G'$ with the generators $a_1, a_2, \dots, a_n$ and the defining relation $R = 1$. Thus, if we can solve the extended identity problem for $G'$, then we can solve it for $G$. Hence, if only a single generator actually appears in $R$, then $G'$ is a cyclic group, and we can thus clearly solve the extended identity problem for $G$.

## §2.
## Solution to the extended word problem in a special case.

In the previous section, we showed that in order to solve the extended word problem in arbitrary groups with one defining relation, it suffices to do so in the case for all groups $G$ with finitely many generators $a_1, \ldots, a_n$ and one defining relation

$$R(a_1, a_2, \ldots, a_n) = 1$$

where each of the generators $a_1, a_2, \ldots, a_n$ actually appears in $R$ written cyclically.

We will now show that the extended identity problem can be solved for the group $G$ when we can solve it for all groups $G^*$ with one defining relation, when certain assumptions hold for $G$ and when the defining relation of $G^*$ contains fewer letters than $R$; if we could do this, then the extended identity problem for all groups with one defining relation will be solved, since it can always be solved for groups whose defining relation only contains a single generator.

First, however, it must be said that the following procedure, i.e. to reduce the solution of the extended identity problem in $G$ to the solution of that in $G^*$, can only be carried out with two restrictions: the relation $R = 1$ of $G$ must be such that in $R$ at least two generators actually appear, and that at least one of those generators has "exponent sum"[10] zero in $R$. (The second assumption presupposes the first). The case when $R$ has only a single generator, is the case covered in §1. However, the case when in $R$ there are more (at least two) generators appearing, all of which have non-zero exponent sum, requires a special treatment. We will cover this later (in §3), and show that this case can be reduced to the one treated in this section.

Finally, it should be noted that all considerations can essentially be reduced to the case when the group $G$ has at most three generators; we will denote $a_1$ by $a$, $a_2$ by $b$, and $a_3$ (which does not necessarily appear) by $c$, and any remaining generators subsumed by an ellipsis $\ldots$, as this will make our notation more manageable.

Let then $a, b, c, \ldots$ be the generators, and $R(a, b, c, \ldots) = 1$ the defining relation of $G$, and suppose that $a$ has exponent sum zero in $R$. All elements in $G$ which have exponent sum zero in $a$ form a (normal) subgroup $H$ of $G$ which we can, as shown elsewhere[11], can present in the following way by generators and relations:

We set

$$a^k b a^{-k} = b_k, \ a^k c a^{-k} = c_k, \ \ldots$$

for all integers $k = 0, \pm 1, \pm 2, \ldots$. Then the $b_k, c_k, \ldots$ generate the group $H$, and the defining relations of $H$, which are written over the $b_k, c_k, \ldots$, are obtained when one expresses $R$ in terms of the $b_k, c_k, \ldots$, which is only possible because of the assumption that $a$ in $R$ has exponent sum zero.[12] Thus $R$ is rewritten into a word $\overline{R}$ over $b_k, c_k, \ldots$; the generators $b_k$ appearing in $\overline{R}$ are denoted $b_\mu$, and the $c_k$ appearing

---

[10]If we have

$$W \equiv a_1^{a_{1,1}} a_2^{a_{1,2}} \cdots a_n^{a_{1,n}} a_1^{a_{2,1}} a_2^{a_{2,2}} \cdots a_n^{a_{2,n}} \cdots a_1^{a_{k,1}} a_2^{a_{k,2}} \cdots a_n^{a_{k,n}},$$

then the sum

$$a_\nu = \sum_{\lambda=1}^{k} a_{\nu,\lambda}$$

is called the *exponent sum* of $a_\nu$ in $W$. The number $a_\nu$ is invariant under identical transformations of $W$.

[11]See footnote 5.

[12]*Proof.* Let $R$ have the form

$$R \equiv a^{\alpha_1} b^{\beta_1} c^{\gamma_1} \cdots a^{\alpha_2} b^{\beta_2} c^{\gamma_2} \cdots a^{\alpha_k} b^{\beta_k} c^{\gamma_k} \cdots,$$

in $\overline{R}$ are denoted $c_\nu$, where hence $\mu$ and $\nu$ are taken from some finite subset of integers in the sequence $0, \pm 1, \pm 2, \ldots$. This gives: if

$$R(a, b, c, \ldots) = \overline{R}(b_\mu, c_\nu, \ldots),$$

then all relations between the $b_k, c_k, \ldots$ follow from the relations

$$\overline{R}_\lambda \equiv \overline{R}(b_{\mu+\lambda}, c_{\nu+\lambda}, \ldots) = 1 \qquad\qquad (\lambda = 0, \pm 1, \ldots)$$

where hence $\lambda$ is for each relation some fixed number, and the $\mu$ and $\nu$ in $\overline{R}_\lambda$ are taken from some finite set of values. Thus we have a representation of $H$ in terms of infinitely many generators and defining relations, but where these generators are distributed across the generators in a particularly simple way, and the relations – in an immediately understandable way – are "isomorphic". Furthermore, for the sequel it will be important to note that *the $\overline{R}_\lambda$ contain fewer letters than $R$, and they contain at least as many letters as there are $a$-letters in $R$.* This is clear e.g. from the process of rewriting $R$ into $\overline{R}$ carried out in footnote 12.

We will now investigate what knowledge about $H$ will be required to solve the extended identity problem in $G$. As mentioned in §1, we must now for an arbitrary word $W(a, b, c, \ldots)$ from $G$ only decide whether it can be transformed into a word $W'$ which omits some fixed of the generators $a, b, c, \ldots$, and furthermore actually produce such a word $W'$. Since in the previous part the generator $a$ was distinguished from the other generators, our investigation is now divided into two parts:

First: when a word $W(a, b, c, \ldots)$ can be transformed into a word $W'$ which does not contain $a$, it follows that in $W$ the generator $a$ has exponent sum zero[13]. Hence, $W$ belongs to $H$, and can thus be written over the $b_k, c_k, \ldots$. Thus $W'$ also belongs to $H$ and, since it does not contain $a$, it can be written over $b_0, c_0, \ldots$. From this, we conclude easily:

*To decide whether a word $W(a, b, c, \ldots)$ can be transformed into a word $W'$ in which $a$ no longer appears, it is necessary and sufficient to be able to decide, for an arbitrary word*

$$\overline{W}(b_k, c_k, \ldots)$$

*from $H$, whether or not it can be transformed into a word $\overline{W}'$ from $H$ which only contains the generators $b_0, c_0, \ldots$.*

Second: if a word $W(a, b, c, \ldots)$ can be transformed into a word $W'$ which omits some generator distinct from $a$, say $b$, then $W a^{-\alpha}$ can also be transformed into such a word $W'$, where $\alpha$ is the exponent sum of $a$ in $W$. Since $a$ has exponent sum zero in $W a^{-\alpha}$ and hence belongs to $H$, we can write it over the $b_k, c_k, \ldots$. A word $W'$ which is equal to $W a^{-\alpha}$ and no longer contains $b$, must hence also have exponent sum zero[13]. From this, it follows that it is sufficient *to find a procedure to decide whether an arbitrary word*

$$\overline{W}(b_k, c_k, \ldots)$$

---

where some of the exponents can also be zero. By assumption, $\alpha_1 + \alpha_2 + \cdots + \alpha_k = 0$. On the other hand,

$$
\begin{aligned}
R \equiv &(a^{\alpha_1} b a^{-\alpha_1})^{\beta_1} (a^{\alpha_1} c a^{-\alpha_1})^{\gamma_1} \cdots (a^{\alpha_1 + \alpha_2} b a^{-\alpha_1 - \alpha_2})^{\beta_2} \\
&\times (a^{\alpha_1 + \alpha_2} c a^{-\alpha_1 - \alpha_2})^{\gamma_2} \cdots (a^{\alpha_1 + \alpha_2 + \cdots + \alpha_k} b a^{-\alpha_1 - \alpha_2 - \cdots - \alpha_k})^{\beta_k} \\
&\times (a^{\alpha_1 + \alpha_2 + \cdots + \alpha_k} c a^{-\alpha_1 - \alpha_2 - \cdots - \alpha_k})^{\gamma_k} \cdots a^{\alpha_1 + \alpha_2 + \cdots + \alpha_k},
\end{aligned}
$$

and since $a^{\alpha_1 + \alpha_2 + \cdots + \alpha_k} = 1$, our claim is proved.

[13]The $W'$ will arise from $W$ by inserting and deleting $R$ and $R^{-1}$ and identical transformations, and since $a$ has exponent sum zero in $R$, this does not change the exponent sum of $a$.

*from $H$ can be transformed into a word $\overline{W}'(c_k, \dots)$ which no longer contains $b_k$. Furthermore, one must also be able to actually write down such a word $\overline{W}'$.* Thus, we have reduced the extended identity problem for $G$ to the solution of two problems in $H$; to treat these, it is necessary to carry out a thorough investigation of $H$.

In the defining relations of $H$:

$$\overline{R}_\lambda \equiv \overline{R}(b_{\mu+\lambda}, c_{\nu+\lambda}, \dots) = 1$$

we denote the largest of the indices $\mu$ of some therein appearing $b_\mu$ by $M_1$, and the smallest by $M_0$. So $b_{M_0+\lambda}$ *and* $b_{M_1+\lambda}$ *actually appear in the cyclically written word* $\overline{R}_\lambda$.

Let now $\mu_0$ and $\mu_1$ denote two integers, and let $\mu_0 \leq \mu_1$. We denote by $\Gamma_{\mu_0,\mu_1}$ the subgroup of $H$ generated by the

$$b_{\mu_0}, b_{\mu_0+1}, \dots, b_{\mu_1}$$

and all $c_k, \dots$ ($k = 0, \pm 1, \pm 2, \dots$). Those groups $\Gamma_{\mu_0,\mu_1}$ for which $\mu_1 - \mu_0 = M_1 - M_0$ will be denoted by $H_\lambda$, where $\lambda = \mu_1 - M_1 = \mu_0 - M_0$. The following schema serves to illustrate this notation:

| Group | Generators | | | | | | | | |
|---|---|---|---|---|---|---|---|---|---|
| $\Gamma_{M_0, M_1+2}$ $\begin{cases} H_0 \\ H_1 \\ H_2 \end{cases}$ | $c_k, \dots;$ $c_k, \dots;$ $c_k, \dots;$ | $b_{M_0},$ | $b_{M_0+1}$ $b_{M_0+1}$ $b_{M_0+2}$ | $\dots$ $b_{M_0+2}$ $b_{M_0+3}$ | $\dots$ $\dots$ $\dots$ | $\dots$ $\dots$ $\dots$ | $\dots$ $\dots$ $\dots$ | $b_{M_1}$ $\dots$ $\dots$ | $b_{M_1+1}$ $b_{M_1+2}$ |

First, we collect some facts $\Gamma_{\mu_0,\mu_1}$ and $H_\lambda$, which we will use later.

I. By a theorem proved elsewhere[14], the *Freiheitssatz* in a suitable formulation, we have:

*The defining relations of the groups $H_\lambda$ consist of a single relation $\overline{R}_\lambda = 1$; the defining relations of $\Gamma_{\mu_0,\mu_1}$ for $\mu_1 - \mu_0 \geq M_1 - M_0$ are:*

$$\overline{R}_{\mu_0-M_0} = 1, \quad \overline{R}_{\mu_0-M_0+1} = 1, \quad \dots, \quad \overline{R}_{\mu_1-M_1} = 1,$$

*and when $\mu_0 - M_0 > \mu_1 - M_1$, the set of defining relations of $\Gamma_{\mu_0,\mu_1}$ is empty.*

II. From the inductive hypothesis at the beginning of this section we can solve the extended identity problem in $H_\lambda$ (and in particular for their subgroups $\Gamma_{\mu_0,\mu_1}$ with $\mu_1 - \mu_0 < M_1 - M_0$), since the defining relation of $H_\lambda$ contains fewer letters than $R$. The fact that the group $H_\lambda$ in general contains generators which do not appear in its defining relation when cyclically written (e.g. some $c_k$) is, as proved in §1, not an obstacle to applying this claim.

III. The groups $H_\lambda$ are, in a sense, the building blocks from which all the groups $\Gamma_{\mu_0,\mu_1}$ are constructed via Schreier's free product with amalgamated subgroups. Indeed, if $\mu_1 - \mu_0 > M_1 - M_0$, then $\Gamma_{\mu_0,\mu_1}$ arises from $\Gamma_{\mu_0,\mu_1-1}$ and $H_{\mu_1-M}$ by forming the free product of these to groups with amalgamated subgroup being that generated by all $c_k, \dots$ and

$$b_{\mu_1-M+M_0}, b_{\mu_1-M+M_0+1}, \dots, b_{\mu_1-1}$$

which generate the subgroup $\Gamma_{\mu_1-M_1+M_0}$ of $\Gamma_{\mu_0,\mu_1-1}$ and $H_{\mu_1-M_1}$. It is important to note that, by I., between the generators in $\Gamma_{\mu_1-M_1+M_0,\mu_1-1}$ there is no relation holding neither in $\Gamma_{\mu_0,\mu_1-1}$ nor in $H_{\mu_1-M_1}$; these generators thus generate isomorphic (free) subgroups of $\Gamma_{\mu_0,\mu_1-1}$ and $H_{\mu_1-M_1}$.

---

[14] See the citation in footnote 5, p. 157.

Let now $\overline{W}$ be an arbitrary word on the generators $b_k, c_k, \ldots$ from $H$. Then we should – in order to solve the extended identity problem in $G$ – decide for $\overline{W}$ whether or not it can be transformed into a word over the $c_k, \ldots$ alone or a word over $b_0, c_0, \ldots$. We can decide both, when we can decide whether $\overline{W}$ can be transformed into a word in the group $H_{-M_0}$ generated by

$$b_0, b_1, \ldots, b_{M_1 - M_0}; c_k, \ldots \qquad\qquad (k = 0, \pm 1, \pm 2, \ldots)$$

since, as noted in II., we can solve the extended identity problem in $H_{-M_0}$, and $H_{-M_0}$ contains both all $c_k, \ldots$ and $b_0$.

There is now, for any given $\overline{W}$, always a subgroup $\Gamma_{\mu_0, \mu_1}$ which contains both all the generators appearing in $\overline{W}$ as well as those of $H_{-M_0}$. Thus, we will be done if we can prove that for an arbitrary word $\overline{W}$, from an arbitrary group $\Gamma_{\mu_0, \mu_1}$, we can decide whether it can be transformed into an arbitrary of the subgroups $H_\lambda$ of the given $H_{\mu_0, \mu_1}$.

This last problem can be solved, when $\Gamma_{\mu_0, \mu_1}$ is one of the groups $H_\lambda$, i.e. when $\mu_1 - \mu_0 = M_1 - M_0$. By using complete induction, we can assume that the question is solved for $\Gamma_{\mu_0, \mu_1 - 1}$; to then solve it for $\Gamma_{\mu_0, \mu_1}$, it suffices to decide for an arbitrary word $\overline{W}$ from $\Gamma_{\mu_0, \mu_1}$ whether or not it can be transformed into a word in $\Gamma_{\mu_0, \mu_1 - 1}$ or one in $H_{\mu_1 - M_1}$; indeed, $\Gamma_{\mu_0, \mu_1}$ is the free product of these two groups, with amalgamated subgroup generated by all $c_k, \ldots$ and

$$b_{\mu_1 + M_0 - M_1}, b_{\mu_1 + M_0 - M_1 + 1}, \ldots, b_{\mu_1 - 1}.$$

We will denote this amalgamated subgroup of $\Gamma_{\mu_0, \mu_1 - 1}$ and $H_{\mu_1 - M_1}$ by $\Delta$. By the Lemma proved in §1, this question can now be easily solved when we can decide in $\Gamma_{\mu_0, \mu_1 - 1}$ and $H_{\lambda_1 - M_1}$ whether an arbitrary word in these groups can be transformed into a word over the generators of $\Delta$. For $H_{\mu_1 - M_1}$ this is clearly possible, since we can solve the extended identity problem in this group. For $\Gamma_{\mu_0, \mu_1 - 1}$ we proceed as follows: the generators of $\Delta$ form a subset of the generators of $H_{\mu_1 - M_1 - 1}$; in $H_{\mu_1 - M_1 - 1}$, we can solve the extended identity problem, and since $H_{\mu_1 - M_1 - 1}$ is a subgroup $H_\lambda$ contained in $\Gamma_{\mu_0, \mu_1 - 1}$, we can by the inductive assumption decide for any word in $\Gamma_{\mu_0, \mu_1 - 1}$ decide whether or not it is an element of $H_{\mu_1 - M_1 - 1}$, and hence also whether it can be transformed into an element of $\Delta$ or not.

## §3.
### Solution to the extended word problem in the general case.

It now only remains to free ourselves of the assumption, made in the course of the investigations in §2, that the defining relation $R = 1$ of the group $G$, for which we wish to solve the extended identity problem, when cyclically written contains at least two generators, one of which has exponent sum zero in $R$.

The case when only one generator actually appears in $R$ is covered in §1. We hence deal with the case when several generators $a, b, c, \ldots$ appear in $R$ when written cyclically, but that none of these generators have exponent sum zero in $R$.

Let $a$ and $b$ have exponent sum $s_1$ resp. $s_2$ in $R$. We will then replace $a$ and $b$ by new generators $\overline{a}$ and $\overline{b}$ by setting

$$a = \overline{a}^{s_2}; \quad b = \overline{b}\,\overline{a}^{-s_1};$$

Then $\overline{a}, \overline{b}, c, \ldots$ generate a group $\overline{G}$ with the single defining relation

$$R(\overline{a}^{s_2}, \overline{b}\,\overline{a}^{-s_1}, c, \ldots) \equiv P(\overline{a}, \overline{b}, c, \ldots) = 1,$$

and now $\overline{a}$ has exponent sum zero in $P$. If we now further assume that we can solve the extended identity problem for all one-relator groups whose defining relation has fewer letters than $R(a, b, c, \dots)$, then we can solve the same problem for $\overline{G}$ by using the results of the previous section. This is because if we look at the subgroup $\overline{H}$ of $\overline{G}$, which consists of all elements which have exponent sum zero in $\overline{a}$, then we find that $\overline{H}$ is constructed from certain subgroups $\overline{H}_\lambda$ with one defining relation in the same manner as in §2, where the defining relation of each group $\overline{H}_\lambda$ contains fewer letters than $R(a, b, c, \dots)$. This is because if we set

$$\overline{a}^k \overline{b} \overline{a}^{-k} = \overline{b}_k, \quad \overline{a}^k c \overline{a}^{-k} = \overline{c}_k, \dots \qquad (k = 0, \pm 1, \dots)$$

as the generators of $\overline{H}$, then the defining relations of $\overline{H}$ are

$$\overline{P}(\overline{b}_{\mu + \lambda}, \overline{c}_{\nu + \lambda}, \dots) \equiv \overline{P}_\lambda = 1 \qquad (\lambda = 0, \pm 1, \dots)$$

where $\overline{P}(\overline{b}_\mu, \overline{c}_\nu, \dots)$ is obtained from $P(\overline{a}, \overline{b}, c, \dots)$ by expressing it as a product of the $\overline{b}_k, \overline{c}_k, \dots$. Now each $\overline{P}_\lambda$ has fewer letters than $R(a, b, c, \dots)$, and indeed as many fewer letters as there are $a$-letters in $R$. Hence the number of $\overline{b}$- and $c$-letters in $P(\overline{a}, \overline{b}, c, \dots)$ corresponds to the number of $b$- and $c$-letters in $R(a, b, c, \dots)$.

By §2 we can now solve the extended identity problem for the group $\overline{G}$ with the generators $\overline{a}, \overline{b}, c, \dots$. This does not immediately tell us that we can solve the same problem in the original group $G$ with generators $a, b, c, \dots$. Indeed, when $|s_2| > 1$, we cannot express $\overline{a}$ and $\overline{b}$ in terms of $a$ and $b$. We have, however, the following result, which is proved elsewhere[15]:

The elements $\overline{a}^{s_2}, \overline{b} \overline{a}^{-s_1}, c, \dots$ in $\overline{G}$ generate a subgroup of $\overline{G}$ isomorphic to $G$ with the single defining relation

$$R(\overline{a}^{s_2}, \overline{b} \overline{a}^{-s_1}, c, \dots) = 1.$$

Furthermore, we have: when a word $\overline{W}(\overline{a}, c, \dots)$ from $\overline{G}$, which does not contain the generator $\overline{b}$, is equal to a word $W(a, c, \dots)$ from $G$, then we have

$$\overline{W}(\overline{a}, c, \dots) \equiv W(\overline{a}^{s_2}, c, \dots)$$

identically in $\overline{a}, c, \dots$. This follows simply from the fact that the group generated by $\overline{a}^{s_2}, c, \dots$ is a free subgroup of the free group generated by $\overline{a}, c, \dots$; indeed, by the *Freiheitssatz* we have that $\overline{a}^{s_2}, c, \dots$ and $\overline{a}, c, \dots$ generate *free* subgroups of $G$ resp. $\overline{G}$.

Thus, if we are given an arbitrary word $W(a, b, c, \dots)$ from $G$, then we can decide whether it can be transformed into a word $W'(a, c, \dots)$ which no longer contains $b$. Indeed, we can in $\overline{G}$ decide for the word

$$W(\overline{a}^{s_2}, \overline{b} \overline{a}^{-s_1}, c, \dots)$$

whether it can be transformed into a word $W(\overline{a}^{s_2}, c, \dots)$ or not.

Thus, we have solved the extended identity problem for $G$; since $b$ is in no way distinguished from the generators $a, c, \dots$ of $G$, we can equally well decide if a word from $G$ can be transformed into one which does not contain $a$ or $c$.

---



%\vspace{0.5cm}

\end{document}